\documentclass[11pt,a4paper]{article} %article文档

\usepackage[ruled,linesnumbered]{algorithm2e} %排版算法
\usepackage{amscd} %简单交换图
\usepackage{amsfonts} %字体
\usepackage{amsmath} %数学公式
\usepackage{amssymb} %更多公式符号
\usepackage{amsthm} %定理环境
\usepackage{bm} %允许公式粗体
\usepackage{booktabs} %改善表格
\usepackage{delarray} %矩阵
\usepackage{extarrows} %箭头
\usepackage{geometry} %页面布局
\usepackage{graphicx} %插图
\usepackage{hyperref} %超链接
\usepackage{longtable} %长表格
\usepackage{tikz-cd} %交换图
\tikzcdset{nfold/.style={double equal sign distance}}
\usepackage[all]{xy} %交换图
\usepackage[affil-it]{authblk} %标题

\usepackage{adjustbox}
\usepackage{booktabs} %美化表格（可选）
\usepackage{caption}
\usepackage{color}
\usepackage{diagbox} %用来画斜线的宏包
\usepackage{enumitem} %编号
\usepackage{float}
\usepackage{hyperref}%创建超文本链接和PDF书签
\usepackage{indentfirst} %首行缩进
\usepackage{makecell} %表格单元格换行
\usepackage{mathrsfs} %数学字体宏包
\usepackage{mathtools} %数学工具宏包
\usepackage{multirow} %跨行表格宏包
\usepackage{relsize} %缩放宏包
\usepackage{subfigure}

\setlist[enumerate]{itemsep=0pt,parsep=0pt}

\title{On the self-extensions of the mod $p$ Jacquet-Langlands correspondence over $\Qp$%
}
\author{Chenhao Tang%
\thanks{Email:
\href{mailto:tangchenhao25@mails.ucas.ac.cn}
{\nolinkurl{tangchenhao25@mails.ucas.ac.cn}}}}
\date{}

\allowdisplaybreaks[4] %多行公式换页

\def\CC{\mathbb{C}}

\def\FF{\mathbb{F}}

\def\PP{\mathbb{P}}

\def\ZZ{\mathbb{Z}}

\def\cC{\mathcal{C}}

\def\cE{\mathcal{E}}
\def\cF{\mathcal{F}}
\def\cG{\mathcal{G}}
\def\cH{\mathcal{H}}

\def\cL{\mathcal{L}}

\def\cO{\mathcal{O}}

\def\cS{\mathcal{S}}

\def\ad{\operatorname{ad}}

\def\Art{\operatorname{Art}}

\def\Coker{\operatorname{Coker}}

\def\d{\operatorname{d}}
\def\deg{\operatorname{deg}} 
 
\def\det{\operatorname{det}}

\def\dim{\operatorname{dim}}
\def\E{\operatorname{E}}
\def\End{\operatorname{End}}
\def\EP{\operatorname{EP}}
\def\Ext{\operatorname{Ext}} %%% new

\def\Gal{\operatorname{Gal}}

\def\GL{\operatorname{GL}}
\def\gr{\operatorname{gr}}

\def\Hom{\operatorname{Hom}}
\def\id{\operatorname{id}} 
\def\Im{\operatorname{Im}} 

\def\Ind{\operatorname{Ind}}

\def\max{\operatorname{max}}

\def\min{\operatorname{min}}

\def\ord{\operatorname{ord}}

\def\ss{\operatorname{ss}}

\def\Tor{\operatorname{Tor}}

\def\wet{\operatorname{wt}}

\def\Fp{\mathbb{F}_p}

\def\Fq{\mathbb{F}_q}

\def\Qp{\mathbb{Q}_p}
\def\Qpbar{\overline{\mathbb{Q}}_p}

\def\rhobar{\overline{\rho}}

\begin{document}

\newtheorem{definition}{Definition}[section] %%编号以section为单位
\newtheorem{remark}[definition]{Remark}
\newtheorem{example}[definition]{Example}
\newtheorem{proposition}[definition]{Proposition}
\newtheorem{lemma}[definition]{Lemma}
\newtheorem{corollary}[definition]{Corollary}
\newtheorem{theorem}[definition]{Theorem}
\newtheorem{conjecture}[definition]{Conjecture}

\maketitle

\begin{abstract}
   Let $p\ge5$ and $D$ be the unique nonsplit quaternion algebra over $\Qp$. We compute, in the category of smooth representations with a fixed central character, the self-extension groups of the $D^{\times}$-representations occurring in the mod $p$ Jacquet-Langlands correspondence in the generic case. The main ingredient is a new weighted filtration on the relevant ring and module. The associated graded module preserves the first nontrivial extension relation, which is not detected by the $\mathfrak m_D$-adic filtration.
\end{abstract}

\tableofcontents

\section{Introduction}

Let $p\ge5$ be a prime. The mod $p$ local Langlands correspondence for $\GL_2(\Qp)$ is well-understood by the work of \cite{Bre03,Col10,Eme11,Pas13}. The mod $p$ Jacquet-Langlands correspondence is partially understood by the work of \cite{Sch18,Lud17,Pas22,HW24a,HW24b}. In this paper, we further study the self-extensions of the representations occurring in the mod $p$ Jacquet-Langlands correspondence.

Let $D$ be the unique nonsplit quaternion algebra over $\Qp$. Let $\FF$ be a sufficiently large finite extension of $\Fp$, which will serve as the coefficient field. Let $G_{\Qp}:=\Gal(\Qpbar/\Qp)$ and
\[
\rhobar:G_{\Qp}\longrightarrow\GL_2(\FF)
\]
be a continuous generic (see Section~\ref{sec:application}) representation. Let $\pi(\rhobar)$ be the smooth admissible representation of $\GL_2(\Qp)$ attached to $\rhobar$ by the mod $p$ local Langlands correspondence. Let
\[
\cS^1(\pi(\rhobar))
:=H^1_{\mathrm{\acute et}}\left(\PP^1_{\CC_p},
\cF_{\pi(\rhobar)}\right),
\]
where $\cS^1$ is Scholze's functor constructed in \cite{Sch18}. It carries an action of $G_{\Qp}\times D^{\times}$. By \cite[Theorem 1.4]{HW24a}, there is a $G_{\Qp}\times D^\times$-equivariant isomorphism
\[
\cS^1(\pi(\rhobar))\cong\rhobar(-1)\boxtimes\mathrm{JL}(\rhobar)
\]
for some admissible smooth representation $\mathrm{JL}(\rhobar)$ of $D^\times$ with central character $\zeta=\det(\rhobar)\omega^{-1}$, where $\omega$ is the mod $p$ cyclotomic character and $\zeta$ is viewed as a character of $\Qp^\times$ via local class field theory. The $D^\times$-representation $\mathrm{JL}(\rhobar)$ is viewed as a candidate for the mod $p$ Jacquet-Langlands correspondence. Here is our main result.

\begin{theorem}[Theorem~\ref{thm:JL-self-extensions}]\label{thm:main}
    Suppose that $\rhobar$ is generic. Then we have
    \[
    \dim_\FF\Ext^i_{D^\times,\zeta}
    \left(\mathrm{JL}(\rhobar),\mathrm{JL}(\rhobar)\right)
    =
    \begin{cases}
    1,&i=0,\\
    3,&i=1,2;\\
    0,&i\ge3.
    \end{cases}
    \]
    Here the extension groups are computed in the category of smooth representations with fixed central character $\zeta$.
\end{theorem}

We explain the significance of the result from the viewpoint of deformation theory. See Section~\ref{sec:deformations} for details. Let $\cO$ be the ring of integers of a finite extension of $\Qp$ with residue field $\FF$. Fix a continuous character $\widetilde\zeta:Z_D\to\cO^\times$ lifting $\zeta$. Let $\Art_{\cO}$ be the category of local (possibly noncommutative) artinian $\cO$-algebras with residue field $\FF$ and central $\cO$-action. Following \cite{Pas13}, we can define the deformation functor $\operatorname{Def}_{\mathrm{JL}(\rhobar),\widetilde\zeta}$ sending $A\in\Art_{\cO}$ to the set of isomorphism classes of deformations of $\mathrm{JL}(\rhobar)^\vee$ to $A$. Since $\End_{D^\times}\left(\mathrm{JL}(\rhobar)\right)=\FF$, there exists a local pseudocompact $\cO$-algebra $R_{\mathrm{JL}(\rhobar),\widetilde\zeta}$ with residue field $\FF$ such that for any $A\in\Art_{\cO}$, we have a natural bijection
\[
\operatorname{Def}_{\mathrm{JL}(\rhobar),\widetilde\zeta}(A)
\cong
\Hom_{\cO,\mathrm{cont}}(R_{\mathrm{JL}(\rhobar),\widetilde\zeta},A)/A^\times,
\]
where $A^\times$ acts by conjugation. Moreover, the maximal commutative quotient $R_{\mathrm{JL}(\rhobar),\widetilde\zeta}^{\mathrm{ab}}$
pro-represents the restriction of $\operatorname{Def}_{\mathrm{JL}(\rhobar),\widetilde\zeta}$ to commutative algebras in the usual sense. By deformation theory, $\Ext^1_{D^\times,\zeta}\left(\mathrm{JL}(\rhobar),\mathrm{JL}(\rhobar)\right)$ corresponds to the tangent space of $R_{\mathrm{JL}(\rhobar),\widetilde\zeta}$, while $\Ext^2_{D^\times,\zeta}\left(\mathrm{JL}(\rhobar),\mathrm{JL}(\rhobar)\right)$ corresponds to the obstructions. As an application, we prove the following criterion for the commutativity of $R_{\mathrm{JL}(\rhobar),\widetilde\zeta}$.

\begin{theorem}[Theorem~\ref{thm:commutativity}]
    Suppose that $\rhobar$ is generic and there is a continuous local $\cO$-algebra surjection
    \[
    R_{\mathrm{JL}(\rhobar),\widetilde\zeta}\twoheadrightarrow
    \cO[[x_1,x_2,x_3]].
    \]
    Then this map is an isomorphism. In particular, $R_{\mathrm{JL}(\rhobar),\widetilde\zeta}$ is commutative.
\end{theorem}

We explain the proof strategy of Theorem~\ref{thm:main}. Let $\cO_D$ be the ring of integers of $D$. Fix a uniformizer $\varpi_D$ of $D$ satisfying $\varpi_D^2=p$. Put $U_D^1:=1+\varpi_D\cO_D$. Let $Z_D\cong\Qp^\times$ be the centre of $D^{\times}$ and put $Z_D^1:=Z_D\cap U_D^1$. Let $\Lambda:=\FF[[U_D^1/Z_D^1]]$ be the complete group algebra of $U_D^1/Z_D^1$, which is a local ring with maximal ideal $\mathfrak{m}_D=\ker(\Lambda\rightarrow \mathbb{F})$. It admits an action of $H:=\FF_{p^2}^\times$. Let $Y,Z\in\mathfrak{m}_D$ be two $H$-eigenvectors defined in \eqref{eq:YZ}, and let $y,z$ denote their symbols in the associated graded ring $\gr_{\mathfrak{m}_D}(\Lambda)$ respectively. Let $M=\mathrm{JL}(\rhobar)^{\vee}:=\Hom_{\FF}(\mathrm{JL}(\rhobar),\FF)$ be the Pontryagin dual of $\mathrm{JL}(\rhobar)$, which is a finitely generated $\Lambda$-module with compatible $D^{\times}$-action. By the work of \cite{Koh17,HW24a,HW24b,HM26}, we have known
\begin{enumerate}[label=\alph*)]
    \item the explicit structure of $\gr_{\mathfrak{m}_D}(M)$ as a $\gr_{\mathfrak{m}_D}(\Lambda)$-module with compatible $H$-action and $\varpi_D$-action. See Proposition~\ref{prop:JL-graded} and assumptions~\ref{it:assum-i}\ref{it:assum-i-red};
    \item $M$ is essentially self-dual. See Proposition~\ref{prop:JL-self-duality} and assumptions~\ref{it:assum-ii}\ref{it:assum-ii-red}.
\end{enumerate}
In fact, we prove Theorem~\ref{thm:main} for any admissible smooth representation $\pi$ of $D^\times$ with central character $\zeta$ satisfying the two assertions above.

Now we focus on the case where $\rhobar$ is irreducible. The case where $\rhobar$ is reducible is more or less similar. We first prove a general criterion (see Lemma~\ref{lem:decomposition-criterion}) about when a direct sum decomposition in the graded level can be lifted to the filtered level. Applying the criterion to $\gr_{\mathfrak{m}_D}(M)$ in \eqref{eq:irr-gr}, we obtain an $\cO_D^{\times}$-stable decomposition $M=M_z\oplus M_y$ (see Proposition~\ref{prop:O_D-decomposition}) where $M_z$ (resp. $M_y$) is free of rank 2 over $\FF[[Z]]$ (resp. $\FF[[Y]]$) and we reduce Theorem~\ref{thm:main} to the computation of $\Ext^i_\Lambda(M_z,M_z)^H$ (see Corollary~\ref{cor:Ext-decomposition} and Proposition~\ref{prop:mixed-ext-dim}). 

Choose an $H$-eigenbasis of $M_z$ as an $\FF[[Z]]$ module. By comparing the $H$-characters, one can write down the leading term \eqref{eq:leading-coefficient} of the action of $Y$, which is the first nontrivial extension relation. On the one hand, by \cite[Proposition 2.21]{HW24a}, the coefficient of the leading term must be nonzero. On the other hand, by our genericity assumption, the leading term is killed in $\gr_{\mathfrak{m}_D}(M_z)$. In order to detect the first nontrivial extension in the graded level, we introduce a weighted filtration $\cF$ (see Section~\ref{sec:weighted-filtration}) on $\Lambda$ and $M_z$, which plays an essential role throughout the whole paper. 

Using the weighted filtration, we define a $\Lambda$-submodule $M_z'\subseteq M_z$ of finite codimension, whose associated graded module (with respect to the weighted filtration) can be decomposed into the direct sum. Applying the decomposition criterion again, we obtain a filtered $(\Lambda,\cF)$-module decomposition $M_z'=N_+\oplus N_-$ compatible with $H'$ (see Proposition~\ref{prop:further-decomposition}), where $H'$ is the subgroup of $H$ with index 2 and $N_+$ (resp. $N_-$) is free of rank 1 over $\FF[[Z]]$. We also reduce the computation of $\Ext^i_\Lambda(M_z',M_z')^H$ to $\Ext^i_\Lambda(N_+,N_+)^{H'}$ (see Corollary~\ref{cor:mix-N-ext}).

Put $N=N_+$. We first convert $\Ext^2_\Lambda(N,N)^{H'}$ to $\Hom_\Lambda(N,N^\circ/N)^{H'}$ via local duality, where $N^\circ=\FF((Z))\otimes_{\FF[[Z]]}N$ (see Propositions~\ref{prop:dual}, \ref{prop:ext2-vanish}). We use the action of the congruence subgroups of $D^\times$ to estimate the order of poles, which gives $\Hom_\Lambda(N,N^\circ/N)^{H'}=0$ (see Proposition~\ref{prop:Hom-vanish}). This shows the vanishing of $\Ext^2_\Lambda(N,N)^{H'}$. To compute $\Ext^i_\Lambda(N,N)^{H'}$ for $i=0,1$, we use the convergent multiplicative filtered spectral sequence associated with the weighted filtration. Its first page is 
\[
\Ext^\bullet_{\gr_\cF(\Lambda)}\left(\gr_\cF(N),\gr_\cF(N)\right)^{H'},
\]
which can be computed explicitly (see Lemma~\ref{lem:E1-C}). The vanishing of $\Ext^2_\Lambda(N,N)^{H'}$ and the nondegenerate pairing between $\Ext^2_\Lambda(N,\Lambda)$ and $N^\circ/N$ induced by local duality restrict the differentials of the spectral sequence and force the multiplicative spectral sequence to have specific structure for every page. See the proof of Proposition~\ref{prop:ext-N} for details. Eventually we obtain $\End_\Lambda(N)^{H'}=\FF$ and $\Ext^1_\Lambda(N,N)^{H'}=0$, which also gives $\Ext^i_\Lambda(M_z',M_z')^H$. Finally we return to $\Ext^i_\Lambda(M_z,M_z)^H$ via Euler-Poincar\'e characteristic and finish the proof of Theorem~\ref{thm:main} in the case where $\rhobar$ is irreducible (see Theorem~\ref{thm:self-extension}).

\subsection*{Acknowledgments}
The author is deeply grateful to his advisor, Professor Yongquan Hu, for suggesting the problem studied in this paper and for his patient guidance and constant encouragement throughout this project.
 
\subsection*{AI Use Disclosure}
Some key ideas of the paper are given by GPT-6 Astra, including
\begin{enumerate}[label=\alph*)]
    \item the decomposition criterion Lemma~\ref{lem:decomposition-criterion}. The author previously expected a decomposition of $\mathrm{JL}(\rhobar)$ using the property of Scholze's functor and $\pi(\rhobar)$, but Lemma~\ref{lem:decomposition-criterion} covers a broader range of cases;
    \item the weighted filtration in Section~\ref{sec:weighted-filtration}. At first, the author wanted to better understand \cite[Proposition 2.21]{HW24a}. After several rounds of chatting, GPT introduced the striking weighted filtration presented in the paper;
    \item the spectral sequence computation of $\Ext^1_\Lambda(N,N)^{H'}$ in the proof of Proposition~\ref{prop:ext-N}.
\end{enumerate}
The author digests all these ideas and rewrites the paper \textbf{by hand}. The author takes full responsibility for this article.

\section{Notation and preliminaries}

Throughout the paper, we fix a prime $p\ge5$. Let $\FF$ be a sufficiently large finite extension of $\FF_p$. All modules are left modules without special emphasis.

For an $R$-module $M$ and $i\ge0$, we write $\E_R^i(M):=\Ext_R^i(M,R)$, which is a right $R$-module. We say that $M$ is Cohen-Macaulay of grade $d$ if $d$ is the unique index such that $\E_R^d(M)\ne0$.

\subsection{Preliminaries on filtered and graded modules}

We use decreasing filtrations throughout the paper. A filtered ring $R$ is equipped with additive subgroups $F^nR~(n\in\ZZ)$ of $R$ satisfying $1\in F^0R$, $F^nR\supseteq F^{n+1}R$ and $F^iR\cdot F^jR\subseteq F^{i+j}R$ for any $i,j\in\mathbb{Z}$. Filtered rings appearing in this paper always satisfy $F^0 R=R$. We write $\gr(R):=\oplus_{n\in\ZZ}F^nR/F^{n+1}R$ for the graded ring.

A filtered $R$-module $X$ is equipped with submodules $F^nX~(n\in\ZZ)$ satisfying $F^nX\supseteq F^{n+1}X$ and $F^iR\cdot F^jX\subseteq F^{i+j}X$ for any $i,j\in\mathbb{Z}$. Then $\gr(X):=\oplus_{n\in\ZZ}F^nX/F^{n+1}X$ is a graded $\gr(R)$-module. We write $\gr^n(X):=F^nX/F^{n+1}X$.
We say $X$ (or the filtration on $X$) is complete if $X\to\varprojlim_{n\in\ZZ} X/F^nX$ is an isomorphism, separated if $\bigcap_nF^nX=0$, exhaustive if $\bigcup_nF^nX=X$ and bounded below if $F^nX=X$ for $n\ll0$. A complete filtration is separated, and a filtration bounded below is exhaustive. If $X$ is separated and exhaustive, and $x\in X$ is a nonzero element, then there is a unique integer $n$ such that $x\in F^nX\setminus F^{n+1}X$. We call $n$ the degree of $x$ and write $\gr^n(x)\in\gr^n(X)$ for its symbol. When $R$ and $X$ admit actions of a group $\Gamma$ preserving their filtrations, we say the actions are compatible if 
\[
\gamma\cdot(rx)=(\gamma\cdot r)(\gamma\cdot x)
\]
for any $\gamma\in\Gamma,~r\in R,~x\in X$. We use analogous conventions in the graded setting.

For a filtered (resp. graded) module $X$ and $k\in\ZZ$, our shift convention is $F^n(X(k)):=F^{n+k}X$ (resp. $X(k)^n:=X^{k+n}$). A filtered $R$-module $X$ is called filt-free if $X\cong\bigoplus_{i=1}^nR(-k_i)$ for some integers $k_i$. A graded $\gr(R)$-module $X$ is called gr-free if $X\cong\bigoplus_{i=1}^n\gr(R)(-k_i)$ for some integers $k_i$. 

Let $X,Y$ be two filtered (resp. graded) modules. We say a filtered (resp. graded) morphism $f:X\to Y$ has degree $k$ if $f(F^nX)\subseteq F^{n+k}Y$ (resp. $f(X^n)\subseteq Y^{n+k}$) for all $n$. We require that the differentials in a filt-free (resp. gr-free) resolution always have degree zero. A filtered morphism $f:X\to Y$ is strict if $f(F^nX)=\Im(f)\cap F^nY$ for all $n$. The filtration on $\Hom_R(X,Y)$ is
\[
F^n\Hom_R(X,Y):=\{f\in\Hom_R(X,Y):f(F^jX)\subseteq F^{j+n}Y \quad\text{for all}\quad j\in\ZZ\}.
\]

Let $C^\bullet$ be a cochain complex of filtered $\FF$-vector spaces. Then it induces a filtration on the cohomology $H^\bullet(C^\bullet)$ and determines a spectral sequence $E_\bullet^{\bullet\bullet}(C)$ starting from $E_0^{i,j}(C)=F^iC^{i+j}/F^{i+1}C^{i+j}$ and $E_1^{i,j}(C)=H^{j}(E_0^{i,\bullet})=H^{i+j}\left(\gr^i(C^\bullet)\right)$, see \cite[\S 5.4]{WC}.

\begin{proposition}\label{prop:spectral-sequence}
    Suppose that for each $n$, the filtration on $C^n$ is complete and bounded below with finite-dimensional graded pieces. Then
    \[
    E_\bullet^{i,j}(C)\Longrightarrow H^{i+j}(C^{\bullet}).
    \]
\end{proposition}
\begin{proof}
    This is a special case of \cite[Theorem 5.5.10 (2)]{WC}. Note that \cite{WC} uses increasing filtrations, so bounded above loc. cit. is bounded below in our convention. A descending chain of finite-dimensional $\FF$-vector spaces must be stable, which implies that the spectral sequence is regular.
\end{proof}

\subsection{A decomposition lemma}

Throughout this subsection, $R$ is a complete filtered ring with $F^0R=R$.

\begin{lemma}\label{lem:resolution-lifting}
    Let $X$ be a complete and exhaustive filtered $R$-module. Assume that $\gr(X)$ admits a bounded gr-free $\gr(R)$-module resolution $(\widetilde{P}_{\bullet},\widetilde{\partial}_{\bullet})$ whose terms are of finite rank. Then $(\widetilde{P}_{\bullet},\widetilde{\partial}_{\bullet})$ can be lifted to a strict bounded filt-free $R$-module resolution $(P_{\bullet},\partial_{\bullet})$ of $X$, i.e., $\gr(P_{\bullet},\partial_{\bullet})\cong (\widetilde{P}_{\bullet},\widetilde{\partial}_{\bullet})$.
\end{lemma}

\begin{proof}
The proof is similar to that of \cite[Corollary I.7.2.9]{LvO96}.
Write 
\[
\widetilde{P}_i
=\bigoplus_s\gr(R)\widetilde{e}_{i,s}\cong\bigoplus_s\gr(R)(-d_{i,s}),
\]
where $\widetilde{e}_{i,s}$ is homogeneous of degree $d_{i,s}$. Let 
\[
P_i=\bigoplus_sR e_{i,s}\cong\bigoplus_sR(-d_{i,s})
\]
be the corresponding filt-free $R$-module, with
$\deg(e_{i,s})=d_{i,s}$. For each $s$, lift $\widetilde\partial_0(\widetilde{e}_{0,s})$ to some $x_{0,s}\in F^{d_{0,s}}X$. The assignment $e_{0,s}\mapsto x_{0,s}$ defines a filtered $R$-morphism 
\[
\partial_0:P_0\longrightarrow X
\]
of degree zero, which is strict and surjective because $\gr(\partial_0)=\widetilde\partial_0$ is surjective and $X$ is complete and exhaustive. Put $K_0=\ker\partial_0$, which is complete and exhaustive under the induced filtration. By strictness we have 
\[
\gr(K_0)\cong\ker\widetilde\partial_0
=\Im\widetilde\partial_1.
\]
Replacing $X$ by $K_0$ and repeating the procedure above, we can lift $\widetilde{P}_{\bullet}$ strictly step by step. Assume the last term of $\widetilde{P}_{\bullet}$ is $\widetilde{P}_n$ and we have lifted $\widetilde\partial_n$ strictly to 
\[
\partial_n:P_n\longrightarrow P_{n-1}.
\]
Put $K_n=\ker\partial_n$. By strictness we have 
\[
\gr(K_n)\cong\ker\widetilde\partial_n=0,
\]
which implies $K_n=0$ because the filtration is separated.
\end{proof}

\begin{remark}\label{rmk:H-compatible}
    If both $X$ and $(\widetilde{P}_{\bullet},\widetilde{\partial}_{\bullet})$ admit compatible actions of some finite group $\Gamma$ whose order is invertible in $R$, then by averaging the lift $(P_{\bullet},\partial_{\bullet})$ can be chosen with a compatible $\Gamma$-action. 
\end{remark}

The following lemma will be repeatedly used throughout the paper.

\begin{lemma}\label{lem:decomposition-criterion}
    Let $X$ be a complete and exhaustive filtered $R$-module. Assume that 
    \begin{enumerate}
        \item $\gr(X)=\widetilde{X}^1\oplus\widetilde{X}^2$ as a graded $\gr(R)$-module;
        \item each $\widetilde{X}^i$ admits a bounded gr-free resolution whose terms are of finite rank;
        \item for any $j>0$, we have
    \begin{equation}\label{eq:graded-Ext-vanish}
        \Ext^1_{\gr(R)}(\widetilde{X}^1,\widetilde{X}^2)^j=\Ext^1_{\gr(R)}(\widetilde{X}^2,\widetilde{X}^1)^j=0
    \end{equation}
    (by \textup{(ii)}, $\Ext^\bullet_{\gr(R)}(\widetilde{X}^1,\widetilde{X}^2)$ and $\Ext^\bullet_{\gr(R)}(\widetilde{X}^2,\widetilde{X}^1)$ are graded abelian groups).
    \end{enumerate}
    Then there exist closed filtered $R$-submodules $X^1,X^2\subseteq X$
    such that 
    \[
    X=X^1\oplus X^2,\quad\gr(X^i)\cong \widetilde{X}^i,\quad i=1,2.
    \]
    Moreover, if for any $j\ge0$, we have
    \begin{equation}\label{eq:graded-Hom-vanish}
        \Hom_{\gr(R)}(\widetilde{X}^1,\widetilde{X}^2)^j=\Hom_{\gr(R)}(\widetilde{X}^2,\widetilde{X}^1)^j=0,
    \end{equation}
    then such a decomposition is unique.
\end{lemma}

\begin{proof}
Choose bounded gr-free $\gr(R)$-module resolutions whose terms are of finite rank
\[
\widetilde{P}^i_\bullet:
\quad
\cdots\xrightarrow{\widetilde{\partial}^i_2}
\widetilde{P}^i_1
\xrightarrow{\widetilde{\partial}^i_1}
\widetilde{P}^i_0
\longrightarrow
\widetilde{X}^i
\longrightarrow 0,\quad i=1,2.
\]
Then $\widetilde{P}_\bullet:= \widetilde{P}^1_\bullet\oplus\widetilde{P}^2_\bullet$ is a gr-free $\gr(R)$-module resolution of $\widetilde{X}^1\oplus\widetilde{X}^2
= \gr(X)$. By Lemma \ref{lem:resolution-lifting} and its proof, we can lift $\widetilde{P}_\bullet$ to a strict bounded filt-free $R$-module resolution
\[
P_\bullet:
\quad
\cdots
\xrightarrow{\partial_2}
P_1
\xrightarrow{\partial_1}
P_0
\longrightarrow
X
\longrightarrow 0
\]
with each $P_n$ a direct sum of filt-free $R$-modules $P^1_n$ and $P^2_n$ such that $\gr(P^i_n)\cong \widetilde{P}^i_n$. The differential $\partial_{\bullet}$ has the block form
\[
    \partial_n=
    \begin{pmatrix}
        \partial^1_n & b_n\\
        c_n & \partial^2_n
    \end{pmatrix}
    :
    P^1_n\oplus P^2_n
    \longrightarrow
    P^1_{n-1}\oplus P^2_{n-1},
\]
where $\gr(\partial^i_n)=\widetilde{\partial}^i_n$ and 
\[
b_n: P^2_n\longrightarrow P^1_{n-1},\quad c_n: P^1_n\longrightarrow P^2_{n-1}
\]
have positive degree.

Expanding the identity $\partial_{n-1}\circ\partial_n=0$ gives
\begin{equation}\label{eq:differential-equation}
    \partial^1_{n-1}\circ b_n+b_{n-1}\circ\partial^2_n=0,\quad c_{n-1}\circ\partial^1_n+\partial^2_{n-1}\circ c_n=0.
\end{equation}
Suppose that $b_n,c_n$ have degree $j>0$ for all $n$. Let
\[
b_{n,j}:\widetilde{P}^2_n\longrightarrow\widetilde{P}^1_{n-1},\quad c_{n,j}:\widetilde{P}^1_n\longrightarrow\widetilde{P}^2_{n-1}
\]
be the graded morphisms of degree $j$ induced by $b_n$ and $c_n$ respectively. Then \eqref{eq:differential-equation} gives
\begin{equation}\label{eq:graded-differential-equation}
    \widetilde{\partial}^1_{n-1}\circ b_{n,j}+b_{n-1,j}\circ\widetilde{\partial}^2_n=0,\quad c_{n-1,j}\circ\widetilde{\partial}^1_n+\widetilde{\partial}^2_n\circ c_{n,j}=0.
\end{equation}
Consider the $\Hom$ complex $\cH_{21}^\bullet:=\Hom_{\gr(R)}^\bullet(\widetilde{P}^2_\bullet,\widetilde{P}^1_\bullet)$ with cohomological degree $r$ part 
\[
\cH_{21}^r=\prod_n\Hom_{\gr(R)}(\widetilde{P}^2_n,\widetilde{P}^1_{n-r})
\]
and differential 
\[
(\delta_{21}^r f)_n=\widetilde{\partial}^1_{n-r}\circ f_n-(-1)^r f_{n-1}\circ\widetilde{\partial}^2_n
\]
for $f=(f_n)_n\in \cH_{21}^r$, see \cite[\S2.7.4, 2.7.5]{WC}. Put $b_j:=(b_{n,j})_n\in\cH_{21}^1$. Then the first equality in \eqref{eq:graded-differential-equation} says precisely $\delta_{21}^1(b_j)=0$. By \cite[Lemma 2.7.3]{WC} and the proof of \cite[Theorem 2.7.2]{WC}, the complex $\cH_{21}^\bullet$ is quasi-isomorphic to $\Hom_{\gr(R)}(\widetilde{P}^2_\bullet,\widetilde{X}^1)$, and consequently
\[
    H^1
    \Hom_{\gr(R)}^\bullet(\widetilde{P}^2_\bullet,\widetilde{P}^1_\bullet)
    \cong
    \Ext^1_{\gr(R)}(\widetilde{X}^2,\widetilde{X}^1).
\]
By \eqref{eq:graded-Ext-vanish}, the degree-$j$ parts of the two graded abelian groups vanish. Hence there exists $u_j=(u_{j,n})_n\in \cH_{21}^0$, with degree $j$, such that
\[
    b_{n,j}=(\delta_{21}^0 u_j)_n
    =
    \widetilde{\partial}^1_n u_{j,n}
    -
    u_{j,n-1}\widetilde{\partial}^2_n.
\]
Similarly, there exists $v_j=(v_{j,n})_n\in \cH_{12}^0$, with degree $j$, such that
\[
    c_{n,j}=(\delta_{12}^0 v_j)_n
    =
    \widetilde{\partial}^2_n v_{j,n}
    -
    v_{j,n-1}\widetilde{\partial}^1_n.
\]

Put 
\[
\widetilde{S}_{j,n}
    :=
    \begin{pmatrix}
        0 & u_{j,n}\\
        v_{j,n} & 0
    \end{pmatrix}
    :\widetilde{P}^1_n\oplus\widetilde{P}^2_n
    \longrightarrow
    \widetilde{P}^1_n\oplus\widetilde{P}^2_n.
\]
Then
\begin{equation}\label{eq:Sjn-equality}
    \begin{pmatrix}
        0 & b_{n,j}\\
        c_{n,j} & 0
    \end{pmatrix}
    =
       \begin{pmatrix}
        \widetilde{\partial}^1_n & 0\\
        0 & \widetilde{\partial}^2_n
    \end{pmatrix}   
    \widetilde{S}_{j,n}-\widetilde{S}_{j,n-1}
    \begin{pmatrix}
        \widetilde{\partial}^1_n & 0\\
        0 & \widetilde{\partial}^2_n
    \end{pmatrix}.
\end{equation}
Lift each $\widetilde{S}_{j,n}$ to a filtered endomorphism $S_{j,n}:P_n\rightarrow P_n$ of degree $j$. Put $U_{j,n}=1+S_{j,n}$. Since $j>0$, by completeness $U_{j,n}$ is invertible with inverse $\sum_{k\ge0}(-S_{j,n})^k$. Put 
\[
P^1_{j,n}=U_{j,n}^{-1}(P_n^1),\quad P^2_{j,n}=U_{j,n}^{-1}(P_n^2).
\]
Then $P^i_{j,n}$ is filt-free equipped with the induced filtration and $P^i_{j,n}\cong P_n^i$ as filt-free $R$-modules. Moreover, with respect to the decomposition $P^1_{j,n}\oplus P^2_{j,n}$, the differential $\partial_n$ has the block form 
\[
    U_{j,n-1}
    \begin{pmatrix}
        \partial^1_n & b_n\\
        c_n & \partial^2_n
    \end{pmatrix}
    U_{j,n}^{-1}=
    \begin{pmatrix}
        {\partial^1_n}' & {b_n}'\\
        {c_n}' & {\partial^2_n}'
    \end{pmatrix},
\]
where $\gr({\partial^i_n}')=\widetilde{\partial}^i_n$ and ${b_n}',{c_n}'$ have degree at least $j+1$ by \eqref{eq:Sjn-equality}.

Repeating the preceding construction, we obtain inductively $U_{k,n}=1+S_{k,n}~(k>j)$ with $S_{k,n}:P_n\rightarrow P_n$ filtered of degree $k$. For $k>j$, put 
\[
U_n^k:=U_{k,n}\circ U_{k-1,n}\circ\cdots\circ U_{j,n}.
\]
Then the sequence $\{U_n^k\}_{k>j}$ converges to some $U_{\infty,n}$ by completeness, which is invertible. Put 
\[
P^1_{\infty,n}=U_{\infty,n}^{-1}(P_n^1),\quad P^2_{\infty,n}=U_{\infty,n}^{-1}(P_n^2).
\]
Then $P^i_{\infty,n}$ is filt-free equipped with the induced filtration and $P^i_{\infty,n}\cong P_n^i$ as filt-free $R$-modules. With respect to the decomposition $P^1_{\infty,n}\oplus P^2_{\infty,n}$, the differential $\partial_n$ has the block form $U_{\infty,n-1}\partial_n U_{\infty,n}^{-1}$, which is diagonal because $P_n$ is separated. Consequently, the filt-free resolution $P_\bullet$ decomposes into $P^1_{\infty,\bullet}\oplus P^2_{\infty,\bullet}$. The first assertion follows by taking $X^i=H_0(P^i_{\infty,\bullet}),~i=1,2$.

Assume that \eqref{eq:graded-Hom-vanish} holds and we have two such decompositions 
\[
X=X^1\oplus X^2={X^1}'\oplus {X^2}'.
\]
Consider the filtered $R$-morphism
\[
X^1\hookrightarrow X\twoheadrightarrow {X^2}'.
\]
If it is nonzero, then it induces a nonzero graded $\gr(R)$-morphism $\gr(X^1)\to\gr({X^2}')$ with some degree $j\ge0$, which contradicts \eqref{eq:graded-Hom-vanish}. Therefore we have $X^1\subseteq{X^1}'$. Similarly we have ${X^1}'\subseteq X^1$, which implies $X^1={X^1}'$. This gives the uniqueness.
\end{proof}

\subsection{Quaternionic notation}

The quaternionic notations used here are from \cite{HW24a} and \cite{HW24b}. Let $D$ be the unique non-split quaternion algebra over $\Qp$ with ring of integers $\cO_D$. Fix a uniformizer $\varpi_D$ of $D$ satisfying $\varpi_D^2=p$ and $\varpi_D^{-1}a\varpi_D=\sigma(a)$ for $a\in\mathbb{Q}_{p^2}$, where $\sigma$ is the Frobenius element on $\mathbb{Q}_{p^2}$. For $n\geq1$, put 
\[
U_D^n:=1+\varpi_D^n\cO_D.
\]
Then $U_D^1$ is the unique pro-$p$ Sylow subgroup of $\cO_D^{\times}$, and we have 
\[
\cO_D^{\times}/U_D^1\cong \Fq^{\times}
\]
where $q:=p^2$. Let $[\cdot]$ denote the Teichm\"uller lifting 
\[
\Fq\hookrightarrow \mathbb{Z}_q\hookrightarrow \cO_D,
\]
and $H<\cO_D^{\times}$ denote the image of $\Fq^{\times}$. Let 
\[
\alpha:H\rightarrow\mathbb{F}^{\times}
\]
denote the character of $H$ sending $x$ to $\overline{x}^{p-1}$, which has order $p+1$. 

Let $Z_D\cong\Qp^\times$ be the centre of $D^{\times}$. Put $Z_D^1:=Z_D\cap U_D^1$. Let 
\[
\Lambda:=\FF[[U_D^1/Z_D^1]]
\]
be the complete group algebra of $U_D^1/Z_D^1$. It is a local ring with maximal ideal 
\[
\mathfrak{m}_D=\ker(\Lambda\rightarrow \mathbb{F}).
\]
The group $H$ and the uniformizer $\varpi_D$ act on $\Lambda$ by conjugation. Put
\begin{equation}\label{eq:YZ}
 Y=\sum_{a\in\FF_q^\times}a^{-1}(1+\varpi_D[a]),\quad
 Z=\sum_{a\in\FF_q^\times}a^{-p}(1+\varpi_D[a]).
\end{equation}
Then $Y,Z\in\mathfrak{m}_D$. The group $H$ acts on $Y$ (resp. $Z$) via $\alpha$ (resp. $\alpha^{-1}$) and $\varpi_D$ permutes $Y$ and $Z$. Let 
\[
y=\gr_{\mathfrak{m}_D}^1(Y),\quad z=\gr_{\mathfrak{m}_D}^1(Z)
\]
be the symbols of $Y$ and $Z$ respectively. Put $h=[y,z]$. Then by \cite[\S2.2, 2.3]{HW24b}, we have 
\[
\gr_{\mathfrak{m}_D}(\Lambda)\cong\FF\langle y,z,h\rangle
\]
with relations $[y,z]=h,~[y,h]=[z,h]=0$ and $\deg (y)=\deg(z)=1,~\deg(h)=2$.

\section{Weighted filtration}\label{sec:weighted-filtration}

Let $\FF\langle\langle Y',Z'\rangle\rangle$ be the complete free algebra of two variables, with $H$ acting on $Y'$ (resp. $Z'$) by $\alpha$ (resp. $\alpha^{-1}$) and $\varpi_D$ permuting $Y'$ and $Z'$. Let $J=(Y',Z')$ be the augmentation ideal of $\FF\langle\langle Y',Z'\rangle\rangle$. Then $\gr_J\FF\langle\langle Y',Z'\rangle\rangle\cong\FF\langle y',z'\rangle$ is the free algebra of two variables. 

Consider the filtered algebra homomorphism 
\[
\Phi:\FF\langle\langle Y',Z'\rangle\rangle\longrightarrow\Lambda
\]
sending $Y'$ to $Y$ and $Z'$ to $Z$. Then $\Phi$ is strict and surjective because 
\[
\gr\Phi:\FF\langle y',z'\rangle\longrightarrow\gr_{\mathfrak{m}_D}(\Lambda)
\]
is surjective. Recall that $\ker(\gr\Phi)$ is the ideal generated by $[y',[y',z']]$ and $[z',[y',z']]$. By strictness of $\Phi$, there exists $f(Y',Z')\in J^4$ such that 
\[
[Y',[Y',Z']]-f(Y',Z')\in\ker\Phi.
\]
The action of $\varpi_D$ gives 
\[
[Z',[Y',Z']]+f(Z',Y')\in\ker\Phi.
\]
By averaging we can assume that $H$ acts on $f(Y',Z')$ by $\alpha$. Separate the monomial $Z'^p$ from $f(Y',Z')$, say, $f(Y',Z')=cZ'^p+g(Y',Z')$ with $g(Y',Z')\in J^4$. By completeness, $\ker\Phi$ is the closed ideal generated by 
\begin{equation}\label{eq:Lambda-relation}
   R_1:= [Y',[Y',Z']]-cZ'^p-g(Y',Z'),\quad R_2:=[Z',[Y',Z']]+cY'^p+g(Z',Y').
\end{equation}

\begin{lemma}\label{lem:c-nonzero}
The coefficient $c$ in \eqref{eq:Lambda-relation} is nonzero.
\end{lemma}

\begin{proof}
    Consider 
    \[
    L=\FF\langle\langle Y',Z'\rangle\rangle/\overline{(Z',{Y'}^{p+1})}\cong\FF[Y']/({Y'}^{p+1}),
    \]
    which is a filtered $\FF\langle\langle Y',Z'\rangle\rangle$-module with compatible $H$-action. If $c=0$, then $R_1$ and $R_2$ will kill $L$ because the degree of any $Y'$-monomial appearing in $g(Y',Z')$ (resp. $g(Z',Y')$) is at least $p+2$ (resp. $2p+1$). Therefore, $L$ is a filtered $\Lambda$-module with compatible $H$-action, satisfying $\dim_{\FF}L/YL=1$. This contradicts \cite[Proposition 2.21]{HW24a}, which says that such a module has dimension at most $p$.
\end{proof}

In order to detect the nontrivial extension in the associated graded module, we define a weighted filtration $\cF^n\Lambda$ on $\Lambda$, which will play an essential role throughout the paper. By the PBW theorem, $\gr_{\mathfrak{m}_D}(\Lambda)$ admits a graded basis $\{z^ih^jy^k\}_{i,j,k\ge0}$. By completeness, $\Lambda$ admits a filtered topological basis $\{Z^i[Y,Z]^jY^k\}_{i,j,k\ge0}$. We set 
\begin{equation}\label{eq:wt}
\wet(Y)=\frac{p-1}{2},\quad\wet(Z)=1,\quad\wet([Y,Z])=\frac{p+1}{2}
\end{equation}
and 
\[
\wet(Z^i[Y,Z]^jY^k)=i\wet(Z)+j\wet([Y,Z])+k\wet(Y),\quad i,j,k\ge0.
\]
The filtration $\cF^n\Lambda$ consists of convergent sums of $Z^i[Y,Z]^jY^k$ with weight at least $n$. Note that the weights in \eqref{eq:wt} are comparable with the $\mathfrak{m}_D$-adic degrees $1,1,2$, so the filtration $\cF$ is complete and defines the same ($\mathfrak{m}_D$-adic) topology on $\Lambda$. 

To show the filtration $\cF$ is multiplicative, i.e., $\cF^i\Lambda\cdot\cF^j\Lambda\subseteq\cF^{i+j}\Lambda$, we define a weighted filtration $\cF'$ on $\FF\langle\langle Y',Z'\rangle\rangle$. Set 
\[
\wet'(Y')=\frac{p-1}{2},\quad\wet'(Z')=1
\]
and for a monomial $X$ with $m_1$ copies of $Y'$ and $m_2$ copies of $Z'$, set 
\[
\wet
'(X)=m_1\wet'(Y')+m_2\wet'(Z').
\]
The filtration ${\cF'}^n$ consists of convergent sums of monomials with weight at least $n$. Then the filtration $\cF'$ is comparable with the $J$-adic filtration, and hence is complete and defines the same topology on $\FF\langle\langle Y',Z'\rangle\rangle$. Moreover, the filtration $\cF'$ is multiplicative. We define $\wet'(X)$ for a monomial $X\in\FF\langle y',z'\rangle$ in the same way.

\begin{lemma}\label{lem:g-weight}
    $g(Y',Z'),g(Z',Y')\in\cF'^{p+2}$.
\end{lemma}

\begin{proof}
    This follows from $g(Y',Z'),g(Z',Y')\in J^4$ and the fact that $H$ acts on $g(Y',Z')$ (resp. $g(Z',Y')$) by $\alpha$ (resp. $\alpha^{-1}$). For example, those monomials possibly appearing in $g(Z',Y')$ with the minimal weight are $Z'^{p+2}$ and words containing three copies of $Z'$ and two copies of $Y'$, which all have weight $p+2$.
\end{proof}

\begin{lemma}\label{lem:weight-multiplicative}
    We have $\Phi({\cF'}^n)=\cF^n\Lambda$ for $n\ge0$. In particular, the filtration $\cF$ is multiplicative.
\end{lemma}

\begin{proof}
We have $\Phi({\cF'}^n)\supseteq\cF^n\Lambda$. Let $f\in\cF'^n$. Assume $\deg_{J}(f)=m$. Let $\bar{f}$ be the symbol of $f$ in $\gr_J\FF\langle\langle Y',Z'\rangle\rangle\cong\FF\langle y',z'\rangle$. Then $\wet'(X)\ge n$ for any monomial $X$ appearing in $\bar{f}$. By the PBW theorem, $\gr\Phi(\bar{f})$ can be uniquely written as 
\[
\sum_{i,j,k\ge0}c_{i,j,k}z^ih^jy^k\in\gr^m_{\mathfrak{m}_D}(\Lambda),
\]
using the relations $[y,z]=h,~[y,h]=[z,h]=0$. Therefore, for each $c_{i,j,k}\ne0$, we have $i+2j+k=m$ and $i+\frac{p+1}{2}j+\frac{p-1}{2}k\ge n$. Put 
\[
P_m=\sum_{i,j,k}c_{i,j,k}Z'^i[Y',Z']^jY'^k\in J^m\cap{\cF'}^n
\]
and $p_m=\sum_{i,j,k}c_{i,j,k}z'^i[y',z']^jy'^k$. Then we have $\Phi(P_m)\in\cF^n\Lambda$ and
\[
\bar{f}-p_m\in\ker\gr\Phi=([y',[y',z']],[z',[y',z']]).
\]
Since $[y',[y',z']]$ and $[z',[y',z']]$ are homogeneous with respect to both $\deg_J$ and $\wet'$, $\bar{f}-p_m$ can be written as 
\[
\sum_i u_{1,i}(y',z')[y',[y',z']]v_{1,i}(y',z')+\sum_i u_{2,i}(y',z')[z',[y',z']]v_{2,i}(y',z')
\]
for some monomials $u_{j,i}$ and $v_{j,i}$ such that 
\[
\deg_J(u_{j,i})+\deg_J(v_{j,i})+3\ge m,\quad\wet'(u_{1,i})+\wet'(v_{1,i})+p\ge n,\quad\wet'(u_{2,i})+\wet'(v_{2,i})+\frac{p+3}{2}\ge n
\]
for any $i$ and $j=1,2$. Put 
\[
Q_m=\sum_i u_{1,i}(Y',Z')R_1v_{1,i}(Y',Z')+\sum_i u_{2,i}(Y',Z')R_2v_{2,i}(Y',Z')\in\ker\Phi.
\]
By Lemma~\ref{lem:g-weight}, we have $Q_m\in J^m\cap{\cF'}^n$, and the symbol of $Q_m$ in $\FF\langle y',z'\rangle$ is equal to $\bar{f}-p_m$. Put $f'=f-P_m-Q_m$. Then $f'\in J^{m+1}\cap{\cF'}^n$. Replacing $f$ by $f'$ and repeating the procedure above, we can construct inductively $P_k\in J^k\cap{\cF'}^n$ and $Q_k\in J^k\cap{\cF'}^n\cap\ker\Phi$ for $k>m$, such that $f=\sum_{k\ge m}P_k+\sum_{k\ge m}Q_k$. Therefore 
\[
\Phi(f)=\sum_{k\ge m}\Phi(P_k)\in\cF^n\Lambda.
\]
\end{proof}

We write
\[
\widetilde{\Lambda}=\gr_{\cF}(\Lambda),\quad\widetilde{y}=\gr_{\cF}(Y),\quad \widetilde{z}=\gr_{\cF}(Z),\quad
 \widetilde{h}=\gr_{\cF}([Y,Z]).
\]

\begin{proposition}\label{prop:graded-algebra}
We have 
\[
\widetilde{\Lambda}=\FF[\widetilde{z},\widetilde{h}][\widetilde{y};\partial=\widetilde{h}\frac{\partial}{\partial\widetilde{z}}+c\widetilde{z}^p\frac{\partial}{\partial\widetilde{h}}
],
\]
where $\FF[\widetilde{z},\widetilde{h}][\widetilde{y};\partial]$ is the skew polynomial ring over the polynomial ring $\FF[\widetilde{z},\widetilde{h}]$, satisfying $[\widetilde{y},f]=\partial(f)$ for $f\in\FF[\widetilde{z},\widetilde{h}]$.
This is a left and right noetherian domain of global dimension at most three.
\end{proposition}

\begin{proof}
    By the relations \eqref{eq:Lambda-relation} and Lemmas~\ref{lem:g-weight}, \ref{lem:weight-multiplicative}, we have $[\widetilde{y},\widetilde{z}]=\widetilde{h},~[\widetilde{y},\widetilde{h}]=c\widetilde{z}^p$ and $[\widetilde{z},\widetilde{h}]=0$. By the universal property of the skew polynomial ring (see \cite[\S1.2.5]{MR01}), we have a surjection 
    \[
    \FF[\widetilde{z},\widetilde{h}][\widetilde{y};\partial]\longrightarrow\widetilde{\Lambda},
    \]
    which is an isomorphism by comparing the bases on both sides. The stated ring-theoretical properties follow by \cite[Theorem 1.2.9, Theorem 7.5.3]{MR01}.
\end{proof}

\begin{remark}
    We see that $c$ is a constant, independent of the choice of $R_1$ and $R_2$. Later we will prove that $c=-2$ in Corollary~\ref{cor:c}.
\end{remark}

We say an operator $T:\Lambda\to\Lambda$ has order at least $d$ if $T$ sends $\cF^n\Lambda$ to $\cF^{n+d}\Lambda$ for every $n$.

\begin{lemma}\label{lem:ad-order}
    The operator $\ad_{Z^p}$ sending $X\in\Lambda$ to $[Z^p,X]$ has order at least $p+\frac{p+1}{2}$.
\end{lemma}

\begin{proof}
    The relation $R_2$ in \eqref{eq:Lambda-relation} and Lemma~\ref{lem:g-weight} give 
    \[
    (\ad_Z)^2(Y)\in\cF^{p+2}\Lambda.
    \]
    Therefore, we have 
    \[
    \ad_{Z^p}(Y)=(\ad_Z)^p(Y)\in\cF^{2p}\Lambda.
    \]
    Moreover, we have 
    \[
    \ad_{Z^p}([Y,Z])=-(\ad_Z)^{p+1}(Y)\in\cF^{2p+1}\Lambda.
    \]
    The result follows by the definition of the filtration $\cF$ and the Leibniz rule
    \[
    \ad_{Z^p}(X_1X_2)=\ad_{Z^p}(X_1)X_2+X_1\ad_{Z^p}(X_2).
    \]
\end{proof}

\begin{proposition}\label{prop:Laurent-extension}
Let $L$ be a filtered $(\Lambda,\cF)$-module, filt-free of finite rank over $\FF[[Z]]$. Let $L^\circ:=\FF((Z))\otimes_{\FF[[Z]]}L$, equipped with the natural tensor-product filtration. Then there is a unique filtered extension of the action of $\Lambda$, given by
\begin{equation}\label{eq:Laurent-limit}
 X^{\circ}(v):=\lim_{j\to\infty}Z^{-p^j}X(Z^{p^j}v),\quad X\in\Lambda,~v\in L^\circ.
\end{equation}
\end{proposition}

\begin{proof}
    For fixed $X\in\Lambda$ and $v\in L^\circ$, choose $j$ sufficiently large such that $Z^{p^{j-1}}v\in L$. Then
    \[
    Z^{-p^{j+1}}X(Z^{p^{j+1}}v)-Z^{-p^j}X(Z^{p^j}v)=Z^{-p^{j+1}}[X,Z^{(p-1)p^j}]Z^{(p-1)p^{j-1}}(Z^{p^{j-1}}v).   
    \]
    Using 
    \[
    [X,Z^{(p-1)p^j}]=\sum_{i=0}^{p-2}Z^{ip^j}[X,Z^{p^j}]Z^{(p-2-i)p^j},\quad\ad_{Z^{p^j}}=(\ad_{Z^p})^{p^{j-1}}
    \]
    and Lemma~\ref{lem:ad-order}, we obtain 
    \[
    [X,Z^{(p-1)p^j}]Z^{(p-1)p^{j-1}}\in\cF^{p^{j+1}+\frac{p-1}{2}p^{j-1}}\Lambda.
    \]
    This bound is uniform in $X\in\Lambda$. This estimate gives 
    \[
    \lim_{j\to\infty}\left(Z^{-p^{j+1}}X(Z^{p^{j+1}}v)-Z^{-p^j}X(Z^{p^j}v)\right)=0,
    \]
    and therefore \eqref{eq:Laurent-limit} is well defined. A similar estimate gives $X^{\circ}(v)=X(v)$ for $v\in L$.
    
    For $X_1,X_2\in\Lambda$ and $v\in L^{\circ}$, we have 
    \begin{align*}
         X_1^{\circ}\left(X_2^{\circ}(v)\right)-Z^{-p^j}X_1X_2(Z^{p^j}v)&=X_1^{\circ}\left(X_2^{\circ}(v)\right)-Z^{-p^j}X_1\left(Z^{p^j}X_2^{\circ}(v)\right)\\
         &+Z^{-p^j}X_1Z^{p^j}\left(X_2^{\circ}(v)-Z^{-p^j}X_2(Z^{p^j}v)\right),
    \end{align*}
    which tends to zero as $j$ tends to infinity. This gives $  X_1^{\circ}\left(X_2^{\circ}(v)\right)=(X_1X_2)^{\circ}(v)$.
    Moreover, by definition the extended action satisfies the relations $R_1,R_2$ in \eqref{eq:Lambda-relation}. This gives the existence of the extension. 
    
    For uniqueness, if $X^{\circ\circ}(v)$ is another extension, then 
    \[
    X^{\circ\circ}(v)-Z^{-p^j}X(Z^{p^j}v)=Z^{-p^j}[Z^{p^j},X]^{\circ\circ}(v),
    \]
    which tends to zero as $j$ tends to infinity by Lemma~\ref{lem:ad-order} and $\ad_{Z^{p^j}}=(\ad_{Z^p})^{p^{j-1}}$. This gives $X^{\circ\circ}(v)=X^{\circ}(v)$.
    \end{proof}
    
\section{The irreducible case}

Let $\pi$ be an admissible smooth representation of $D^{\times}$ over $\FF$ with central character $\zeta$. Let $M=\pi^{\vee}:=\Hom_{\FF}(\pi,\FF)$ be the Pontryagin dual of $\pi$, which is a finitely generated $\Lambda$-module with compatible $D^{\times}$-action. We make the following assumptions throughout the section:
\begin{enumerate}
\item\label{it:assum-i} There exist an integer $2\le r\le p-3$, a character $\chi$ of $H$ satisfying $\chi^\sigma=\chi\alpha^r$, and an isomorphism of graded $\gr_{\mathfrak{m}_D}(\Lambda)$-modules with compatible $H$-action 
\begin{equation}\label{eq:irr-gr}
\gr_{\mathfrak{m}_D}(M)\cong(\chi^\vee\otimes\FF[z])\oplus\left((\chi^\sigma)^\vee\otimes\FF[y]\right)\oplus\left((\chi\alpha^{-1})^\vee\otimes\FF[y]\right)\oplus\left((\chi^\sigma\alpha)^\vee\otimes\FF[z]\right),
\end{equation}
where on the right $\varpi_D$ permutes both the first two summands and the last two summands.
\item\label{it:assum-ii} There is a $D^{\times}$-equivariant isomorphism of $\Lambda$-modules
\begin{equation}\label{eq:self-duality}
\E^2_\Lambda(M)\otimes(\zeta^{-1}\circ\mathrm{Nr})\cong M,
\end{equation}
where $\E^2_\Lambda(M)=\Ext^2_\Lambda(M,\Lambda)$ is endowed with the $D^{\times}$-action defined in \cite[Proposition 3.2]{Koh17} and $\mathrm{Nr}$ is the reduced norm of $D^{\times}$. We view $\E^2_\Lambda(M)$ as a left $\Lambda$-module via the anti-involution $g\mapsto g^{-1}$ on $U_D^1/Z_D^1$.
\end{enumerate}

\subsection{An \texorpdfstring{$\cO_D^{\times}$}{O\textasciicircum{}x}-decomposition}

In this subsection, we show that the decomposition at the graded level in assumption~\ref{it:assum-i} can be lifted to the filtered level after forgetting the $\varpi_D$-action.

\begin{lemma}\label{lem:y-z-resolution}
    Let $\mu$ and $\nu$ be characters of $H$. Then $\mu\otimes\FF[z]$ and $\nu\otimes\FF[y]$ admit the following gr-free $\gr_{\mathfrak{m}_D}(\Lambda)$-module resolutions with compatible $H$-action respectively
\begin{equation}\label{eq:resolution-Fz}
\begin{split}
    0\longrightarrow\mu\alpha\otimes\gr_{\mathfrak{m}_D}(\Lambda)(-3)\xrightarrow{(-h\ \ y)}\left(\mu\alpha\otimes\gr_{\mathfrak{m}_D}(\Lambda)(-1)\right)\oplus\left(\mu\otimes\gr_{\mathfrak{m}_D}(\Lambda)(-2)\right)\\
    \xrightarrow{\binom{y}{h}}
 \mu\otimes\gr_{\mathfrak{m}_D}(\Lambda)
 \longrightarrow\mu\otimes\FF[z]\longrightarrow0,
\end{split}
\end{equation}
\begin{equation}\label{eq:resolution-Fy}
\begin{split}
    0\longrightarrow\nu\alpha^{-1}\otimes\gr_{\mathfrak{m}_D}(\Lambda)(-3)\xrightarrow{(h\ \ z)}\left(\nu\alpha^{-1}\otimes\gr_{\mathfrak{m}_D}(\Lambda)(-1)\right)\oplus\left(\nu\otimes\gr_{\mathfrak{m}_D}(\Lambda)(-2)\right)\\
    \xrightarrow{\binom{z}{-h}}
 \nu\otimes\gr_{\mathfrak{m}_D}(\Lambda)
 \longrightarrow\nu\otimes\FF[y]\longrightarrow0.
\end{split}
\end{equation}
Here we use row vectors and the matrix multiplies on the right.
\end{lemma}

\begin{proof}
This is a direct check. We take the middle term of \eqref{eq:resolution-Fz} as an example. Suppose that $(a_1,a_2)\in\gr_{\mathfrak{m}_D}(\Lambda)^2$ satisfies $a_1y+a_2h=0$. Using that $\gr_{\mathfrak{m}_D}(\Lambda)/(h)\cong\FF[y,z]$ is a domain, we obtain $a_1=ah$ for some $a\in\gr_{\mathfrak{m}_D}(\Lambda)$. Using $yh=hy$ and that $\gr_{\mathfrak{m}_D}(\Lambda)$ is a domain, we obtain $a_2=-ay$, which gives $(a_1,a_2)=-a(-h,y)$.

%When $\nu=\mu^\sigma$, we define the action of $\varpi_D$ as follow:
%\begin{enumerate}
    %\item On $\left(\mu\otimes\FF[z]\right)\oplus\left(\mu^\sigma\otimes\FF[y]\right)$, the action of $\varpi_D$ permutes $y$ and $z$.
    %\item On $\left(\mu\otimes\gr_{\mathfrak{m}_D}(\Lambda)\right)\oplus\left(\mu^\sigma\otimes\gr_{\mathfrak{m}_D}(\Lambda)\right)$, the action of $\varpi_D$ permutes $y$ (resp. $z$) in $\mu\otimes\gr_{\mathfrak{m}_D}(\Lambda)$ and $z$ (resp. $y$) in $\mu^\sigma\otimes\gr_{\mathfrak{m}_D}(\Lambda)$. Similar for $\left(\mu\alpha\otimes\gr_{\mathfrak{m}_D}(\Lambda)(-1)\right)\oplus\left(\mu^\sigma\alpha^{-1}\otimes\gr_{\mathfrak{m}_D}(\Lambda)(-1)\right)$, $\left(\mu\otimes\gr_{\mathfrak{m}_D}(\Lambda)(-2)\right)\oplus\left(\mu^\sigma\otimes\gr_{\mathfrak{m}_D}(\Lambda)(-2)\right)$ and $\left(\mu\alpha\otimes\gr_{\mathfrak{m}_D}(\Lambda)(-3)\right)\oplus\left(\mu^\sigma\alpha^{-1}\otimes\gr_{\mathfrak{m}_D}(\Lambda)(-3)\right)$.
%\end{enumerate}
%One checks that this action is compatible with everything.
\end{proof}

\begin{lemma}\label{lem:mixed-gr-Ext}
    Let $\mu$ (resp. $\nu$) be a character of $H$. Then for $i\ge0,~j\in\mathbb{Z}$, we have
    \begin{equation}\label{eq:mixed-gr-Ext-1}
         \Ext_{\gr_{\mathfrak{m}_D}(\Lambda)}^i(\mu\otimes\FF[z],\nu\otimes\FF[y])^j=
         \begin{cases}
             \nu\mu^{-1}\alpha^{-1},&(i,j)=(1,-1),~(2,-3);\\
             0,&\text{otherwise}            
         \end{cases}
    \end{equation}
and
\begin{equation}\label{eq:mixed-gr-Ext-2}
         \Ext_{\gr_{\mathfrak{m}_D}(\Lambda)}^i(\nu\otimes\FF[y],\mu\otimes\FF[z])^j=
         \begin{cases}
             \mu\nu^{-1}\alpha,&(i,j)=(1,-1),~(2,-3);\\
             0,&\text{otherwise}.            
         \end{cases}
    \end{equation}
\end{lemma}
\begin{proof}
Applying $\Hom_{\gr_{\mathfrak{m}_D}(\Lambda)}(\cdot,\nu\otimes\FF[y])$ to \eqref{eq:resolution-Fz} gives
\[
\nu\mu^{-1}\otimes\FF[y]\xrightarrow{(y\ \ 0)}
\left(\nu\mu^{-1}\alpha^{-1}\otimes\FF[y](1)\right)\oplus\left(\nu\mu^{-1}\otimes\FF[y](2)\right)
\xrightarrow{\binom{0}{y}}\nu\mu^{-1}\alpha^{-1}\otimes\FF[y](3).
\]
Taking cohomology gives \eqref{eq:mixed-gr-Ext-1}. The second formula is proved similarly.
\end{proof}

\begin{proposition}\label{prop:O_D-decomposition}
    There is a unique closed filtered $\Lambda$-module decomposition 
    \[
    M\cong M_y\oplus M_z
    \]
    such that  
    \[   
    \gr_{\mathfrak{m}_D}(M_y)\cong\left((\chi^\sigma)^\vee\otimes\FF[y]\right)\oplus\left((\chi\alpha^{-1})^\vee\otimes\FF[y]\right),\quad \gr_{\mathfrak{m}_D}(M_z)\cong(\chi^\vee\otimes\FF[z])\oplus\left((\chi^\sigma\alpha)^\vee\otimes\FF[z]\right).
    \]
    Moreover, the decomposition is compatible with the $H$-action, and the action of $\varpi_D$ permutes $M_y$ and $M_z$.
\end{proposition}
\begin{proof}
    This follows by Lemma~\ref{lem:decomposition-criterion}, assumption~\ref{it:assum-i} and Lemma~\ref{lem:mixed-gr-Ext}. The assertions for the $H$-action and $\varpi_D$-action follow by the uniqueness of the decomposition.
\end{proof}

Let $\pi_y:=M_y^{\vee},~\pi_z:=M_z^{\vee}$. Then $\pi_y$ and $\pi_z$ are smooth admissible representations of $\cO_D^{\times}Z_D$ over $\FF$ with central character $\zeta$.

\begin{corollary}\label{cor:pi-induction}
    $\pi\cong\Ind_{\cO_D^{\times}Z_D}^{D^{\times}}\pi_z$.
\end{corollary}

\begin{corollary}\label{cor:Ext-decomposition}
    For $i\ge0$, we have 
    \[
    \Ext^i_{D^\times,\zeta}(\pi,\pi)\cong\Ext^i_{\Lambda}(M_z,M_z)^H\oplus\Ext^i_{\Lambda}(M_z,M_y)^H.
    \]
    In particular, for $i\ge3$, we have $\Ext^i_{D^\times,\zeta}(\pi,\pi)=0$.
\end{corollary}

\begin{proof}
    By Corollary~\ref{cor:pi-induction} and Shapiro's lemma, we have 
    \[
    \Ext^i_{D^\times,\zeta}(\pi,\pi)\cong\Ext^i_{\cO_D^{\times}Z_D,\zeta}(\pi|_{\cO_D^{\times}Z_D},\pi_z).
    \]
    Since the extension group is computed in the category of smooth representations with the fixed central character $\zeta$, we have 
    \[
    \Ext^i_{\cO_D^{\times}Z_D,\zeta}(\pi|_{\cO_D^{\times}Z_D},\pi_z)\cong\Ext^i_{(U_D^1/Z_D^1)\rtimes H}(\pi_y\oplus\pi_z,\pi_z).
    \]
    Since the order of $H$ is prime to $p$, we have 
    \[
    \Ext^i_{(U_D^1/Z_D^1)\rtimes H}(\pi_y\oplus\pi_z,\pi_z)\cong\Ext^i_{U_D^1/Z_D^1}(\pi_y,\pi_z)^H\oplus\Ext^i_{U_D^1/Z_D^1}(\pi_z,\pi_z)^H.
    \]
    Pontryagin duality gives the first assertion. By Lemma~\ref{lem:resolution-lifting}, Lemma~\ref{lem:y-z-resolution} and Proposition~\ref{prop:O_D-decomposition}, $M_z$ admits a filt-free $\Lambda$-module resolution of length 2. This gives the second assertion.
\end{proof}

\begin{proposition}\label{prop:mixed-ext-dim}
    We have $\Hom_{\Lambda}(M_z,M_y)^H=0$ and 
    \[
    \dim_\FF \Ext^1_{\Lambda}(M_z,M_y)^H=\dim_\FF \Ext^2_{\Lambda}(M_z,M_y)^H=2.
    \]
\end{proposition}
\begin{proof}
    By Lemma~\ref{lem:resolution-lifting}, Remark~\ref{rmk:H-compatible}, Lemma~\ref{lem:y-z-resolution} and Proposition~\ref{prop:O_D-decomposition}, $M_z$ admits a filt-free $\Lambda$-module resolution $P_{\bullet}$ with compatible $H$-action, whose terms are of finite rank. Then $\Hom_{\Lambda}(P_{\bullet},M_y)^H$ is a cochain complex of filtered $\FF$-vector spaces whose filtrations are complete and bounded below with finite-dimensional graded pieces. By Proposition~\ref{prop:spectral-sequence}, we have a convergent spectral sequence $(E_\bullet^{i,j})\Longrightarrow\Ext^{i+j}_{\Lambda}(M_z,M_y)^H$ with 
    \[
    E_1^{i,j}\cong\Ext_{\gr_{\mathfrak{m}_D}(\Lambda)}^{i+j}\left(\gr_{\mathfrak{m}_D}(M_z),\gr_{\mathfrak{m}_D}(M_y)\right)^{i,H}.
    \]
    By Lemma \ref{lem:mixed-gr-Ext} and Proposition~\ref{prop:O_D-decomposition}, we have
    \[
    \Ext_{\gr_{\mathfrak{m}_D}(\Lambda)}^i\left(\gr_{\mathfrak{m}_D}(M_z),\gr_{\mathfrak{m}_D}(M_y)\right)^{j,H}=
    \begin{cases}
         \FF^2,&(i,j)=(1,-1),~(2,-3);\\
             0,&\text{otherwise}.     
    \end{cases}
    \]
    By analyzing the differentials of the spectral sequence, we have $E_{\infty}^{i,j}\cong E_1^{i,j}$ for all $i,j$, which gives the result.
\end{proof}

\subsection{Leading coefficients and self-duality}

Lift the $\FF[z]$-basis of $\gr_{\mathfrak{m}_D}(M_z)$ to $e_0,e_1\in M_z$ with $H$-characters $\chi^\vee$ and $(\chi^\sigma\alpha)^\vee$ respectively. Then 
\[
M_z\cong\FF[[Z]]e_0\oplus\FF[[Z]]e_1
\]
as an $\FF[[Z]]$-module and the $Z$-adic filtration coincides with the $\mathfrak{m}_D$-adic filtration. If $(e_0',e_1')$ is another $\FF[[Z]]$-basis of $M_z$, with $H$-characters $\chi^\vee$ and $(\chi^\sigma\alpha)^\vee$ respectively, then we must have
\begin{equation}\label{eq:base-change}
     e_0'=A(Z^{p+1})e_0+Z^{p-r}B(Z^{p+1})e_1,\quad e_1'=Z^{r+1}C(Z^{p+1})e_0+D(Z^{p+1})e_1
\end{equation}
for some $A,D\in\FF[[Z]]^{\times}$ and $B,C\in\FF[[Z]]$. Therefore, the $Z$-adic filtration on $\FF[[Z]]e_0'\oplus\FF[[Z]]e_1'$ automatically coincides with the $\mathfrak{m}_D$-adic filtration.

Since $y$ kills $\gr_{\mathfrak{m}_D}(M_z)$, the action of $Y$ on $M_z$ raises the $\mathfrak{m}_D$-adic filtration by at least two. By comparing $H$-actions, we have
\begin{equation}\label{eq:leading-coefficient}
Y(Z^ke_0)\equiv h_0(k)Z^{k+p-1-r}e_1\pmod{Z^{k+p}M_z},\quad Y(Z^ke_1)\equiv h_1(k)Z^{k+r}e_0\pmod{Z^{k+p}M_z}
\end{equation}
for every integer $k\ge0$.

Put $\wet(e_0)=0,~\wet(e_1)=r-\frac{p-1}{2}$ and define
\[
 \cF^nM_z=Z^{\max\{n,0\}}\FF[[Z]]e_0\oplus Z^{\max\{n-\wet(e_1),0\}}\FF[[Z]]e_1.
\]
By formula \eqref{eq:base-change}, the filtration $\cF^nM_z$ is independent of the choice of $e_0$ and $e_1$. By \eqref{eq:leading-coefficient}, for any $n$ we have
\[
Y(\cF^nM_z)\subseteq\cF^{n+(p-1)/2}M_z.
\]
Therefore $(M_z,\cF)$ is a filtered $(\Lambda,\cF)$-module. Moreover, $\cF^nM_z$ is a good filtration in the sense of \cite[\S I.5]{LvO96}.

Write 
\[
\widetilde{M}_z=\gr_{\cF}(M_z),\quad\widetilde{e}_0=\gr_{\cF}(e_0),\quad\widetilde{e}_1=\gr_{\cF}(e_1).
\]
By formula \eqref{eq:leading-coefficient}, we have
\[
\widetilde{y}(\widetilde{z}^k\widetilde{e}_0)=h_0(k)\widetilde{z}^{k+p-1-r}\widetilde{e}_1,\quad\widetilde{y}(\widetilde{z}^k\widetilde{e}_1)=h_1(k)\widetilde{z}^{k+r}\widetilde{e}_0.
\]

\begin{lemma}\label{lem:leading-linear}
    There are $b_0,b_1\in\FF^\times$ and $\beta_0,\beta_1\in\FF$, satisfying $b_0b_1=2c,~\beta_0-\beta_1=r-\frac{p-1}{2},$ such that $h_0(k)=b_0(k-\beta_0),~h_1(k)=b_1(k-\beta_1)$ for every integer $k\ge0$. Moreover, $\beta_0,\beta_1$ are independent of the choice of $e_0,e_1$.
\end{lemma}

\begin{proof}
    By relation $R_2$ in \eqref{eq:Lambda-relation} and Lemma~\ref{lem:g-weight}, for every integer $k\ge0$, the coefficient of $Z^{k+p+1-r}e_1$ (resp. $Z^{k+r+2}e_0$) in $[Z,[Y,Z]]Z^ke_0$ (resp. $[Z,[Y,Z]]Z^ke_1$) is zero, which gives
    \[
    h_0(k+2)-2h_0(k+1)+h_0(k)=0,\quad h_1(k+2)-2h_1(k+1)+h_1(k)=0.
    \]
    Thus there exist $b_0,b_1,d_0,d_1\in\FF$, such that $h_0(k)=b_0k+d_0,~h_1(k)=b_1k+d_1$ for every integer $k\ge0$. By relation $R_1$ in \eqref{eq:Lambda-relation} and Lemma~\ref{lem:g-weight}, for every integer $k\ge0$, the coefficient of $Z^{k+p}e_0$ (resp. $Z^{k+p}e_1$) in $([Y,[Y,Z]]-cZ^p)Z^ke_0$ (resp. $([Y,[Y,Z]]-cZ^p)Z^ke_1$) is zero, which gives
    \[
    b_0b_1(p-r)+b_0d_1-b_1d_0=c,\quad b_0b_1(r+1)+b_1d_0-b_0d_1=c.
    \]
    Their sum is $b_0b_1=2c$, which is nonzero by Lemma~\ref{lem:c-nonzero}. Take $\beta_0=-d_0/b_0,~\beta_1=-d_1/b_1$. Then the difference of the equalities above gives $ \beta_0-\beta_1=r-\frac{p-1}{2}$.

    Suppose that $(e_0',e_1')$ is another $H$-eigenbasis and write $\widetilde{e}_0'=\gr_{\cF}(e_0'),~\widetilde{e}_1'=\gr_{\cF}(e_1')$. By \eqref{eq:base-change}, we have 
    \[
    \widetilde{e}_0'=A(0)\widetilde{e}_0,\quad\widetilde{e}_1'=D(0)\widetilde{e}_1.
    \]
    Therefore, for every integer $k\ge0$, we have 
    \[
    \widetilde{y}(\widetilde{z}^k\widetilde{e}_0')=\frac{A(0)}{D(0)}b_0(k-\beta_0)\widetilde{z}^{k+p-1-r}\widetilde{e}_1',\quad\widetilde{y}(\widetilde{z}^k\widetilde{e}_1')=\frac{D(0)}{A(0)}b_1(k-\beta_1)\widetilde{z}^{k+r}\widetilde{e}_0'.
    \]
    This shows that $\beta_0,\beta_1$ are independent of the choice of $e_0,e_1$.
\end{proof}

Next we use our self-duality assumption~\ref{it:assum-ii} to determine the leading coefficients $\beta_0$ and $\beta_1$.

\begin{lemma}\label{lem:grM_z-dual}
    The graded $\widetilde{\Lambda}$-module $\widetilde{M_z}$ is Cohen-Macaulay of grade $2$. Moreover, $\E^2_{\widetilde{\Lambda}}(\widetilde{M_z})$
    is a graded right $\widetilde{\Lambda}$-module, free of rank two over $\FF[\widetilde z]$, with homogeneous generators $\widetilde{e_0}^{\d},\widetilde{e_1}^{\d}$ of degrees
    $-p,-r-\frac{p+1}{2}$ and $H$-characters $\chi\alpha^{-1},\chi^\sigma$ respectively. The action of $\widetilde{\Lambda}$ is given by
    \[
    \widetilde{e_0}^{\d}\widetilde{z}^k\widetilde{h}=b_1\widetilde{e_1}^{\d}\widetilde{z}^{k+r+1},\quad \widetilde{e_0}^{\d}\widetilde{z}^k\widetilde{y}=-b_1(k+\beta_1+r+1)\widetilde{e_1}^{\d}\widetilde{z}^{k+r},
    \]
    \[
    \widetilde{e_1}^{\d}\widetilde{z}^k\widetilde{h}=b_0\widetilde{e_0}^{\d}\widetilde{z}^{k+p-r},\quad \widetilde{e_1}^{\d}\widetilde{z}^k\widetilde{y}=-b_0(k+\beta_0+p-r)\widetilde{e_0}^{\d}\widetilde{z}^{k+p-1-r}.
    \]
\end{lemma}

\begin{proof}
    Put $R=\FF[\widetilde{z},\widetilde{h}]$. Consider the $R$-linear surjection $\varepsilon:R^2\to\widetilde{M_z}$ sending the basis of $R^2$ to $\widetilde{e_0},\widetilde{e_1}$ respectively. By Lemma~\ref{lem:leading-linear}, we have 
    \[
    (\widetilde{h},-b_0\widetilde{z}^{p-r}),(-b_1\widetilde{z}^{r+1},\widetilde{h})\in\ker\varepsilon.
    \]
    Put 
    \[
    K=R(\widetilde{h},-b_0\widetilde{z}^{p-r})+R(-b_1\widetilde{z}^{r+1},\widetilde{h})\subseteq R^2.
    \]
    Then $\varepsilon$ induces $\overline{\varepsilon}:R^2/K\to\widetilde{M_z}$. Let $\iota:\widetilde{M_z}\to R^2/K$ be the $\FF[\widetilde{z}]$-linear map sending $\widetilde{e_0},\widetilde{e_1}$ to the image of the basis of $R^2$. Then $\iota$ is surjective and $\overline{\varepsilon}\circ\iota=\id$, which implies $\overline{\varepsilon}$ is an isomorphism. Thus, we have a gr-free $R$-module resolution of $\widetilde{M_z}$ with compatible $H$-action
    \begin{equation}\label{eq:R-resolution-Mz}
        \begin{split}
             0\longrightarrow\left(\chi^\vee\otimes R(-\frac{p+1}{2})\right)\oplus\left((\chi^\sigma\alpha)^\vee\otimes R(-r-1)\right)\xrightarrow{\left(\begin{smallmatrix} \widetilde{h} & -b_0\widetilde{z}^{p-r} \\ -b_1\widetilde{z}^{r+1} & \widetilde{h} \end{smallmatrix}\right)}\\
             (\chi^\vee\otimes R)\oplus\left((\chi^\sigma\alpha)^
            \vee\otimes R(\frac{p-1}{2}-r)\right)\xrightarrow{\,\varepsilon\,}\widetilde{M_z}\longrightarrow0.
        \end{split}
    \end{equation}

    Applying $\Hom_R(\cdot,R)$ to \eqref{eq:R-resolution-Mz}, we obtain $\E_R^i(\widetilde{M_z})=0$ for $i\ne1$ and a graded right $R$-module isomorphism with compatible $H$-action 
    \[
    \E_R^1(\widetilde{M_z})\cong q_0\FF[\widetilde{z}]\oplus q_1\FF[\widetilde{z}],
    \]
    where 
    \[
    \deg(q_0)=-\frac{p+1}{2},\quad\deg(q_1)=-r-1,\quad q_0\widetilde{h}=b_1q_1\widetilde{z}^{r+1},\quad q_1\widetilde{h}=b_0q_0\widetilde{z}^{p-r}
    \]
    and $H$ acts on $q_0$ (resp. $q_1$) by $\chi$ (resp. $\chi\alpha^{r+1}$). Since $\widetilde{\Lambda}=\bigoplus_{k\ge0}R\widetilde{y}^k$ we have $\Ext_R^i(\widetilde{M_z},\widetilde{\Lambda})=0$ for $i\ne1$ and 
    \[
    \Ext_R^1(\widetilde{M_z},\widetilde{\Lambda})\cong\bigoplus_{k\ge0}(q_0\FF[\widetilde{z}]\oplus q_1\FF[\widetilde{z}])\widetilde{y}^k
    \]
    as a graded right $\widetilde{\Lambda}$-module with compatible $H$-action.

    Consider the $\widetilde{\Lambda}$-linear morphism 
    \[
    \theta:\widetilde{\Lambda}\otimes_R\widetilde{M_z}\longrightarrow\widetilde{\Lambda}\otimes_R\widetilde{M_z},~a\otimes v\mapsto a\widetilde{y}\otimes v-a\otimes\widetilde{y}v.
    \]
    This is well-defined by the definition of the skew polynomial ring $\widetilde{\Lambda}$. Since $\widetilde{\Lambda}=\bigoplus_{k\ge0}\widetilde{y}^kR$ is a free right $R$-module with basis $\{\widetilde{y}^k:k\ge0\}$, the morphism $\theta$ is injective
    by considering the highest degree of $\widetilde{y}$, and its cokernel is $\widetilde{M_z}$ by successively eliminating positive powers of $\widetilde{y}$. Therefore, we have a short exact sequence of $\widetilde{\Lambda}$-modules 
    \begin{equation}\label{eq:standard-ses}
    0\longrightarrow \widetilde{\Lambda}\otimes_R\widetilde{M_z}
    \xrightarrow{\,\theta\,}\widetilde{\Lambda}\otimes_R\widetilde{M_z}
    \longrightarrow\widetilde{M_z}\longrightarrow0.
    \end{equation}
    The long exact sequence associated with \eqref{eq:standard-ses}, the adjunction formula 
    \[
    \Ext_{\widetilde{\Lambda}}^i(\widetilde{\Lambda}\otimes_R\widetilde{M_z},\widetilde{\Lambda})
    \cong\Ext_R^i(\widetilde{M_z},\widetilde{\Lambda})
    \]
    and the vanishing of $\Ext_R^i(\widetilde{M_z},\widetilde{\Lambda})$ for $i\ne1$ give $\E^j_{\widetilde{\Lambda}}(\widetilde{M_z})=0$ for $j\ne1,2$ and
    \[
    0\longrightarrow\E^1_{\widetilde{\Lambda}}(\widetilde{M_z})\longrightarrow\Ext_R^1(\widetilde{M_z},\widetilde{\Lambda})
    \xrightarrow{\,\Theta\,}\Ext_R^1(\widetilde{M_z},\widetilde{\Lambda})\longrightarrow\E^2_{\widetilde{\Lambda}}(\widetilde{M_z})\longrightarrow0,
    \]
    where $\Theta$ is induced by $\theta$ and the adjunction formula.

    We compute $\Theta$ explicitly. Applying $\widetilde{\Lambda}\otimes_R(\cdot)$ to the resolution \eqref{eq:R-resolution-Mz} and extending $\theta$ to the resolution, we obtain a $\widetilde{\Lambda}$-linear chain morphism (omitting degree shifts and $H$-actions)
    \[
    \begin{tikzcd}
	0 & {\widetilde{\Lambda}^2} & {\widetilde{\Lambda}^2} & {\widetilde{\Lambda}\otimes_R \widetilde{M_z}} & 0 \\
	0 & {\widetilde{\Lambda}^2} & {\widetilde{\Lambda}^2} & {\widetilde{\Lambda}\otimes_R \widetilde{M_z}} & 0
	\arrow[from=1-1, to=1-2]
	\arrow[from=1-2, to=1-3]
	\arrow["{\theta_1}"', from=1-2, to=2-2]
	\arrow[from=1-3, to=1-4]
	\arrow["{\theta_0}"', from=1-3, to=2-3]
	\arrow[from=1-4, to=1-5]
	\arrow["\theta"', from=1-4, to=2-4]
	\arrow[from=2-1, to=2-2]
	\arrow[from=2-2, to=2-3]
	\arrow[from=2-3, to=2-4]
	\arrow[from=2-4, to=2-5]
    \end{tikzcd}
    \]
    where 
    \[
    \theta_0=
    \begin{pmatrix}
    \widetilde{y} & b_0\beta_0\widetilde{z}^{p-1-r}\\
    b_1\beta_1\widetilde{z}^r & \widetilde{y}
    \end{pmatrix},
    \quad
    \theta_1=
    \begin{pmatrix}
    \widetilde{y} & b_0(\beta_0+p-r)\widetilde{z}^{p-1-r}\\
    b_1(\beta_1+r+1)\widetilde{z}^r & \widetilde{y}
    \end{pmatrix}
    \]
    (here we use the equalities $b_0b_1=2c,~\beta_0-\beta_1=r-\frac{p-1}{2}$). Applying $\Hom_{\widetilde{\Lambda}}(\cdot,\widetilde{\Lambda})$ and the adjunction formula and taking cohomology, we obtain a morphism of right $\widetilde{\Lambda}$-modules
    \[
    \Theta:\Ext_R^1(\widetilde{M_z},\widetilde{\Lambda})\to\alpha^{-1}\otimes\Ext_R^1(\widetilde{M_z},\widetilde{\Lambda})(\frac{p-1}{2})
    \]
    sending $q_0$ to $q_0\widetilde{y}+b_1(\beta_1+r+1)q_1\widetilde{z}^r$ and $q_1$ to $b_0(\beta_0+p-r)q_0\widetilde{z}^{p-1-r}+q_1\widetilde{y}$. Here the twist and shift come from $\theta$. By considering the highest degree of $\widetilde{y}$, we obtain that $\Theta$ is injective, which implies $\E^1_{\widetilde{\Lambda}}(\widetilde{M_z})=0$. Moreover, we have 
    \[
    \E^2_{\widetilde{\Lambda}}(\widetilde{M_z})\cong\Coker\Theta=\widetilde{e_0}^{\d}\FF[\widetilde{z}]\oplus\widetilde{e_1}^{\d}\FF[\widetilde{z}],
    \]
    where 
    \[
    \deg(\widetilde{e_0}^{\d})=-p,\quad\deg(\widetilde{e_1}^{\d})=-r-\frac{p+1}{2},\quad\widetilde{e_0}^{\d}\widetilde{h}=b_1\widetilde{e_1}^{\d}\widetilde{z}^{r+1},\quad\widetilde{e_1}^{\d}\widetilde{h}=b_0\widetilde{e_0}^{\d}\widetilde{z}^{p-r},
    \]
    \[
    \widetilde{e_0}^{\d}\widetilde{y}=-b_1(\beta_1+r+1)\widetilde{e_1}^{\d}\widetilde{z}^r,\quad\widetilde{e_1}^{\d}\widetilde{y}=-b_0(\beta_0+p-r)\widetilde{e_0}^{\d}\widetilde{z}^{p-1-r}
    \]
    and $H$ acts on $\widetilde{e_0}^{\d}$ (resp. $\widetilde{e_1}^{\d}$) by $\chi\alpha^{-1}$ (resp. $\chi\alpha^r=\chi^\sigma$).
\end{proof}

Recall that for a right $\Lambda$-module, we can view it as a left $\Lambda$-module via the anti-involution $g\to g^{-1}$ on $U_D^1/Z_D^1$. The anti-involution sends $Y$ to $-Y$ and $Z$ to $-Z$, and hence preserves both the $\mathfrak{m}_D$-adic filtration and the weighted filtration. By assumption~\ref{it:assum-ii}, we have an isomorphism of $\Lambda$-modules with compatible $H$-action 
\[
\left(\E_\Lambda^2(M_y)\otimes(\zeta^{-1}\circ\mathrm{Nr})\right)\oplus\left(\E_\Lambda^2(M_z)\otimes(\zeta^{-1}\circ\mathrm{Nr})\right)\cong M,
\]
which is automatically $\mathfrak{m}_D$-adic filtered.

\begin{lemma}\label{lem:z-dual-z}
    The isomorphism \eqref{eq:self-duality} identifies $M_z$ with $\E_\Lambda^2(M_z)\otimes(\zeta^{-1}\circ\mathrm{Nr})$.
\end{lemma}

\begin{proof}
    Choose a strict $\mathfrak{m}_D$-adic filt-free $\Lambda$-module resolution of $M_z$ with compatible $H$-action. By \cite[Lemma III.2.2.3(1)]{LvO96}, this induces a good filtration $G^n$ on $\E_\Lambda^2(M_z)$, which hence is complete and separated. %By Proposition~\ref{prop:spectral-sequence}, we have a converged spectral sequence $(E_r^{p,q})\Longrightarrow G^p\left(\Ext^{p+q}_{\Lambda}(M_z,\Lambda)\right)$ with $E_1^{p,q}\cong\Ext_{\gr_{\mathfrak{m}_D}(\Lambda)}^{p+q}\left(\gr_{\mathfrak{m}_D}(M_z),\gr_{\mathfrak{m}_D}(\Lambda)\right)^p$. By the resolution \eqref{eq:resolution-Fz} and Proposition~\ref{prop:O_D-decomposition}, we have $\Ext^i_{\gr_{\mathfrak{m}_D}(\Lambda)}\left(\gr_{\mathfrak{m}_D}(M_z),\gr_{\mathfrak{m}_D}(\Lambda)\right)=0$ for $i\ne2$ and a graded $\gr_{\mathfrak{m}_D}(\Lambda)$-module isomorphism
    Using the resolution \eqref{eq:resolution-Fz} and Proposition~\ref{prop:O_D-decomposition}, one checks that $\gr_{\mathfrak{m}_D}(M_z)$ is Cohen-Macaulay of grade $2$ and
    \[
    \E^2_{\gr_{\mathfrak{m}_D}(\Lambda)}\left(\gr_{\mathfrak{m}_D}(M_z)\right)\cong\gr_{\mathfrak{m}_D}(M_z)(3)\otimes(\zeta\circ\mathrm{Nr}|_H)
    \]
    (here we use $\zeta\circ\mathrm{Nr}|_H=\chi\chi^{\sigma}$). This can also be obtained by \cite[Corollary 4.14]{Zhao}, which is based on \cite[Lemma 3.3.1.9, Proposition 3.3.1.10]{BHHMS2}. By \cite[Proposition 3.3.4.6]{BHHMS2}, we have
    \[
    \gr_G\left(\E_\Lambda^2(M_z)\right)\otimes(\zeta^{-1}\circ\mathrm{Nr}|_H)\cong\E^2_{\gr_{\mathfrak{m}_D}(\Lambda)}\left(\gr_{\mathfrak{m}_D}(M_z)\right)\otimes(\zeta^{-1}\circ\mathrm{Nr}|_H)\cong\gr_{\mathfrak{m}_D}(M_z)(3)
    \]
    as a graded $\gr_{\mathfrak{m}_D}(\Lambda)$-module with compatible $H$-action. Lifting the $H$-eigenbasis, we obtain 
    \[
    \E_\Lambda^2(M_z)\otimes(\zeta^{-1}\circ\mathrm{Nr}|_H)\cong\FF[[Z]]q_0\oplus\FF[[Z]]q_1
    \]
    with 
    \[
    G^n\left(\E_\Lambda^2(M_z)\right)\otimes(\zeta^{-1}\circ\mathrm{Nr}|_H)=Z^{\max\{n+3,0\}}\FF[[Z]]q_0\oplus Z^{\max\{n+3,0\}}\FF[[Z]]q_1,
    \]
    where $H$ acts on $q_0$ (resp. $q_1$) by $\chi^\vee$ (resp. $(\chi^\sigma\alpha)^\vee$). Since $y$ kills $\gr_G\left(\E_\Lambda^2(M_z)\right)$, we have 
    \[
    YG^n\left(\E_\Lambda^2(M_z)\right)\subseteq G^{n+2}\left(\E_\Lambda^2(M_z)\right),
    \]
    and therefore 
    \[
    \mathfrak{m}_DG^n\left(\E_\Lambda^2(M_z)\right)=G^{n+1}\left(\E_\Lambda^2(M_z)\right)
    \]
    for any $n\ge-3$. We obtain that the $\mathfrak{m}_D$-adic filtration on $\E_\Lambda^2(M_z)$ is precisely the filtration $G$ shifted by $-3$. Consequently, we have 
    \[
    \gr_{\mathfrak{m}_D}\left(\E_\Lambda^2(M_z)\right)\otimes(\zeta^{-1}\circ\mathrm{Nr}|_H)\cong\gr_{\mathfrak{m}_D}(M_z)
    \]
    as a graded $\gr_{\mathfrak{m}_D}(\Lambda)$-module with compatible $H$-action. The $\varpi_D$-action gives 
    \[
    \gr_{\mathfrak{m}_D}\left(\E_\Lambda^2(M_y)\right)\otimes(\zeta^{-1}\circ\mathrm{Nr}|_H)\cong\gr_{\mathfrak{m}_D}(M_y).    
    \]
    The result follows by the uniqueness of the decomposition in Proposition~\ref{prop:O_D-decomposition}.
\end{proof}

\begin{proposition}\label{prop:leading-coefficient}
    $\beta_0=\frac{r}{2},~\beta_1=-\frac{r+1}{2}$.
\end{proposition}

\begin{proof}
    Choose a strict filt-free $(\Lambda,\cF)$-module resolution of $M_z$, which induces a complete and separated filtration $\cG^n$ on $\E_\Lambda^2(M_z)$. By Lemma~\ref{lem:grM_z-dual} and \cite[Proposition 3.3.4.6]{BHHMS2}, we have a graded $\widetilde{\Lambda}$-module isomorphism with compatible $H$-action
    \[
    \gr_\cG\left(\E_\Lambda^2(M_z)\right)\otimes(\zeta^{-1}\circ\mathrm{Nr}|_H)\cong\E^2_{\widetilde{\Lambda}}(\widetilde{M_z})\otimes(\zeta^{-1}\circ\mathrm{Nr}|_H)\cong(\chi\chi^\sigma)^{-1}\otimes(\FF[\widetilde{z}]\widetilde{e_1}^{\d}\oplus\FF[\widetilde{z}]\widetilde{e_0}^{\d}).
    \]
    Here we use the involution to view $\E^2_{\widetilde{\Lambda}}(\widetilde{M_z})$ as a left $\widetilde{\Lambda}$-module, and the action is given by 
    \[
    \widetilde{y}\widetilde{z}^k\widetilde{e_1}^{\d}=(-1)^{p-1-r}b_0(k+\beta_0+p-r)\widetilde{z}^{k+p-1-r}\widetilde{e_0}^{\d},\quad\widetilde{y}\widetilde{z}^k\widetilde{e_0}^{\d}=(-1)^rb_1(k+\beta_1+r+1)\widetilde{z}^{k+r}\widetilde{e_1}^{\d}.
    \]
    Lifting the $H$-eigenbasis, we obtain 
    \[
    \E^2_{\Lambda}(M_z)\otimes(\zeta^{-1}\circ\mathrm{Nr}|_H)\cong\FF[[Z]]e_1^{\d}\oplus\FF[[Z]]e_0^{\d},
    \]
    where
    \[
    \cG^n\left(\E_\Lambda^2(M_z)\right)\otimes(\zeta^{-1}\circ\mathrm{Nr}|_H)=Z^{\max\{n+r+\frac{p+1}{2},0\}}\FF[[Z]]e_1^{\d}\oplus Z^{\max\{n+p,0\}}\FF[[Z]]e_0^{\d}
    \]
    and $H$ acts on $e_1^{\d}$ (resp. $e_0^{\d}$) by $\chi^\vee$ (resp. $(\chi^\sigma\alpha)^\vee$). By the discussion at the beginning of this section, the filtration $\cF$ on $\E_\Lambda^2(M_z)$ is precisely the filtration $\cG$ shifted by $-r-\frac{p+1}{2}$. Consequently, we have 
    \[
    \gr_{\cF}\left(\E_\Lambda^2(M_z)\right)\otimes(\zeta^{-1}\circ\mathrm{Nr}|_H)\cong\gr_{\cF}(M_z)
    \]
    as a graded $\widetilde{\Lambda}$-module with compatible $H$-action, with the leading coefficient $\beta_0^{\d}=-\beta_0+r$ and $\beta_1^{\d}=-\beta_1-r-1$. By Lemma~\ref{lem:leading-linear}, we have $\beta_0^{\d}=\beta_0$ and $\beta_1^{\d}=\beta_1$, which give the result.
\end{proof}

\subsection{Balancing lattice and a further decomposition}

In this subsection, we define a $\Lambda$-submodule $M_z'$ of $M_z$ with finite codimension and decompose it into a direct sum of two $\Lambda$-submodules free of rank-$1$ over $\FF[[Z]]$, using the weighted filtration introduced in the previous subsection.

Define
\[
M_z':=\cF^{\max\{r-\frac{p-1}{2},0\}}M_z=
\begin{cases}
 Z^{r-\frac{p-1}{2}}\FF[[Z]]e_0\oplus\FF[[Z]]e_1,
       &r\ge(p-1)/2,\\
 \FF[[Z]]e_0\oplus Z^{(p-1)/2-r}\FF[[Z]]e_1,
       &r<(p-1)/2.
 \end{cases}
\]
This is a closed $\Lambda$-submodule with compatible $H$-action. Let 
\[
f_0=Z^{\max\{r-\frac{p-1}{2},0\}}e_0,~f_1=Z^{\max\{\frac{p-1}{2}-r,0\}}e_1
\]
be the two generators. Endow $M_z'=\FF[[Z]]f_0\oplus\FF[[Z]]f_1$ with the $Z$-adic filtration. This is the filtration induced by $\cF^nM_z$ up to a shift, and we still denote it by $\cF^nM_z'$. Put 
\[
\vartheta=\beta_1+\min\left\{r-\frac{p-1}{2},0\right\}.
\]
Then by \eqref{eq:leading-coefficient} and Lemma~\ref{lem:leading-linear}, the action of $\Lambda$ on $M_z'$ is given by
\[
Y(Z^kf_0)\equiv b_0(k-\vartheta)Z^{k+(p-1)/2}f_1\pmod{\cF^{k+p}M_z'},
\]
\[Y(Z^kf_1)\equiv b_1(k-\vartheta)Z^{k+(p-1)/2}f_0\pmod{\cF^{k+p}M_z'}.
\]
This implies $M_z'$ is a filtered $(\Lambda,\cF)$-module, and the action of $\widetilde{\Lambda}$ on 
\[
\widetilde{M_z'}:=\gr_{\cF}(M_z')=\FF[\widetilde{z}]\widetilde{f_0}\oplus\FF[\widetilde{z}]\widetilde{f_1}
\]
is given by
\[
\widetilde{y}(\widetilde{z}^k\widetilde{f_0})=b_0(k-\vartheta)\widetilde{z}^{k+(p-1)/2}\widetilde{f_1},\quad\widetilde{y}(\widetilde{z}^k\widetilde{f_1})=b_1(k-\vartheta)\widetilde{z}^{k+(p-1)/2}\widetilde{f_0},
\]
where $\widetilde{f}_0:=\gr_{\cF}(f_0),~\widetilde{f}_1:=\gr_{\cF}(f_1)$.

Define 
\[
H'=\ker\alpha^{(p+1)/2}\subseteq H.
\]
It has index two in $H$, contains $H\cap Z_D$, and
$\alpha|_{H'}$ has order $(p+1)/2$. Then $H'$ acts on $f_0$ and $f_1$ via the same character, say $\psi$. Choose $b\in\FF^{\times}$ such that $b^2=b_0b_1=2c$. Put 
\[
e_+:=f_0+\frac{b_0}{b}f_1,\quad e_-:=f_0-\frac{b_0}{b}f_1.
\]
Then the action of $\widetilde{\Lambda}$ on $\widetilde{M_z'}=\FF[\widetilde{z}]\widetilde{e_+}\oplus\FF[\widetilde{z}]\widetilde{e_-}$ is given by
\begin{equation}\label{eq:leading-coefficient-N}
    \widetilde{y}(\widetilde{z}^k\widetilde{e_+})=b(k-\vartheta)\widetilde{z}^{k+(p-1)/2}\widetilde{e_+},\quad\widetilde{y}(\widetilde{z}^k\widetilde{e_-})=-b(k-\vartheta)\widetilde{z}^{k+(p-1)/2}\widetilde{e_-}.
\end{equation}
Therefore, we have a graded $\widetilde{\Lambda}$-module decomposition 
\[
\widetilde{M_z'}=\widetilde{N_+}\oplus\widetilde{N_-}
\]
with compatible $H'$-action, where $\widetilde{N_+}:=\FF[\widetilde{z}]\widetilde{e_+}$ and $\widetilde{N_-}:=\FF[\widetilde{z}]\widetilde{e_-}$. Moreover, elements in $H\setminus H'$ permute $\widetilde{N_+}$ and $\widetilde{N_-}$.

\begin{lemma}\label{lem:N+-rsolution}
    $\widetilde{N_+}$ admits a gr-free $\widetilde{\Lambda}$-module resolution with compatible $H'$-action
    \begin{equation}\label{eq:N+-resolution}
        0\longrightarrow\psi\alpha\otimes\widetilde{\Lambda}(-p)\xrightarrow{\,d_2\,}\left(\psi\alpha\otimes\widetilde{\Lambda}(-\frac{p-1}{2})\right)\oplus\left(\psi\otimes\widetilde{\Lambda}(-\frac{p+1}{2})\right)\xrightarrow{d_1}\psi\otimes\widetilde{\Lambda}\xrightarrow{\,\varepsilon\,}\widetilde{N_+}\longrightarrow0,
    \end{equation}
    where $\varepsilon$ sends the basis of $\widetilde{\Lambda}$ to $\widetilde{e_+}$, 
    \[
    d_1=\left(\widetilde{y}+b\vartheta\widetilde{z}^{(p-1)/2},\widetilde{h}-b\widetilde{z}^{(p+1)/2}\right)^\top,\quad d_2=(\widetilde{h}-b\widetilde{z}^{(p+1)/2},-\widetilde{y}-b\vartheta\widetilde{z}^{(p-1)/2}-b\frac{p+1}{2}\widetilde{z}^{(p-1)/2}).
    \]
    Moreover, the resolution of $\widetilde{N_-}$ is obtained by replacing $b$ by $-b$.
\end{lemma}

\begin{proof}
    One can check that this is a complex. By Proposition~\ref{prop:graded-algebra} $\widetilde{\Lambda}$ is a domain, and thus $d_2$ is injective. Arguments similar to the proof of Lemma~\ref{lem:grM_z-dual} imply that $\varepsilon$ induces an isomorphism $\widetilde{\Lambda}/\Im d_1\cong\widetilde{
    N_+}$. It suffices to prove $\ker d_1=\Im d_2$.

    Put $u=\widetilde{y}+b\vartheta\widetilde{z}^{(p-1)/2}$ and $v=\widetilde{h}-b\widetilde{z}^{(p+1)/2}$. Then $[v,\widetilde{z}]=0,~[u,\widetilde{z}]=\widetilde{h}=v+b\widetilde{z}^{(p+1)/2}$ and 
    \[
    [u,v]=[\widetilde{y},v]=c\widetilde{z}^p-b\frac{p+1}{2}\widetilde{z}^{(p-1)/2}\widetilde{h}=-b\frac{p+1}{2}\widetilde{z}^{(p-1)/2}v.
    \]
    This implies 
    \[
    \widetilde{\Lambda}=\FF[\widetilde{z},v][u;(v+b\widetilde{z}^{(p+1)/2})\frac{\partial}{\partial\widetilde{z}}-b\frac{p+1}{2}\widetilde{z}^{(p-1)/2}v\frac{\partial}{\partial v}]
    \]
    and $v$ is normal in the sense of $v\widetilde{\Lambda}=\widetilde{\Lambda}v$. Moreover, we have 
    \[
    \widetilde{\Lambda}/(v)\cong\FF[\widetilde{z}][u;b\widetilde{z}^{(p+1)/2}\frac{\partial}{\partial\widetilde{z}}],
    \]
    which is a domain by \cite[Theorem 1.2.9]{MR01}. If $a_1u+a_2v=0$ for some $(a_1,a_2)\in\widetilde{\Lambda}^2$, reduction modulo $v$ gives $a_1=wv$ for some $w\in\widetilde{\Lambda}$. Then 
    \[
    0=wvu+a_2v=w(uv+b\frac{p+1}{2}\widetilde{z}^{(p-1)/2}v)+a_2v,
    \]
    which implies $a_2=-w(u+b\frac{p+1}{2}\widetilde{z}^{(p-1)/2})$ because $\widetilde{\Lambda}$ is also a domain. This gives $(a_1,a_2)\in\Im d_2$.
\end{proof}

\begin{lemma}\label{lem:mixed-N-Ext}
For $i\ge0,~j\in\mathbb{Z}$, we have
\[
\Ext_{\widetilde{\Lambda}}^i(\widetilde{N_+},\widetilde{N_-})^j=
         \begin{cases}
             \alpha^{-j},&i=1,~-\frac{p-1}{2}\le j<0;\\
             \alpha^{-j},&i=2,~-p\le j<-\frac{p+1}{2};\\
             \FF,&i=1,2,~j=-\frac{p+1}{2},~r=\frac{p-1}{2};\\
             0,&\text{otherwise}.            
         \end{cases}
\]
The same assertion holds after interchanging $\widetilde{N_+}$ and $\widetilde{N_-}$.
\end{lemma}

\begin{proof}
    Applying $\Hom_{\widetilde{\Lambda}}(\cdot,\widetilde{N_-})$ to the resolution \eqref{eq:N+-resolution}, we obtain
    \[
    \FF[\widetilde{z}]\xrightarrow{d^1_-}\left(\alpha^{-1}\otimes\FF[\widetilde{z}](\frac{p-1}{2})\right)\oplus\FF[\widetilde{z}](\frac{p+1}{2})\xrightarrow{d^2_-}\alpha^{-1}\otimes\FF[\widetilde{z}](p),
    \]
    where 
    \[
    d^1_-(f)=(-b\widetilde{z}^{(p+1)/2}\frac{\partial f}{\partial\widetilde{z}}+2b\vartheta\widetilde{z}^{(p-1)/2}f,-2b\widetilde{z}^{(p+1)/2}f),
    \]
    \[
    d^2_-(U,V)=-2b\widetilde{z}^{(p+1)/2}U+b\widetilde{z}^{(p+1)/2}\frac{\partial V}{\partial\widetilde{z}}-b(2\vartheta+\frac{p+1}{2})\widetilde{z}^{(p-1)/2}V.
    \]
    The cochain complex decomposes with respect to the degree. For $j\ge0$, the degree $j$ part is
    \[
    \FF\widetilde{z}^j\xrightarrow{d^1_-}(\alpha^{-1}\otimes\FF{\widetilde{
    z}}^{j+(p-1)/2})\oplus\FF{\widetilde{
    z}}^{j+(p+1)/2}\xrightarrow{d^2_-}\alpha^{-1}\otimes\FF{\widetilde{
    z}}^{j+p},
    \]
    which is exact. For $-\frac{p-1}{2}\le j<0$, the degree $j$ part is
    \[
    0\longrightarrow(\alpha^{-1}\otimes\FF{\widetilde{
    z}}^{j+(p-1)/2})\oplus\FF{\widetilde{
    z}}^{j+(p+1)/2}\xrightarrow{d^2_-}\alpha^{-1}\otimes\FF{\widetilde{
    z}}^{j+p},
    \]
    which gives $\Ext_{\widetilde{\Lambda}}^1(\widetilde{N_+},\widetilde{N_-})^j=\alpha^{-j}$ and $\Ext_{\widetilde{\Lambda}}^2(\widetilde{N_+},\widetilde{N_-})^j=0$. For $j=-\frac{p+1}{2}$, the degree $j$ part is $d^2_-:\FF\to\alpha^{-1}\otimes\FF{\widetilde{
    z}}^{(p-1)/2}$, sending $1$ to $-b(2\vartheta+\frac{p+1}{2})\widetilde{z}^{(p-1)/2}$. Note that $2\vartheta+\frac{p+1}{2}=0$ if and only if $r=\frac{p-1}{2}$. We obtain 
    \[
    \Ext_{\widetilde{\Lambda}}^1(\widetilde{N_+},\widetilde{N_-})^{-(p+1)/2}=\Ext_{\widetilde{\Lambda}}^2(\widetilde{N_+},\widetilde{N_-})^{-(p+1)/2}=0
    \]
    if $r\ne\frac{p-1}{2}$ and 
    \[
    \Ext_{\widetilde{\Lambda}}^1(\widetilde{N_+},\widetilde{N_-})^{-(p+1)/2}=\Ext_{\widetilde{\Lambda}}^2(\widetilde{N_+},\widetilde{N_-})^{-(p+1)/2}=\FF
    \]
    if $r=\frac{p-1}{2}$. For $-p\le j<-\frac{p+1}{2}$, we have $\Ext_{\widetilde{\Lambda}}^2(\widetilde{N_+},\widetilde{N_-})^j=\alpha^{-j}$.
\end{proof}

\begin{proposition}\label{prop:further-decomposition}
There is a unique closed filtered $(\Lambda,\cF)$-module decomposition 
\[
M_z'\cong N_+\oplus N_-
\]
such that 
\[
\gr_{\cF}(N_+)\cong\widetilde{N_+},\quad\gr_{\cF}(N_-)\cong\widetilde{N_-}.
\]
Moreover, the decomposition is compatible with the $H'$-action, and the elements in $H\setminus H'$ permute $N_+$ and $N_-$.
\end{proposition}

\begin{proof}
    This follows by Lemma~\ref{lem:decomposition-criterion} and Lemma~\ref{lem:mixed-N-Ext}. The assertions of the $H'$-action and elements in $H\setminus H'$ follow by the uniqueness of the decomposition.
\end{proof}

\begin{corollary}\label{cor:mix-N-ext}
    For $i\ge0$, we have 
    \[
    \Ext^i_{\Lambda}(M_z',M_z')^H\cong\Ext^i_{\Lambda}(N_+,N_+)^{H'}\oplus\Ext^i_{\Lambda}(N_+,N_-)^{H'}.
    \]
    Moreover, we have 
    \[
    \dim\Ext^i_{\Lambda}(N_+,N_-)^{H'}= 
    \begin{cases}
        1,&i=1,2,~r=\frac{p-1}{2};\\
        0,&\text{otherwise}.            
    \end{cases}
    \]
\end{corollary}

\begin{proof}
    This follows by Lemma~\ref{lem:mixed-N-Ext} and Proposition~\ref{prop:further-decomposition}, with the similar proof to that of Corollary~\ref{cor:pi-induction}, Corollary~\ref{cor:Ext-decomposition} and Proposition~\ref{prop:mixed-ext-dim}.
\end{proof}

\begin{lemma}\label{lem:N-dual}
    The right $(\Lambda,\cF)$-module $\E^2_\Lambda(N_+)$ is isomorphic to $\psi^{-1}\alpha^{-1}\otimes\FF[[Z]](p)$ as a filtered $\FF[[Z]]$-module with compatible $H'$-action.
\end{lemma}

\begin{proof}
    Applying $\Hom_{\widetilde{\Lambda}}(\cdot,\widetilde{\Lambda})$ to \eqref{eq:N+-resolution}, we obtain
    \[
    \psi^{-1}\otimes\widetilde{\Lambda}\longrightarrow\left(\psi^{-1}\alpha^{-1}\otimes\widetilde{\Lambda}(\frac{p-1}{2})\right)\oplus\left(\psi^{-1}\otimes\widetilde{\Lambda}(\frac{p+1}{2})\right)\longrightarrow\psi^{-1}\alpha^{-1}\otimes\widetilde{\Lambda}(p).
    \]
    An argument similar to the proof of Lemma~\ref{lem:N+-rsolution} gives that $\widetilde{N_+}$ is Cohen-Macaulay of grade $2$ with 
    \[
    E^2_{\widetilde{\Lambda}}(\widetilde{N_+})\cong\psi^{-1}\alpha^{-1}\otimes\FF[\widetilde{z}](p)
    \]
    as a right $\widetilde{\Lambda}$-module with compatible $H'$-action, with the $\widetilde{\Lambda}$-action on the right-hand side 
    \[
    1\cdot\widetilde{h}=b\widetilde{z}^{(p+1)/2},\quad1\cdot\widetilde{y}=-b(\vartheta+\frac{p+1}{2})\widetilde{z}^{(p-1)/2}.
    \]
    The same argument as in the proof of Lemma~\ref{lem:z-dual-z} and Proposition~\ref{prop:leading-coefficient} implies the result.
\end{proof}

\subsection{Computation of \texorpdfstring{$\Ext^i_\Lambda(N,N)^{H'}$}{ExtiLambda(N,N)^{H'}}}

Throughout this subsection, let $N=N_+$, which is free of rank $1$ over $\FF[[Z]]$. Lift \eqref{eq:N+-resolution} to a strict filt-free $(\Lambda,\cF)$-module resolution $P_{N,\bullet}$ of $N$ with compatible $H'$-action of length 2. Let $N^\circ=\FF((Z))\otimes_{\FF[[Z]]}N$ equipped with the tensor-product filtration $\cF^nN^\circ$. By Proposition~\ref{prop:Laurent-extension}, the filtered action of $(\Lambda,\cF)$ extends to $N^\circ=\FF((Z))$ uniquely. Write $\widetilde{N^\circ}=\gr_{\cF}(N^\circ)$.

Our first goal is to prove the following proposition.

\begin{proposition}\label{prop:Hom-vanish}
    $\Hom_{\Lambda}(N,N^\circ/N)^{H'}=0$.
\end{proposition}
We first recall some group-theoretic input from \cite[\S~2.1, Proposition~2.1]{HW24b}. Put 
\[
G_n=U_D^nZ_D^1/Z_D^1
\]
for $n\ge1$. Choose $\xi\in\FF_q^\times$ with $\xi^p=-\xi$. Write 
\[
g_1=1+\varpi_D,\quad g_2=1+\varpi_D[\xi],\quad g_3=1+p[\xi]
\]
for their images in $G_1=U_D^1/Z_D^1$. The images of $g_1^{p^j},g_2^{p^j}$ form an $\FF_p$-basis of $G_{2j+1}/G_{2j+2}$, and the image of $g_3^{p^j}$ forms an $\FF_p$-basis of $G_{2j+2}/G_{2j+3}$. Moreover, the $p$-power map identifies $G_{2j+1}/G_{2j+2}$ (resp. $G_{2j+2}/G_{2j+3}$) with $G_{2j+3}/G_{2j+4}$ (resp. $G_{2j+4}/G_{2j+5}$). 

\begin{lemma}\label{lem:evaluation-1}
    For $s,v\ge0$, we have
    \[
    [g_2^{p^s},g_1^{p^v}]_{\rm grp}\equiv g_3^{-2p^{s+v}}\pmod{G_{2(s+v)+3}},\quad [g_3^{p^s},g_1^{p^v}]_{\rm grp}\equiv g_2^{-2p^{s+v+1}}\pmod{G_{2(s+v)+4}}.
    \]
    Here $[g,g']_{\rm grp}:=gg'g^{-1}g'^{-1}$.
\end{lemma}

\begin{proof}
    Using binomial expansion, we have
    \[
    g_1^{p^v}=1+p^v\varpi_D\pmod{\varpi_D^{2v+2}},\quad g_2^{p^s}=1+p^s\varpi_D[\xi]\pmod{\varpi_D^{2s+2}},
    \]
    \[
    g_3^{-2p^{s+v}}=1-2p^{s+v+1}[\xi]\pmod{\varpi_D^{2s+2v+3}}
    \]
    (here we need $p\ge5$). Therefore, we have
    \begin{align*}
            [g_2^{p^s},g_1^{p^v}]_{\rm grp}&=1+(g_2^{p^s}g_1^{p^v}-g_1^{p^v}g_2^{p^s})g_2^{-p^s}g_1^{-p^v}\\
            &=1+\left((g_2^{p^s}-1)(g_1^{p^v}-1)-(g_1^{p^v}-1)(g_2^{p^s}-1)\right)g_2^{-p^s}g_1^{-p^v}\\
            &\equiv1+p^{s+v}(\varpi_D[\xi]\varpi_D-p[\xi])\pmod{\varpi_D^{2(s+v)+3}}\\
            &\equiv g_3^{-2p^{s+v}}\pmod{\varpi_D^{2(s+v)+3}}.
    \end{align*}
    The second congruence can be obtained similarly.
\end{proof}

Lift the $\FF[\widetilde{z}]$-basis of $\widetilde{N}$ to $e\in N$ with $H'$-character $\psi$, such that $N=\FF[[Z]]e$ and $\cF^n N=Z^{\max\{n,0\}}\FF[[Z]]e$. Then $N^\circ=\FF((Z))e$ and $\cF^n N^\circ=Z^n\FF[[Z]]e$.

\begin{lemma}\label{lem:evaluation-2}
    For any $k\in\ZZ$, we have
    \[
    (g_1-1)Z^ke\equiv-Z^{k+1}e\pmod{Z^{k+2}N},\quad
    (g_2-1)Z^ke\equiv\xi Z^{k+1}e\pmod{Z^{k+2}N},
    \]
    \[
    (g_3-1)Z^ke\equiv-\xi bZ^{k+(p+1)/2}e\pmod{Z^{k+(p+3)/2}N}.
    \]
\end{lemma}

\begin{proof}
    The map from $\FF_q\cong U_D^1/U_D^2$ to $\mathfrak m_D/\mathfrak m_D^2$ sending $a\in\FF_q$ to $(1+\varpi_D[a])-1$ is $\FF_p$-linear, so it must have the form $aA+a^pB$. By \eqref{eq:YZ}, we have
    \[
    y=\sum_{a\in\FF_q^\times}a^{-1}\left((1+\varpi_D[a])-1\right)=\sum_{a\in\FF_q^\times}a^{-1}(aA+a^pB)=-A.
    \]
    Similarly we have $z=-B$. Therefore, for $a\in\FF_q$, we have
    \[
    \gr_{\mathfrak m_D}^1((1+\varpi_D[a])-1)=-ay-a^pz.
    \]
    In particular, we have
    \[
    \gr_{\mathfrak m_D}^1(g_1-1)=-y-z,\quad\gr_{\mathfrak m_D}^1(g_2-1)=-\xi y+\xi z
    \]
    and therefore
    \[
    g_1-1\equiv-Z\pmod{\cF^2\Lambda},\quad g_2-1\equiv\xi Z\pmod{\cF^2\Lambda}.
    \]
    This gives the first two congruences.

    By \cite[Lemma~2.4]{HW24b}, we have $g_3-1\in\mathfrak{m}_D^2$ and $G_3-1\subseteq\mathfrak{m}_D^3$. The first congruence in Lemma~\ref{lem:evaluation-1} with $s=v=0$ gives
    \[
    g_3^{-2}-1\equiv [g_2,g_1]_{\rm grp}-1\equiv(g_2-1)(g_1-1)-(g_1-1)(g_2-1)\pmod{\mathfrak{m}_D^3},
    \]
    which implies $\gr_{\mathfrak m_D}^2(g_3-1)=-\xi h$.
    Since $g_3$ commutes with $H$, its PBW expansion consists of $H$-invariant monomials. Note that $ZY$ and $[Y,Z]$ are the only two nonconstant $H$-invariant PBW monomials of weight at most $\frac{p+1}{2}$, and $zy$ does not occur in $\gr_{\mathfrak m_D}^2(g_3-1)$. Therefore
    \[
    g_3-1\equiv-\xi[Y,Z]\pmod{\cF^{(p+3)/2}\Lambda}.
    \]
    The last congruence follows by the $\widetilde{h}$ action on $\widetilde{N}$ and Proposition~\ref{prop:Laurent-extension}.
\end{proof}

\begin{corollary}\label{cor:c}
    The constant $c$ in Proposition~\ref{prop:graded-algebra} is $-2$.
\end{corollary}

\begin{proof}
    By the proof of Lemma~\ref{lem:evaluation-2}, we have
    \[
    g_1-1\equiv-Y+ZP(Z)\pmod{\cF^{(p+1)/2}\Lambda}
    \]
    for some polynomial $P$. By \eqref{eq:Lambda-relation} and Lemma~\ref{lem:g-weight}, we have
    \[
    [g_3,g_1]_{\rm grp}-1\equiv-\xi[Y,[Y,Z]]\equiv-c\xi Z^p\pmod{\cF^{p+1}\Lambda}.
    \]
    Since $G_4$ is topologically generated by $g_3^p$, $g_1^{p^2}$ and $g_2^{p^2}$, we have $G_4-1\subseteq\cF^{p+1}\Lambda$. Therefore, by Lemma~\ref{lem:evaluation-1}, we have
    \[
    [g_3,g_1]_{\rm grp}-1\equiv g_2^{-2p}-1\equiv2\xi Z^p\pmod{\cF^{p+1}\Lambda}.
    \]
    This gives $c=-2$.
\end{proof}

We say an operator $A$ on $N^\circ$ has order $s$ if $A(Z^ke)\in Z^{k+s}N$ and the coefficient of $Z^{k+s}e$ is nonzero for all $k\in\ZZ$. Write $\ord(A)=s$. For example, by Lemma~\ref{lem:evaluation-2}, we have
\[
\ord(g_1-1)=\ord(g_2-1)=1,\quad\ord(g_3-1)=\frac{p+1}{2}.
\]

\begin{lemma}\label{lem:evaluation-3}
    If $g\in G_{2j+1}\setminus G_{2j+2}$, then $\ord(g-1)=p^j$. If $g\in G_{2j+2}\setminus G_{2j+3}$, then $\ord(g-1)=\frac{p+1}{2} p^j$.
\end{lemma}

\begin{proof}
    As the images of $g_1^{p^j},g_2^{p^j}$ form an $\FF_p$-basis of $G_{2j+1}/G_{2j+2}$ and that of $g_3^{p^j}$ forms an $\FF_p$-basis of $G_{2j+2}/G_{2j+3}$, we can expand $g$ with respect to the basis. If $g\in G_{2j+2}\setminus G_{2j+3}$, then 
    \[
    g\equiv g_3^{ap^j}\pmod{G_{2j+3}}
    \]
    for some $a\in\FF_p^{\times}$ and by Lemma~\ref{lem:evaluation-2}, we have
    \[
    \ord(g_3^{ap^j}-1)=\ord(g_3^{p^j}-1)=\ord\left((g_3-1)^{p^j}\right)=\frac{p+1}{2} p^j.
    \]
    If $g\in G_{2j+1}\setminus G_{2j+2}$, then 
    \[
    g\equiv g_1^{ap^j}g_2^{bp^j}\pmod{G_{2j+2}}
    \]
    for some $a,b\in\FF_p$ which are not both zero, and by Lemma~\ref{lem:evaluation-2}, we have
    \[
    \ord(g_1^{ap^j}g_2^{bp^j}-1)=\ord(g_1^{ap^j}-1+g_2^{bp^j}-1)=p^j
    \]
    because the leading coefficient $-a+\xi^{p^j}b$ is nonzero. The terms in deeper group degrees do not change the order because $p^j<\frac{p+1}{2}p^j<p^{j+1}$, which gives the result. 
\end{proof}

Put $T=g_1-1$ as an operator on $N^\circ$. By Lemma~\ref{lem:evaluation-2}, $T$ is invertible on $N^\circ$, and for any $k\in\ZZ$, we have
\[
T^ke\equiv(-1)^kZ^ke\pmod{Z^{k+1}N}.
\]
Consequently, $\{T^ke:k\ge0\}$ form a topological basis of $N$, and the classes of $\{T^{-m}e:m\ge1\}$ form a basis of $N^\circ/N$.

\begin{lemma}\label{lem:evaluation-4}
    For $m=ap^v>0$ with $p\nmid a$ and $s\ge0$, we have
    \[
    \ord[g_2^{p^s},T^{-m}]=(\frac{p+1}{2}p^s-1)p^v-m,\quad
    \ord[g_3^{p^s},T^{-m}]=(p^{s+1}-1)p^v-m.
    \]
\end{lemma}

\begin{proof}
    By Lemma~\ref{lem:evaluation-1}, Lemma~\ref{lem:evaluation-2}, Lemma~\ref{lem:evaluation-3} and the identities
    \[
    T^{p^v}=g_1^{p^v}-1,\quad[g,g']=([g,g']_{\rm grp}-1)g'g,
    \]
    we have
    \[
    \ord[g_2^{p^s},T^{p^v}]=\frac{p+1}{2}p^{s+v},\quad
    \ord[g_3^{p^s},T^{p^v}]=p^{s+v+1}
    \]
    for $s,v\ge0$. For $m=ap^v>0$ with $p\nmid a$, we have
    \[
    [g,(T^{p^v})^{-a}]
    =-\sum_{i=1}^a(T^{p^v})^{-i}[g,T^{p^v}](T^{p^v})^{-(a+1-i)}.
    \]
    Note that
    \[
    \ord\left((T^{p^v})^{-i}[g,T^{p^v}](T^{p^v})^{-(a+1-i)}\right)=\ord([g,T^{p^v}])-(a+1)p^v
    \]
    with the same leading coefficient for all $1\le i\le a$, and $a\ne0$ in $\FF$ because $p\nmid a$. This gives 
    \[
    \ord( [g,(T^{p^v})^{-a}])=\ord([g,T^{p^v}])-(a+1)p^v
    \]
    which implies the result.
\end{proof}

\begin{proof}[Proof of Proposition~\ref{prop:Hom-vanish}]
    Suppose that 
    \[
    0\ne\varphi\in\Hom_\Lambda(N,N^\circ/N)^{H'}.
    \]
    Then $\varphi$ commutes with $T$. Since $\{T^ke:k\ge0\}$ form a topological basis of $N$ and the classes of $\{T^{-m}e:m\ge1\}$ form a basis of $N^\circ/N$, there exists a unique finite sum
    \[
    f(T)=\sum_{m\ge1}a_mT^{-m}
    \]
    such that 
    \[
    \varphi(u)=f(T)u\pmod{N}
    \]
    for all $u\in N$. The $\Lambda$-equivariance implies
    \begin{equation}\label{eq:Lambda-equi}
        [g_2^{p^s},f(T)](N)\subseteq N,\quad [g_3^{p^s},f(T)](N)\subseteq N
    \end{equation}
    for any $s\ge0$.
    
    Let $S$ be the finite set of valuations $v_p(m)$ for which $a_m\ne0$. For $v\in S$, put
    \[
    m_v=\max\{m:a_m\ne0,\ v_p(m)=v\},\quad L_v(x)=m_v-xp^v,\quad
    F(x)=\max_{v\in S}L_v(x).
    \]
    For $s\ge0$, define
    \[
    \tau_{2s}=\frac{p+1}{2}p^s-1,\quad \tau_{2s+1}=p^{s+1}-1.
    \]
    If $F(\tau_j)>0$ for some $j\ge0$, then there exist at least two distinct indices $v\in S$ such that $F(\tau_j)=L_v(\tau_j)$. Otherwise, by Lemma~\ref{lem:evaluation-4}, this contradicts \eqref{eq:Lambda-equi}.

    We claim that $F(\tau_0)\le0$. Otherwise, let $v_{\max}\in S$ be the largest index such that $F(\tau_0)=L_{v_{\max}}(\tau_0)$. Then $v_{\max}\ge1$. Since $p^{v_{\max}}\mid m_{v_{\max}}$, we have
    \[
    F(\tau_0)=m_{v_{\max}}-\frac{p-1}{2}p^{v_{\max}}\ge p^{v_{\max}}.
    \]
    Assume $v_0\in S$ is another index such that $F(\tau_0)=L_{v_0}(\tau_0)$. Then $v_0\le v_{\max}-1$ and
    \[
    F(\tau_0)-F(\tau_1)\le L_{v_0}(\tau_0)-L_{v_0}(\tau_1)=p^{v_0}(\tau_1-\tau_0)\le\frac{p-1}{2}p^{v_{\max}-1},
    \]
    which implies $F(\tau_1)>0$. Note that for $v\in S$ and $v\ge v_{\max}$, we have
    \[
    L_v(\tau_1)=m_v-\tau_1p^v=L_v(\tau_0)+(\tau_0-\tau_1)p^v<L_{v_0}(\tau_0)+(\tau_0-\tau_1)p^{v_0}=L_{v_0}(\tau_1).
    \]
    This implies that if $F(\tau_1)=L_v(\tau_1)$ for some $v\in S$, then we must have $v\le v_{\max}-1$. Repeat the analysis above. Using
    \[
    \tau_{j+1}-\tau_j=\frac{p-1}{2}p^{\lceil j/2\rceil}
    \]
    and
    \[
    \frac{p-1}{2}\sum_{j=0}^{v_{\max}-1}p^{v_{\max}-1-j+\lceil j/2\rceil}
    <(p-1)\sum_{a\ge0}p^{v_{\max}-1-a}=p^{v_{\max}},
    \]
    we obtain inductively that for $0\le j\le v_{\max}-1$
    \[
    F(\tau_j)-F(\tau_{j+1})\le\frac{p-1}{2}p^{v_{\max}-1-j+\lceil j/2\rceil},\quad F(\tau_{j+1})>0
    \]
    and if $F(\tau_{j+1})=L_v(\tau_{j+1})$ for some $v\in S$, then $v\le v_{\max}-j-1$. In particular, we have $F(\tau_{v_{\max}})>0$ and $v=0$ is the only possible index such that $F(\tau_{v_{\max}})=L_v(\tau_{v_{\max}})$. This is a contradiction. Consequently, we have $F(\tau_0)\le0$, which implies
    \[
    a_m\ne0\quad\implies\quad
    m\le\frac{p-1}{2}p^{v_p(m)}.
    \]

    Finally, let $m$ be the largest index such that $a_m\ne0$. Rewrite 
    \[
    \varphi(e)=f(T)e\pmod{N}
    \]
    in the $Z$-coordinate. Then $m$ is still the largest index such that the coefficient of $Z^{-m}$ is nonzero. However, since $\varphi$ is $H'$-invariant and the order of $\alpha|_{H'}$ is $\frac{p+1}{2}$, we must have 
    \[
    \frac{p+1}{2}\mid m,
    \]
    which is a contradiction.
\end{proof}

Next we compute $\Ext_{\Lambda}^2(N,N)^{H'}$.

\begin{proposition}\label{prop:dual}
    We have a $\Lambda$-module isomorphism with compatible $H'$-action 
    \[
    \E_{\Lambda}^2(N)^\vee\cong N^\circ/N.
    \]
\end{proposition}

\begin{proof}
    Since $U_D^1/Z_D^1$ is a torsion-free $p$-adic analytic pro-$p$ group of dimension $3$, it is a Poincaré group of dimension $3$. Applying the local duality theorem \cite[Theorem 7.2]{TV}, we obtain
    \[
    \E_{\Lambda}^2(N)^\vee\cong H_{\mathfrak{m}_D}^1(N),
    \]
    where $H_{\mathfrak{m}_D}^i(N):=\varinjlim_n\Ext_{\Lambda}^i(\Lambda/\mathfrak{m}_D^n,N)$ is the local cohomology. 

    We show $H_{\mathfrak{m}_D}^i(N^\circ)=0$ for $i\ge0$. It suffices to show 
    \[
    \Ext_{\Lambda}^i(\Lambda/\mathfrak{m}_D^n,N^\circ)=0
    \]
    for $n\ge0$. Endow $\Lambda/\mathfrak{m}_D^n$ with the quotient filtration of the weighted filtration $\cF^n\Lambda$, so that $\Lambda/\mathfrak{m}_D^n$ becomes a finite-dimensional filtered $(\Lambda,\cF)$-module. Choose a strict filt-free $(\Lambda,\cF)$-module resolution $P_\bullet$ of $\Lambda/\mathfrak{m}_D^n$ with compatible $H'$-action. Then $\Hom_\Lambda(P_\bullet,N^\circ)$ determines a filtered spectral sequence $E_\bullet^{\bullet\bullet}$ with 
    \[
    E_1^{i,j}\cong\Ext_{\widetilde{\Lambda}}^{i+j}\left(\gr_{\cF}(\Lambda/\mathfrak{m}_D^n),\widetilde{N^\circ}\right)^i
    \]
    (note that we cannot apply Proposition~\ref{prop:spectral-sequence} to conclude that $E_\bullet^{\bullet\bullet}$ converges to $\Ext_{\Lambda}^\bullet(\Lambda/\mathfrak{m}_D^n,N^\circ)$ because the filtration on $N^\circ$ is not bounded below). Consider the central element $\widetilde{z}^p\in\widetilde{\Lambda}$. It acts invertibly on $\widetilde{N^\circ}$ and nilpotently on $\gr_{\cF}(\Lambda/\mathfrak{m}_D^n)$, and thus 
    \[
    \Ext_{\widetilde{\Lambda}}^{i+j}\left(\gr_{\cF}(\Lambda/\mathfrak{m}_D^n),\widetilde{N^\circ}\right)=0.
    \]
    Since the filtration on $N^\circ$ is complete and exhaustive, by the Eilenberg-Moore Comparison Theorem \cite[Theorem 5.5.11]{WC}, we have $\Ext_{\Lambda}^{i+j}(\Lambda/\mathfrak{m}_D^n,N^\circ)=0$.

    Consider the short exact sequence of $\Lambda$-modules with compatible $H'$-action
    \[
    0\longrightarrow N\longrightarrow N^\circ\longrightarrow N^\circ/N\longrightarrow0.
    \]
    The long exact sequence associated with it gives $H_{\mathfrak{m}_D}^1(N)\cong H_{\mathfrak{m}_D}^0(N^\circ/N)$. The result follows from $H_{\mathfrak{m}_D}^0(N^\circ/N)=N^\circ/N$ because $N^\circ/N$ is $\mathfrak{m}_D$-torsion.
\end{proof}

\begin{proposition}\label{prop:ext2-vanish}
    $\Ext_{\Lambda}^2(N,N)^{H'}=0$.
\end{proposition}

\begin{proof}
    Since $N$ admits a free resolution $P_{N,\bullet}$ of length 2, we have
    \[
    \Ext_{\Lambda}^2(N,N)\cong\E_{\Lambda}^2(N)\otimes_\Lambda N.
    \]
    The adjunction formula and Proposition~\ref{prop:dual} give
    \[
    \Ext_{\Lambda}^2(N,N)^\vee\cong\Hom_\Lambda\left(N,\E_{\Lambda}^2(N)^\vee\right)\cong\Hom_\Lambda(N,N^\circ/N).
    \]
    The result follows by Proposition~\ref{prop:Hom-vanish}.
\end{proof}

Now we turn to the computation of $\End_{\Lambda}(N)^{H'}$ and $\Ext_{\Lambda}^1(N,N)^{H'}$. We first introduce the spectral sequences used in the computation. Let 
\[
C^\bullet:=\Hom_\Lambda(P_{N,\bullet},N)^{H'},\quad{C^\circ}^\bullet:=\Hom_\Lambda(P_{N,\bullet},N^\circ)^{H'}.
\]
Then $C^\bullet$ and ${C^\circ}^\bullet$ determine filtered spectral sequences $\left(E_\bullet^{\bullet\bullet}(C),d_\bullet^C\right)$ and $\left(E_\bullet^{\bullet\bullet}(C^\circ),d_\bullet^{C^\circ}\right)$ respectively with 
\[
E_1^{i,j}(C)\cong\Ext_{\widetilde{\Lambda}}^{i+j}(\widetilde{N},\widetilde{N})^{i,H'},\quad E_1^{i,j}(C^\circ)\cong\Ext_{\widetilde{\Lambda}}^{i+j}(\widetilde{N},\widetilde{N^\circ})^{i,H'}.
\]
Moreover, by Proposition~\ref{prop:spectral-sequence}, $E_\bullet^{\bullet\bullet}(C)$ converges to $\Ext^\bullet_{\Lambda}(N,N)^{H'}$. The inclusion
\[
\iota:N\hookrightarrow N^\circ
\]
induces a filtered chain map from $C^\bullet$ to ${C^\circ}^\bullet$, which further induces a morphism of spectral sequences 
\[
\iota_r:E_r^{\bullet\bullet}(C)\longrightarrow E_r^{\bullet\bullet}(C^\circ).
\]

%We want to define a multiplicative structure on $E_r^{i,j}(C)$. Let $D^\bullet:=\Hom_{\Lambda}^\bullet(P_{N,\bullet},P_{N,\bullet})^{H'}$ be the $\Hom$ complex with $D^m:=\prod_n\Hom_{\Lambda}(P_{N,n},P_{N,n-m})^{H'}$ and differential $\partial(f)=d_N\circ f-(-1)^mf\circ d_N$. Then $D^\bullet$ determines a filtered spectral sequence $E_r^{i,j}(D)$ with $E_1^{i,j}(D)\cong H^{i+j}\left(\Hom_{\widetilde{\Lambda}}^\bullet(P_{\widetilde{N},\bullet},P_{\widetilde{N},\bullet})^{H'}\right)_i$. The composition of morphism gives a multiplication $D^i\times D^j\to D^{i+j}$, which satisfies the Leibniz rule
%\[
%\partial(f\circ g)=\partial(f)\circ g+(-1)^if\circ\partial(g).
%\]
%Therefore, $D^\bullet$ is a filtered differential graded algebra, and by \cite[\S~5.4.8]{WC} the multiplication above induces a multiplicative structure $E_r^{i,j}(D)\times E_r^{i',j'}(D)\to E_r^{i+i',j+j'}(D)$ for every $r\ge1$, satisfying the Leibniz rule. Moreover, by the proof of Lemma~\ref{lem:decomposition-criterion}, $D^\bullet$ is quasi-isomorphic to $C^\bullet$, and it induces a quasi-isomorphism between $\Hom_{\widetilde{\Lambda}}^\bullet(P_{\widetilde{N},\bullet},P_{\widetilde{N},\bullet})^{H'}$ and $\Hom_{\widetilde{\Lambda}}(P_{\widetilde{N},\bullet},\widetilde{N})^{H'}$. This induces $E_1^{i,j}(D)\cong E_1^{i,j}(C)$ and further $E_r^{i,j}(D)\cong E_r^{i,j}(C)$ for any $r\ge1$. This gives a multiplicative structure $E_r^{i,j}(C)\times E_r^{i',j'}(C)\to E_r^{i+i',j+j'}(C)$. Moreover, by \cite[\S~2.6]{BK}, the multiplication on $E_1^{\bullet\bullet}(C)$ is the usual Yoneda product.

The Yoneda product
\[
\Ext^i_{\Lambda}(N,N)^{H'}\times\Ext^j_{\Lambda}(N,N)^{H'}\longrightarrow\Ext^{i+j}_{\Lambda}(N,N)^{H'}
\]
induces a multiplicative structure
\[
E_r^{i,j}(C)\times E_r^{i',j'}(C)\longrightarrow E_r^{i+i',j+j'}(C)
\]
for every $r$, satisfying the Leibniz rule
\[
d_r^C(\lambda_1\lambda_2)=d_r^C(\lambda_1)\lambda_2+(-1)^{i+j}\lambda_1d_r^C(\lambda_2)
\]
and compatible with taking cohomology. Then $E_r^{\bullet\bullet}(C)$ becomes a differential graded algebra. For $r=1$, it coincides with the Yoneda product
\[
\Ext^i_{\widetilde{\Lambda}}(\widetilde{N},\widetilde{N})^{H'}\times\Ext^j_{\widetilde{\Lambda}}(\widetilde{N},\widetilde{N})^{H'}\longrightarrow\Ext^{i+j}_{\widetilde{\Lambda}}(\widetilde{N},\widetilde{N})^{H'}.
\]
Similarly, the Yoneda product
\[
\Ext^i_{\Lambda}(N,N^\circ)^{H'}\times\Ext^j_{\Lambda}(N,N)^{H'}\longrightarrow\Ext^{i+j}_{\Lambda}(N,N^\circ)^{H'}
\]
induces a multiplicative structure
\[
E_r^{i,j}(C^\circ)\times E_r^{i',j'}(C)\longrightarrow E_r^{i+i',j+j'}(C^\circ)
\]
for every $r$, satisfying the Leibniz rule
\[
d_r^{C^\circ}(\lambda_1\lambda_2)=d_r^{C^\circ}(\lambda_1)\lambda_2+(-1)^{i+j}\lambda_1d_r^C(\lambda_2)
\]
and compatible with taking cohomology. Then $E_r^{\bullet\bullet}(C^\circ)$ becomes a differential graded right $E_r^{\bullet\bullet}(C)$-module. For $r=1$, it coincides with the Yoneda product
\[
\Ext^i_{\widetilde{\Lambda}}(\widetilde{N},\widetilde{N^\circ})^{H'}\times\Ext^j_{\widetilde{\Lambda}}(\widetilde{N},\widetilde{N})^{H'}\longrightarrow\Ext^{i+j}_{\widetilde{\Lambda}}(\widetilde{N},\widetilde{N^\circ})^{H'}.
\]
Moreover, for every $r$,
\[
\iota_r:E_r^{\bullet\bullet}(C)\longrightarrow E_r^{\bullet\bullet}(C^\circ)
\]
is a morphism of differential graded right $E_r^{\bullet\bullet}(C)$-modules.

\begin{lemma}\label{lem:E1-C}
As a differential graded algebra, we have
    \[
    E_1^{\bullet\bullet}(C)=\FF[u]\otimes\bigwedge(\FF\eta_1\oplus\FF\eta_2),
    \]
    where
    \[
    u\in E_1^{p(p+1)/2,-p(p+1)/2}(C),\quad\eta_1\in E_1^{(p+1)(p-2)/2,-(p+1)(p-2)/2+1}(C),\quad\eta_2\in E_1^{0,1}(C).
    \]
    Moreover, as a differential graded right $E_1^{\bullet\bullet}(C)$-module, we have
    \[
     E_1^{\bullet\bullet}(C^\circ)\cong E_1^{\bullet\bullet}(C)[u^{-1}]=\FF[u,u^{-1}]\otimes\bigwedge(\FF\eta_1\oplus\FF\eta_2).
    \]
\end{lemma}

\begin{proof}
    Applying $\Hom_{\widetilde{\Lambda}}(\cdot,\widetilde{N})$ to the resolution \eqref{eq:N+-resolution}, we obtain
    \[
    \FF[\widetilde{z}]\xrightarrow{d^1_+}\left(\alpha^{-1}\otimes\FF[\widetilde{z}](\frac{p-1}{2})\right)\oplus\FF[\widetilde{z}](\frac{p+1}{2})\xrightarrow{d^2_+}\alpha^{-1}\otimes\FF[\widetilde{z}](p),
    \]
    where 
    \[
    d^1_+(f)=(b\widetilde{z}^{(p+1)/2}\frac{\partial f}{\partial\widetilde{z}},0),\quad d^2_+(U,V)=-b\frac{\partial(\widetilde{z}^{(p+1)/2}V)}{\partial\widetilde{z}}.
    \]
    The complex decomposes with respect to the degree. An argument similar to that in the proof of Lemma~\ref{lem:mixed-N-Ext} gives
    \[
    \Ext_{\widetilde{\Lambda}}^i(\widetilde{N},\widetilde{N})^{j,H'}=
         \begin{cases}
             \FF,&i=0,1,~j\ge0,~\frac{p(p+1)}{2}\mid j;\\
             \FF,&i=1,2,~j\ge0,~p\mid j+1,~\frac{p+1}{2}\mid j;\\
             0,&\text{otherwise}.            
         \end{cases}
    \]
    
    Put 
    \[
    u=[\widetilde{z}^{p(p+1)/2}]\in\Hom_{\widetilde{\Lambda}}(\widetilde{N},\widetilde{N})^{p(p+1)/2,H'},
    \] 
    which is a central element in the Yoneda algebra $\Ext^\bullet_{\widetilde{\Lambda}}(\widetilde{N},\widetilde{N})^{H'}$ because $\widetilde{z}^{p(p+1)/2}$ is a central element in $\widetilde{\Lambda}$. Put 
    \[
    \eta_1=[(0,\widetilde{z}^{(p^2-1)/2})]\in\Ext_{\widetilde{\Lambda}}^1(\widetilde{N},\widetilde{N})^{(p+1)(p-2)/2,H'},\quad\eta_2=[(-\widetilde{z}^{(p-1)/2},0)]\in\Ext_{\widetilde{\Lambda}}^1(\widetilde{N},\widetilde{N})^{0,H'},
    \]
    \[
    \eta_3=[\widetilde{z}^{p(p+1)/2-1}]\in\Ext_{\widetilde{\Lambda}}^2(\widetilde{N},\widetilde{N})^{(p+1)(p-2)/2,H'}.
    \]
    Then 
    \[
    \Ext^\bullet_{\widetilde{\Lambda}}(\widetilde{N},\widetilde{N})^{H'}=\FF[u]\oplus\FF[u]\eta_1\oplus\FF[u]\eta_2\oplus\FF[u]\eta_3.
    \]
    To compute the Yoneda product, we need to lift $\eta_i$ to the chain level. Write the resolution \eqref{eq:N+-resolution} as
    \[
    0\longrightarrow\widetilde{P}_2\longrightarrow\widetilde{P}_1\longrightarrow\widetilde{P}_0\longrightarrow\widetilde{N}\longrightarrow0.
    \]
    Put
    \[
    s_1=
    \begin{pmatrix}
        0 \\
        \widetilde{z}^{(p^2-1)/2}
    \end{pmatrix}
    :\widetilde{P}_1\rightarrow\widetilde{P}_0,
    \quad
    s_1'=(\widetilde{z}^{(p^2-1)/2},\frac{p-1}{2}\widetilde{z}^{(p^2-3)/2}):\widetilde{P}_2\rightarrow\widetilde{P}_1,
    \]
    \[
    s_2=
    \begin{pmatrix}
        -\widetilde{z}^{(p-1)/2} \\
        0
    \end{pmatrix}
    :\widetilde{P}_1\rightarrow\widetilde{P}_0,
    \quad
    s_2'=(0,\widetilde{z}^{(p-1)/2}):\widetilde{P}_2\rightarrow\widetilde{P}_1.
    \]
    Then one can check
    \[
    \varepsilon\circ s_i=\widetilde{\eta}_i,\quad s_i\circ d_2+d_1\circ s_i'=0,\quad i=1,2.
    \]
    By the definition of the Yoneda product, we have
    \[
    \eta_1\eta_2=[\widetilde{\eta}_1\circ s_2']=[\widetilde{z}^{p(p+1)/2-1}]=\eta_3,\quad\eta_2\eta_1=[\widetilde{\eta}_2\circ s_1']=-\eta_3,\quad\eta_1^2=\eta_2^2=0.
    \]
    This gives the assertion for $E_1^{\bullet\bullet}(C)$. 
    
    Applying $\Hom_{\widetilde{\Lambda}}(\cdot,\widetilde{N^\circ})$ to the resolution \eqref{eq:N+-resolution}, we obtain
    \[
    \FF[\widetilde{z},\widetilde{z}^{-1}]\xrightarrow{d^1_\circ}\left(\alpha^{-1}\otimes\FF[\widetilde{z},\widetilde{z}^{-1}](\frac{p-1}{2})\right)\oplus\FF[\widetilde{z},\widetilde{z}^{-1}](\frac{p+1}{2})\xrightarrow{d^2_\circ}\alpha^{-1}\otimes\FF[\widetilde{z},\widetilde{z}^{-1}](p),
    \]
    where $d^i_\circ$ are given by the same formulas as $d^i_+$ by Proposition~\ref{prop:Laurent-extension}. We have     
    \[
    \Ext_{\widetilde{\Lambda}}^i(\widetilde{N},\widetilde{N^\circ})^{H'}_j=
         \begin{cases}
             \FF,&i=0,1,~\frac{p(p+1)}{2}\mid j;\\
             \FF,&i=1,2,~p\mid j+1,~\frac{p+1}{2}\mid j;\\
             0,&\text{otherwise}.            
         \end{cases}
    \]
    The remaining computation is similar.
\end{proof}

Let $\FF(p+1)$ be the filtered $\FF$-vector space with the filtration equal to $\FF$ for degree $\le-p-1$ and $0$ for degree $>-p-1$. Recall that $\Ext^2_{\Lambda}(N,N^\circ)$ can be computed using the resolution $P_{N,\bullet}$ and hence is a quotient of $N^\circ$.
\begin{lemma}\label{lem:residue}
    We have a filtered map 
    \[
    \tau:\Ext^2_{\Lambda}(N,N^\circ)^{H'}\longrightarrow\FF(p+1)
    \]
    which is nonzero on the class $[Z^{-1}]$.
\end{lemma}

\begin{proof}
    Define $\tau$ by taking $H'$-invariants in the composition
    \[
    \Ext^2_{\Lambda}(N,N^\circ)\xrightarrow{\pi_*}\Ext^2_{\Lambda}(N,N^\circ/N)\cong\E^2_{\Lambda}(N)\otimes_\Lambda (N^\circ/N)\xrightarrow{\langle\cdot,\cdot\rangle}\FF(p+1),
    \]
    where $\pi_*$ is induced by $\pi:N^\circ\to N^\circ/N$ and $\langle\cdot,\cdot\rangle$ is induced by the isomorphism in Proposition~\ref{prop:dual}. By Lemma~\ref{lem:N-dual}, we have 
    \[
    \E^2_{\Lambda}(N)\cong\psi^{-1}\alpha^{-1}\otimes\FF[[Z]](p),
    \]
    which implies that $\tau$ is filtered. The image of the class $[Z^{-1}]\in\Ext^2_{\Lambda}(N,N^\circ)$ in $\E^2_{\Lambda}(N)\otimes_\Lambda (N^\circ/N)$ is $1\otimes Z^{-1}$. Since $Z^{-1}\in N^\circ/N$ is killed by $Z$, the map
    \[
    \langle\cdot,[Z^{-1}]\rangle:\E^2_{\Lambda}(N)\longrightarrow\FF
    \]
    factors through $\E^2_{\Lambda}(N)/\E^2_{\Lambda}(N)Z$, which is $1$-dimensional. Therefore $\langle1,[Z^{-1}]\rangle\ne0$ because the pairing $\langle\cdot,\cdot\rangle$ is nondegenerate.
\end{proof}

Let $K^\bullet$ be the filtered cochain complex of $\FF$-vector spaces with $K^n=0$ for $n\ne2$ and $K^2=\FF(p+1)$. Then $\tau$ defined in Lemma~\ref{lem:residue} can be upgraded to a filtered chain map from ${C^\circ}^\bullet$ to $K^\bullet$, which induces a morphism of spectral sequences $E_\bullet^{\bullet\bullet}(C^\circ)\to E_\bullet^{\bullet\bullet}(K)$. Note that on every page $r$, the only nonzero term of $E_r^{\bullet\bullet}(K)$ is $E_r^{-p-1,p+3}(K)=\FF$. Write 
\[
\tau_r:E_r^{-p-1,p+3}(C^\circ)\longrightarrow E_r^{-p-1,p+3}(K)=\FF.
\]
It satisfies
\begin{enumerate}
    \item $\tau_r\circ d_r^{C^\circ}:E_r^{-p-1-r,p+2+r}(C^\circ)\to\FF$ is zero.
    \item If $\omega_{r+1}\in E_{r+1}^{-p-1,p+3}(C^\circ)$ is the image of the cocycle $\omega_r\in E_r^{-p-1,p+3}(C^\circ)$, then we have $\tau_{r+1}(\omega_{r+1})=\tau_r(\omega_r)$.
\end{enumerate}

\begin{lemma}\label{lem:residue-nonzero}
    $\tau_1(\frac{\iota_1(\eta_1\eta_2)}{u})$ is nonzero.
\end{lemma}

\begin{proof}
    By Lemma~\ref{lem:E1-C}, we have 
    \[
    \frac{\iota_1(\eta_1\eta_2)}{u}=[\widetilde{z}^{-1}]\in\Ext_{\widetilde{\Lambda}}^2(\widetilde{N},\widetilde{N^\circ})^{-p-1,H'}.
    \]
    Therefore, by Lemma~\ref{lem:residue}
    \[
    \tau_1(\frac{\iota_1(\eta_1\eta_2)}{u})=\tau([Z^{-1}])\ne0.
    \]
\end{proof}

\begin{proposition}\label{prop:ext-N}
    We have 
    \[
    \Hom_\Lambda(N,N)^{H'}=\FF,\quad\Ext^1_{\Lambda}(N,N)^{H'}=0.
    \]
\end{proposition}

\begin{proof}
    We begin with $E_1^{\bullet\bullet}(C)$. First suppose that $d_1^C$ is zero. By Lemma~\ref{lem:E1-C}, any element in $E_1^{\bullet\bullet}(C^\circ)$ can be expressed as $\frac{\iota_1(x)}{u^k}$ for some $x\in E_1^{\bullet\bullet}(C)$ and $k\ge0$. Then by the Leibniz rule,
    \[
    0=\iota_1\circ d_1^{C}(x)=d_1^{C^\circ}(\frac{\iota_1(x)}{u^k}\cdot u^k)=d_1^{C^\circ}(\frac{\iota_1(x)}{u^k})\cdot u^k,
    \]
    which implies $d_1^{C^\circ}(\frac{\iota_1(x)}{u^k})=0$ because $u$ is not a zero divisor on $E_1^{\bullet\bullet}(C^\circ)$. Therefore, passing to the pages $E_2^{\bullet\bullet}(C)$ and $E_2^{\bullet\bullet}(C^\circ)$ preserves all structures established in Lemma~\ref{lem:E1-C}. The differential $d_r^C$ cannot be zero for all $r\ge1$. Otherwise the spectral sequence $E_r^{\bullet\bullet}(C)$ degenerates at page $E_1^{\bullet\bullet}(C)$, which contradicts Proposition~\ref{prop:ext2-vanish}. Therefore, there exists a minimal $r_1\ge1$ such that $d_{r_1}^C$ is nonzero.

    By Lemma~\ref{lem:E1-C} and assumption, 
    \[ 
    E_{r_1}^{\bullet\bullet}(C)\cong\FF[u]\otimes\bigwedge(\FF\eta_1\oplus\FF\eta_2)
    \]
    is a differential graded algebra generated by $u,\eta_1,\eta_2$. We claim that $d_{r_1}^C(u)$ is nonzero. Otherwise one of $d_{r_1}^C(\eta_1),~d_{r_1}^C(\eta_2)$ must be nonzero. We assume $d_{r_1}^C(\eta_1)$ is nonzero. Then 
    \[
    d_{r_1}^C(\eta_1)=\mu u^j\eta_1\eta_2
    \]
    for some $\mu\in\FF^\times$ and $j\ge0$. Consider 
    \[
    x=\frac{\iota_{r_1}(\eta_1)}{u^{j+1}}\in E_{r_1}^{\bullet\bullet}(C^\circ).
    \]
    Using $d_{r_1}^C(u)=0$ and the Leibniz rule, we have
    \[
    d_{r_1}^{C^\circ}(x)=\frac{\iota_{r_1}\left(d_{r_1}^C(\eta_1)\right)}{u^{j+1}}=\mu\frac{\iota_{r_1}(\eta_1\eta_2)}{u}\in E_{r_1}^{-p-1,p+3}(C^\circ).
    \]
    On the one hand, we have $\tau_{r_1}\circ d_{r_1}^{C^\circ}(x)=0$. On the other hand, we have
    \[
    \tau_{r_1}(\frac{\iota_{r_1}(\eta_1\eta_2)}{u})=\tau_1(\frac{\iota_1(\eta_1\eta_2)}{u})\ne0
    \]
    by Lemma~\ref{lem:residue-nonzero}, which is a contradiction. Therefore, we have $d_{r_1}^C(\eta_1)=0$. Similarly, we have $d_{r_1}^C(\eta_2)=0$. This is a contradiction.

    Next we determine the differential $d_{r_1}^C$. We know 
    \[
    d_{r_1}^C(u)\in\FF[u]\eta_1\oplus\FF[u]\eta_2,\quad\deg(u)=\frac{p(p+1)}{2},\quad\deg(\eta_1)=\frac{(p+1)(p-2)}{2},\quad\deg(\eta_2)=0.
    \]
    As $\deg(\eta_1)$ is not congruent to $\deg(\eta_2)$ modulo $\deg(u)$, we must have
    \[
    d_{r_1}^C(u)=\lambda u^{n+1}\eta_l
    \]
    for some $\lambda\in\FF^\times,~n\ge0$ and $l\in\{1,2\}$. Comparing degrees, we have
    \[
    \deg(u)+r_1=(n+1)\deg(u)+\deg(\eta_l)\implies r_1=n\deg(u)+\deg(\eta_l).
    \]
    We claim $n=0$. Otherwise, we have
    \[
    r_1\ge\deg(u)>\deg(\eta_1)+\deg(\eta_2).
    \]
    By Lemma~\ref{lem:E1-C}, any nonzero class in $E_r^{\bullet\bullet}(C)$ for $r\ge r_1$ has nonnegative degree, and therefore the class $\eta_1\eta_2$ cannot be the image of a differential $d_r^C$ for any $r\ge r_1$. Since $E_r^{i,j}(C)$ converges to $\Ext^{i+j}_{\Lambda}(N,N)^{H'}$, the class $\eta_1\eta_2$ survives to a nonzero class in $\Ext^2_{\Lambda}(N,N)^{H'}$, which contradicts Proposition~\ref{prop:ext2-vanish}. This also gives $l=1$ because $r_1\ge1$.
    
    Using $d_{r_1}^C\circ d_{r_1}^C(u)=0$, the Leibniz rule and $\eta_1^2=0$, we have
    \[
    0=\lambda d_{r_1}^C(u\eta_1)=\lambda d_{r_1}^C(u)\eta_1+\lambda u d_{r_1}^C(\eta_1)=\lambda u d_{r_1}^C(\eta_1),
    \]
    which implies $d_{r_1}^C(\eta_1)=0$. Assume
    \[
    d_{r_1}^C(\eta_2)=\lambda'\eta_1\eta_2
    \]
    for some $\lambda'\in\FF$. Consider
    \[
    x'=\frac{\iota_{r_1}(\eta_2)}{u}\in E_{r_1}^{\bullet\bullet}(C^\circ).
    \]
    Using the Leibniz rule, we have
    \[
    d_{r_1}^{C^\circ}(x')=\frac{\iota_{r_1}\left(d_{r_1}^C(\eta_2)\right)}{u}+\frac{\iota_{r_1}(\eta_2)}{u^2}\cdot d_{r_1}^C(u)=(\lambda'-\lambda)\frac{\iota_{r_1}(\eta_1\eta_2)}{u}.
    \]
    A similar argument using $\tau_{r_1}$ gives $\lambda'=\lambda$. 
    
    To conclude, we have
    \[
    d_{r_1}^C(u)=\lambda u\eta_1,\quad d_{r_1}^C(\eta_1)=0,\quad d_{r_1}^C(\eta_2)=\lambda\eta_1\eta_2
    \]
    for some $\lambda\in\FF^\times$. Using the Leibniz rule, we have
    \[
    d_{r_1}^C(u^k)=k\lambda u^k\eta_1,\quad d_{r_1}^C(u^k\eta_1)=k\lambda u^k\eta_1^2=0,\quad d_{r_1}^C(u^k\eta_2)=\lambda(k+1)u^k\eta_1\eta_2
    \]
    for any $k\ge0$. The same formulas hold for $d_{r_1}^{C^\circ}$ and $k\in\ZZ$ by the Leibniz rule. Taking cohomology, we obtain
    \[
    E_{r_1+1}^{\bullet\bullet}(C)\cong\FF[u_2]\otimes\bigwedge(\FF\eta_1^{(2)}\oplus\FF\eta_2^{(2)})
    \]
    as a differential graded algebra and
    \[
    E_{r_1+1}^{\bullet\bullet}(C^\circ)\cong\FF[u_2,u_2^{-1}]\otimes\bigwedge(\FF\eta_1^{(2)}\oplus\FF\eta_2^{(2)})
    \]
    as a differential graded right $E_{r_1+1}^{\bullet\bullet}(C)$-module, where
    \[
    u_2=[u^p],\quad\eta_1^{(2)}=[\eta_1],\quad\eta_2^{(2)}=[u^{p-1}\eta_2].
    \]
    We have
    \[
    \frac{\iota_{r_1+1}(\eta_1^{(2)}\eta_2^{(2)})}{u_2}=[\frac{\iota_{r_1}(\eta_1\eta_2)}{u}],\quad\deg(u_2)=p\deg(u)>\deg(\eta_1^{(2)})+\deg(\eta_2^{(2)})
    \]
    and $\deg(\eta_1^{(2)})$ is not congruent to $\deg(\eta_2^{(2)})$ modulo $\deg(u_2)$. 
    
    One checks that we can repeat the procedure above, finding the next page $E_{r_2}^{\bullet\bullet}(C)~(r_2>r_1)$ with nonzero differential $d_{r_2}^C$ and deducing $r_2=\deg(\eta_2^{(2)})$. Inductively, we have 
    \[
    E_{r_m}^{\bullet\bullet}(C)\cong\FF[u_m]\otimes\bigwedge(\FF\eta_1^{(m)}\oplus\FF\eta_2^{(m)})
    \]
    as a differential graded algebra and
    \[
    E_{r_m}^{\bullet\bullet}(C^\circ)\cong\FF[u_m,u_m^{-1}]\otimes\bigwedge(\FF\eta_1^{(m)}\oplus\FF\eta_2^{(m)})
    \]
    as a differential graded right $E_{r_m}^{\bullet\bullet}(C)$-module, where
    \[
    u_m=[u_{m-1}^p],\quad\eta_1^{(m)}=[\eta_1^{(m-1)}],\quad\eta_2^{(m)}=[u_{m-1}^{p-1}\eta_2^{(m-1)}],\quad r_m=\deg(\eta_2^{(m)})
    \]
    if $m$ is even and
    \[
    u_m=[u_{m-1}^p],\quad\eta_1^{(m)}=[u_{m-1}^{p-1}\eta_1^{(m-1)}],\quad\eta_2^{(m)}=[\eta_2^{(m-1)}],\quad r_m=\deg(\eta_1^{(m)})
    \]
    if $m$ is odd. Therefore, we have
    \[
    \lim_{m\to\infty}\deg(u_m)=\lim_{m\to\infty}\deg(\eta_1^{(m)})=\lim_{m\to\infty}\deg(\eta_2^{(m)})=\infty,
    \]
    which implies $E_{\infty}^{i,j}=0$ except for $(i,j)=(0,0)$. Since the spectral sequence $E_\bullet^{i,j}(C)$ converges to $\Ext^{i+j}_{\Lambda}(N,N)^{H'}$, we obtain
    \[
    \Hom_\Lambda(N,N)^{H'}=\FF,\quad\Ext^1_{\Lambda}(N,N)^{H'}=0.
    \]
\end{proof}

\begin{corollary}\label{cor:Mz'-ext}
    We have
    \[
    \dim_\FF\Ext^i_{\Lambda}(M_z',M_z')^{H}= 
    \begin{cases}
        1,&i=0;\\
        1,&i=1,2,~r=\frac{p-1}{2};\\
        0,&\text{otherwise}.            
    \end{cases}
    \]
\end{corollary}

\begin{proof}
    This follows by Corollary~\ref{cor:mix-N-ext} and Proposition~\ref{prop:ext-N}.
\end{proof}

\subsection{The main result}

In this subsection, we use the results obtained in the previous subsections to compute $\Ext^i_\Lambda(M_z,M_z)^H$ and $\Ext^i_{D^\times,\zeta}(\pi,\pi)$. Put $T_z:=M_z/M_z'$ and $M_z^\circ:=\FF((Z))\otimes_{\FF[[Z]]}M_z=\FF((Z))\otimes_{\FF[[Z]]}M_z'$.

\begin{lemma}\label{lem:Mz-dual}
    We have $H$-equivariant isomorphisms
    \[
    \Ext^2_\Lambda(M_z,M_z)^\vee
\cong\Hom_\Lambda(M_z,M_z^\circ/M_z)
    \]
and
\[
\Ext^i_\Lambda(T_z,M_z')^\vee\cong\Ext^{3-i}_\Lambda(M_z',T_z),
\quad
\Ext^i_\Lambda(T_z,T_z)^\vee\cong\Ext^{3-i}_\Lambda(T_z,T_z)
\]
for $0\leq i\leq3$.
\end{lemma}

\begin{proof}
    Arguments similar to those in the proofs of Proposition~\ref{prop:dual} and Proposition~\ref{prop:ext2-vanish} give the first assertion and
    \[
    \E^3_\Lambda(T_z)^\vee\cong T_z,\quad\E^i_\Lambda(T_z)=0,\quad i=0,1,2.
    \]
    The short exact sequence
    \[
    0\longrightarrow M_z'\longrightarrow M_z\longrightarrow T_z\longrightarrow0
    \]
    implies that the projective dimension of $T_z$ is at most 3. Since $\Lambda$ is local and left-noetherian, $T_z$ admits a free $\Lambda$-module resolution of finite rank
    \[
    0\longrightarrow P_3\longrightarrow P_2\longrightarrow P_1\longrightarrow P_0\longrightarrow T_z\longrightarrow0.
    \]
    Applying $\Hom_\Lambda(\cdot,\Lambda)$, we obtain an exact sequence
    \[
    0\longrightarrow \E^0_\Lambda(P_0)\longrightarrow \E^0_\Lambda(P_1)\longrightarrow \E^0_\Lambda(P_2)\longrightarrow \E^0_\Lambda(P_3)\longrightarrow \E^3_\Lambda(T_z)\longrightarrow0,
    \]
    which is a free right $\Lambda$-module resolution of finite rank of $\E^3_\Lambda(T_z)$. Applying $(\cdot)\otimes_{\Lambda}M_z'$ to the exact sequence and using 
    \[
    \E^0_\Lambda(P_j)\otimes_{\Lambda}M_z'\cong\Hom_\Lambda(P_j,M_z'),
    \]
    we obtain
    \[
    \Ext^i_\Lambda(T_z,M_z')\cong\Tor^\Lambda_{3-i}\left(\E^3_\Lambda(T_z),M_z'\right).
    \]
    Taking Pontryagin duals and using the adjunction formula, we obtain
    \[
     \Ext^i_\Lambda(T_z,M_z')^\vee\cong\Tor^\Lambda_{3-i}\left(\E^3_\Lambda(T_z),M_z'\right)^\vee\cong\Ext^{3-i}_\Lambda(M_z',T_z).
    \]
    The assertion for $\Ext^i_\Lambda(T_z,T_z)$ is similar.
\end{proof}

For two $\Lambda$-modules $M_1$ and $M_2$ with compatible $H$-actions, define
\[
\EP_\Lambda^H(M_1,M_2):=\sum_{i=0}^{\infty}(-1)^i\dim_\FF\Ext^i_\Lambda(M_1,M_2)^H
\]
when the summation makes sense.

\begin{corollary}\label{cor:ep}
    We have
    \[
    \EP_\Lambda^H(T_z,M_z')=-\EP_\Lambda^H(M_z',T_z),\quad\EP_\Lambda^H(T_z,T_z)=0.
    \]
\end{corollary}

\begin{proof}
    This follows immediately by Lemma~\ref{lem:Mz-dual} and the definition of $\EP_\Lambda^H(\cdot,\cdot)$.
\end{proof}

\begin{lemma}\label{lem:T}
    We have
    \[
    \Hom_\Lambda(T_z,M_z^\circ/M_z)^H=0,\quad\dim_\FF\Ext^1_\Lambda(M_z',T_z)^H\leq1.
    \]
\end{lemma}
\begin{proof}
    Put $\ell=|r-\frac{p-1}{2}|$. We can assume $\ell\ne0$. Then
    \[
    T_z=
    \begin{cases}
        \bigoplus_{k=0}^{\ell-1}\FF Z^ke_0, &r>\frac{p-1}{2};\\
        \bigoplus_{k=0}^{\ell-1}\FF Z^ke_1, &r<\frac{p-1}{2}
    \end{cases}
    \]
    and
    \[
    T_z|_{H'}\cong\bigoplus_{j=1}^{\ell}\psi\alpha^j.
    \]
    Since $\alpha|_{H'}$ has order $(p+1)/2$ and
    $\ell<(p+1)/2$, the complex
    $\Hom_\Lambda(P_{N,\bullet},T_z)^{H'}$ has the shape
    \[
    0\longrightarrow\FF\longrightarrow\FF
    \]
    (see \eqref{eq:N+-resolution}). Moreover, by Shapiro's lemma we have
    \[
    \Ext^1_\Lambda(M_z',T_z)^H\cong\Ext^1_\Lambda(N_+,T_z)^{H'}.
    \]
    This gives the second assertion.

    The image of any element in $\Hom_\Lambda(T_z,M_z^\circ/M_z)^H$ is contained in $(M_z^\circ/M_z)[Z^\ell]$, which has basis
    \[
    [Z^{-j}e_0],\ [Z^{-j}e_1],\quad 1\leq j\leq\ell.
    \]
    If $r>(p-1)/2$, the generator of $T_z$ has $H$-character $\chi^\vee$. But there does not exist an integer $1\leq j\leq\ell$ satisfying $j\equiv0\pmod{p+1}$ or $j\equiv r+1\pmod{p+1}$.
    When $r<(p-1)/2$, the generator of $T_z$ has $H$-character
    $(\chi^\sigma\alpha)^\vee$. But there does not exist an integer $1\leq j\leq\ell$ satisfying
    $j\equiv0\pmod{p+1}$ or $j\equiv p-r\pmod{p+1}$. This gives the first assertion.
\end{proof}

\begin{lemma}\label{lem:Hom-1-dimensional}
    We have
    \[
    \dim_\FF\Hom_\Lambda(M_z,M_z^\circ/M_z)^H=1.
    \]
\end{lemma}

\begin{proof}
    If $\ell=|r-\frac{p-1}{2}|=0$, this follows from
    \[
    \Hom_\Lambda(M_z,M_z^\circ/M_z)^H=\Hom_\Lambda(M_z',M_z^\circ/M_z')^H\cong\Ext_\Lambda^2(M_z',M_z')^{\vee,H}
    \]
    and Corollary~\ref{cor:Mz'-ext}. 
    
    Assume $\ell\ne0$. We first prove that 
    \[
    \dim_\FF\Hom_\Lambda(M_z,M_z^\circ/M_z)^H\le1.
    \] 
    Applying $\Hom_\Lambda(\cdot,M_z^\circ/M_z)^H$ to the short exact sequence
    \[
    0\longrightarrow M_z'\longrightarrow M_z\longrightarrow T_z\longrightarrow0,
    \]
    we obtain an injection
    \[
    \Hom_\Lambda(M_z,M_z^\circ/M_z)^H\hookrightarrow\Hom_\Lambda(M_z',M_z^\circ/M_z)^H
    \]
    by Lemma~\ref{lem:T}. Applying $\Hom_\Lambda(M_z',\cdot)^H$ to the short exact sequence
    \[
    0\longrightarrow T_z\longrightarrow M_z^\circ/M_z'
    \longrightarrow M_z^\circ/M_z\longrightarrow0,
    \]
    we obtain an injection
    \[
    \Hom_\Lambda(M_z',M_z^\circ/M_z)^H\hookrightarrow\Ext^1_\Lambda(M_z',T_z)^H
    \]
    because 
    \[
    \Hom_\Lambda(M_z',M_z^\circ/M_z')^H\cong\Ext^2_{\Lambda}(M_z',M_z')^{\vee,H}=0
    \]
    by Corollary~\ref{cor:Mz'-ext}. By Lemma~\ref{lem:T}, we have
    \[
    \dim_\FF\Hom_\Lambda(M_z,M_z^\circ/M_z)^H\le1.
    \]

    Now we construct an explicit nonzero element of $\Hom_\Lambda(M_z,M_z^\circ/M_z)^H$. Define an $\FF[[Z]]$-linear map
    \[
    \varphi:M_z\longrightarrow M_z^\circ/M_z,\quad e_0\mapsto b_0r[Z^{-r-1}e_1],\quad e_1\mapsto b_1(r+1)[Z^{r-p}e_0].
    \]
    We need to show $\varphi$ is also $Y$-linear, i.e.,
    \[
    Y\varphi(Z^ke_i)=\varphi(YZ^ke_i)
    \]
    for any $k\ge0$ and $i=0,1$. Using \eqref{eq:leading-coefficient}, it suffices to show
    \[
    rb_0b_1(-r-1-\beta_1)Z^{-1}e_0-b_0b_1(r+1)(-\beta_0)Z^{-1}e_0\in M_z
    \]
    and
    \[
    b_0b_1(r+1)(r-p-\beta_0)Z^{-1}e_1-rb_0b_1(-\beta_1)Z^{-1}e_1\in M_z,
    \]
    which follow by Proposition~\ref{prop:leading-coefficient}.
\end{proof}

\begin{proposition}\label{prop:z-z-ext}
    For $0\le i\le2$, we have
    \[
    \dim_\FF\Ext^i_{\Lambda}(M_z,M_z)^{H}=1.
    \]
\end{proposition}

\begin{proof}
    The assertion for $i=2$ follows from Lemma~\ref{lem:Mz-dual} and Lemma~\ref{lem:Hom-1-dimensional}. Any element of $\Hom_{\Lambda}(M_z,M_z)^{H}$ preserves the filtration $\cF^nM_z$, and therefore the image of
    \[
    \Hom_{\Lambda}(M_z,M_z)^{H}\longrightarrow\Hom_{\Lambda}(M_z',M_z)^{H}
    \]
    is contained in $\Hom_{\Lambda}(M_z',M_z')^{H}=\FF$. As $T_z$ is $Z$-torsion and $M_z$ is $Z$-torsion-free, the map
    \[
    \Hom_{\Lambda}(M_z,M_z)^{H}\longrightarrow\Hom_{\Lambda}(M_z',M_z')^{H}
    \]
    is injective, which gives the assertion for $i=0$.

    By Corollary~\ref{cor:ep}, we have
    \begin{align*}
    \EP_\Lambda^H(M_z,M_z)&=\EP_\Lambda^H(M_z',M_z')+\EP_\Lambda^H(M_z',T_z)+\EP_\Lambda^H(T_z,M_z')+\EP_\Lambda^H(T_z,T_z)\\
    &=\EP_\Lambda^H(M_z',M_z')=1,
    \end{align*}
    which gives
    \[
    \dim_\FF\Ext^1_{\Lambda}(M_z,M_z)^{H}=\dim_\FF\Ext^2_{\Lambda}(M_z,M_z)^{H}=1.
    \]
\end{proof}

\begin{theorem}\label{thm:self-extension}
    We have
        \[
    \dim_\FF\Ext^i_{D^\times,\zeta}(\pi,\pi)= 
    \begin{cases}
        1,&i=0;\\
        3,&i=1,2;\\
        0,&i\ge3.
    \end{cases}
    \]
\end{theorem}

\begin{proof}
    This follows by Corollary~\ref{cor:Ext-decomposition}, Proposition~\ref{prop:mixed-ext-dim} and Proposition~\ref{prop:z-z-ext}. 
\end{proof}

We give an explicit characterization for $\Ext^1_{D^\times,\zeta}(\pi,\pi)$. Put $M_y^\circ:=\FF((Y))\otimes_{\FF[[Y]]} M_y$ and extend the $\Lambda$-action to $M_y^\circ$ by the $Y$-version of Proposition~\ref{prop:Laurent-extension}.

\begin{proposition}\label{prop:ext-explicit}
    The natural connecting homomorphisms induce isomorphisms
    \[
    \Ext^1_{\Lambda}(M_z,M_y)^H\cong\Hom_{\Lambda}(M_z,M_y^\circ/M_y)^H=\FF\varphi_{1,0}\oplus\FF\varphi_{0,1},
    \]
    \[
    \Ext^1_{\Lambda}(M_z,M_z)^H\cong\Hom_{\Lambda}(M_z,M_z^\circ/M_z)^H=\FF\varphi,
    \]
    where $\varphi$ is defined in Lemma~\ref{lem:Hom-1-dimensional} and $\varphi_{1,0},\varphi_{0,1}$ are determined by
    \[
    \varphi_{1,0}(e_0)=[Y^{-1}\varpi_D(e_1)],\quad\varphi_{1,0}(e_1)=0;\quad\varphi_{0,1}(e_0)=0,\quad\varphi_{0,1}(e_1)=[Y^{-1}\varpi_D(e_0)].
    \]
\end{proposition}

\begin{proof}
    Put $F^nM_y^\circ=Y^nM_y$ for $n\in\ZZ$. This is a complete and exhaustive filtration, and
    \[
    \gr_F(M_y^\circ)\cong\left((\chi^\sigma)^\vee\otimes\FF[y,y^{-1}]\right)\oplus\left((\chi\alpha^{-1})^\vee\otimes\FF[y,y^{-1}]\right).
    \]
    Using the resolution \eqref{eq:resolution-Fz} and the proof of Lemma~\ref{lem:mixed-gr-Ext}, for any $i\ge0$ we have
    \[
    \Ext_{\gr_{\mathfrak{m}_D}(\Lambda)}^i\left(\gr_{\mathfrak{m}_D}(M_z),\gr_F(M_y^\circ)\right)=0.
    \]
    The same argument as in the proof of Proposition~\ref{prop:dual} gives $\Ext^i_{\Lambda}(M_z,M_y^\circ)=0$ for $i\ge0$. Applying $\Hom_{\Lambda}(M_z,\cdot)^H$ to the short exact sequence
    \[
    0\longrightarrow M_y\longrightarrow M_y^\circ\longrightarrow M_y^\circ/M_y\longrightarrow0,
    \]
    we obtain 
    \[   
    \Ext^1_{\Lambda}(M_z,M_y)^H\cong\Hom_{\Lambda}(M_z,M_y^\circ/M_y)^H.
    \]
    By Proposition~\ref{prop:mixed-ext-dim} (conjugated by $\varpi_D$) and the proofs of Proposition~\ref{prop:dual} and Proposition~\ref{prop:ext2-vanish}, we have
    \[
    \dim_\FF\Hom_{\Lambda}(M_z,M_y^\circ/M_y)^H=\dim_\FF\Ext^2_{\Lambda}(M_y,M_z)^H=2.
    \]
    One checks that $\varphi_{1,0},\varphi_{0,1}$ are $\Lambda$-linear and linearly independent. This gives the first assertion.

    The same argument as in the proof of Proposition~\ref{prop:dual} gives $\Ext^i_{\Lambda}(T_z,M_z^\circ)=0$ for $i\ge0$. Therefore we have
    \[
    \Ext^i_{\Lambda}(M_z,M_z^\circ)^H\cong\Ext^i_{\Lambda}(M_z',M_z^\circ)^H\cong\Ext^i_{\Lambda}(N_+,N^\circ_+)^{H'}\oplus\Ext^i_{\Lambda}(N_+,N_-^\circ)^{H'}.
    \]
    Using the resolution \eqref{eq:N+-resolution} and the proof of Lemma~\ref{lem:mixed-N-Ext}, for any $i\ge0$ we have $\Ext_{\widetilde{\Lambda}}^i\left(\widetilde{N_+},\widetilde{N_-^\circ}\right)=0$. The same argument as in the proof of Proposition~\ref{prop:dual} gives $\Ext^i_{\Lambda}(N_+,N_-^\circ)^{H'}=0$ for $i\ge0$. Applying $\Hom_{\Lambda}(N_+,\cdot)^{H'}$ to the short exact sequence
    \[
    0\longrightarrow N_+\longrightarrow N_+^\circ\longrightarrow N_+^\circ/N_+\longrightarrow0,
    \]
    we obtain a long exact sequence
    \[
    \begin{split}
         0\longrightarrow\Hom_{\Lambda}(N_+,N_+)^{H'}\longrightarrow\Hom_{\Lambda}(N_+,N_+^\circ)^{H'}\longrightarrow\Hom_{\Lambda}(N_+,N_+^\circ/N_+)^{H'}\\
         \longrightarrow\Ext^1_{\Lambda}(N_+,N_+)^{H'}\longrightarrow\Ext^1_{\Lambda}(N_+,N_+^\circ)^{H'}\longrightarrow\Ext^1_{\Lambda}(N_+,N_+^\circ/N_+)^{H'}.
    \end{split}
    \]
    By Proposition~\ref{prop:Hom-vanish}, we have $\Hom_{\Lambda}(N_+,N_+^\circ/N_+)^{H'}=0$. By Proposition~\ref{prop:ext-N}, we have $\Hom_{\Lambda}(N_+,N_+)^{H'}=\FF$ and $\Ext^1_{\Lambda}(N_+,N_+)^{H'}=0$. Since $N_+$ is Cohen-Macaulay of grade $2$, the same argument as in the proof of Lemma~\ref{lem:Mz-dual} gives
    \[
    \Ext^1_{\Lambda}(N_+,N_+^\circ/N_+)^{H'}\cong\Tor^\Lambda_1\left(\E^2_\Lambda(N_+),N_+\right)^{H',\vee}\cong\Ext^1_{\Lambda}(N_+,N_+)^{H',\vee}=0.
    \]
    Therefore, we have
    \[
    \Hom_{\Lambda}(N_+,N_+^\circ)^{H'}\cong\Hom_{\Lambda}(N_+,N_+)^{H'}=\FF,\quad\Ext^1_{\Lambda}(N_+,N_+^\circ)^{H'}=0.
    \]
    This gives
    \[
    \Hom_{\Lambda}(M_z,M_z^\circ)^H=\FF,\quad\Ext^1_{\Lambda}(M_z,M_z^\circ)^H=0,
    \]
    Applying $\Hom_{\Lambda}(M_z,\cdot)^H$ to the short exact sequence
    \[
    0\longrightarrow M_z\longrightarrow M_z^\circ\longrightarrow M_z^\circ/M_z\longrightarrow0,
    \]
    we obtain 
    \[   
    \Ext^1_{\Lambda}(M_z,M_z)^H\cong\Hom_{\Lambda}(M_z,M_z^\circ/M_z)^H=\FF\varphi.
    \]
\end{proof}

For $\lambda_0,\lambda_1\in\FF$, put $\varphi_{\lambda_0,\lambda_1}=\lambda_0\varphi_{1,0}+\lambda_1\varphi_{0,1}\in\Hom_{\Lambda}(M_z,M_y^\circ/M_y)^H$. The extension corresponding to $\varphi_{\lambda_0,\lambda_1}$ can be identified with
\[
E_{\lambda_0,\lambda_1}=\{(t_1,t_2)\in M_y^\circ\oplus M_z:\varphi_{\lambda_0,\lambda_1}(t_2)=[t_1]\}.
\]
It splits if and only if $(\lambda_0,\lambda_1)=(0,0)$. Similarly, let
\[
E_\varphi=\{(t_1,t_2)\in M_z^\circ\oplus M_z:\varphi(t_2)=[t_1]\}
\]
be the nonsplit extension corresponding to $\varphi\in\Hom_{\Lambda}(M_z,M_z^\circ/M_z)^H$. By fixing the central character $\zeta$, the Pontryagin dual $E_{\lambda_0,\lambda_1}^\vee$ (resp. $E_\varphi^\vee$) is a smooth $\cO_D^{\times}Z_D$-module. Write
\[
\pi_{\lambda_0,\lambda_1}:=\Ind_{\cO_D^{\times}Z_D}^{D^\times}E_{\lambda_0,\lambda_1}^\vee,\quad\pi_\varphi:=\Ind_{\cO_D^{\times}Z_D}^{D^\times}E_\varphi^\vee.
\]

\begin{proposition}
      We have 
      \[
      \pi_\varphi^{U_D^1}=\pi_\varphi[\mathfrak{m}_D]=\chi\oplus\chi^\sigma\oplus\chi^\sigma\alpha\oplus\chi\alpha^{-1},
      \]
      \[
      \pi_{1,0}[\mathfrak{m}_D]=\chi^{\oplus2}\oplus(\chi^\sigma)^{\oplus2}\oplus\chi^\sigma\alpha\oplus\chi\alpha^{-1},\quad\pi_{0,1}[\mathfrak{m}_D]=\chi\oplus\chi^\sigma\oplus(\chi^\sigma\alpha)^{\oplus2}\oplus(\chi\alpha^{-1})^{\oplus2}.
      \]
\end{proposition}

\begin{proof}
    By the definition of $\varphi$, the quotient $E_\varphi/\mathfrak{m}_DE_\varphi$ factors through
    \[
    \overline{E}:=\left\{(t_1,t_2)\in (M_z^\circ/ZM_z)\oplus\left(\FF[[Z]]e_0/(Z^{r+2})\oplus\FF[[Z]]e_1/(Z^{p-r+1})\right):\varphi(t_2)=[t_1]\right\}.
    \]
    Consider the four elements in $\overline{E}$ 
    \[
    \left(b_0r[Z^{-r}e_1],[Ze_0]\right),\quad\left(b_1(r+1)[Z^{r+1-p}e_0],[Ze_1]\right),
    \]
    \[
    \left(b_0r[Z^{-1}e_1],[Z^re_0]\right),\quad\left(b_1(r+1)[Z^{-1}e_0],[Z^{p-r-1}e_1]\right).
    \]
    By \eqref{eq:leading-coefficient} and Proposition~\ref{prop:leading-coefficient} we have
    \[
    Y\left(b_0r[Z^{-r}e_1],[Ze_0]\right)=\left(-cr(r-1)[e_0],b_0(1-\frac{r}{2})[Z^{p-r}e_1]\right),
    \]
    \[
    Z\left(b_1(r+1)[Z^{-1}e_0],[Z^{p-r-1}e_1]\right)=\left(b_1(r+1)[e_0],[Z^{p-r}e_1]\right),
    \]
    \[
    Y\left(b_1(r+1)[Z^{r+1-p}e_0],[Ze_1]\right)=\left(c(r+1)(r+2)[e_1],b_1\frac{r+3}{2}[Z^{r+1}e_0]\right),
    \]
    \[
    Z\left(b_0r[Z^{-1}e_1],[Z^re_0]\right)=\left(b_0r[e_1],[Z^{r+1}e_0]\right).
    \]
    Since the matrices $
    \begin{pmatrix}
    -cr(r-1) & b_0(1-\frac{r}{2})\\
    b_1(r+1) & 1
    \end{pmatrix} $ and $
    \begin{pmatrix}
    c(r+1)(r+2) & b_1\frac{r+3}{2}\\
    b_0r & 1
    \end{pmatrix}$ are invertible, the quotient $E_\varphi/\mathfrak{m}_DE_\varphi$ further factors through
    \[
    \left\{(t_1,t_2)\in (M_z^\circ/M_z)\oplus\left(\FF[[Z]]e_0/(Z^{r+1})\oplus\FF[[Z]]e_1/(Z^{p-r})\right):\varphi(t_2)=t_1\right\},
    \]
    which is isomorphic to $\FF[[Z]]e_0/(Z^{r+1})\oplus\FF[[Z]]e_1/(Z^{p-r})$. Therefore, we have
    \[
    E_\varphi/\mathfrak{m}_DE_\varphi=\FF[e_0]\oplus\FF[e_1]=\chi^\vee\oplus(\chi^\sigma\alpha)^\vee.
    \]
    A similar argument gives
    \[
    E_{1,0}/\mathfrak{m}_DE_{1,0}=\FF[e_0]\oplus\FF[e_1]\oplus\FF[\varpi_D(e_0)]=\chi^\vee\oplus(\chi^\sigma\alpha)^\vee\oplus(\chi^\sigma)^\vee,
    \]
    \[
    E_{0,1}/\mathfrak{m}_DE_{0,1}=\FF[e_0]\oplus\FF[e_1]\oplus\FF[\varpi_D(e_1)]=\chi^\vee\oplus(\chi^\sigma\alpha)^\vee\oplus(\chi\alpha^{-1})^\vee.
    \]
    The result follows from taking Pontryagin duals.
\end{proof}

\begin{proposition}\label{prop:application}
    Let $\pi_{\lambda_0,\lambda_1}$ be as above. Suppose $(\lambda_0,\lambda_1)\ne(0,0)$. Then
    \[
    \dim_\FF\Ext^1_{D^\times,\zeta}(\pi,\pi_{\lambda_0,\lambda_1})=3.
    \]
\end{proposition}

\begin{proof}
    The same argument as in the proof of Corollary~\ref{cor:Ext-decomposition} gives
    \[
    \Ext_{D^\times,\zeta}^1(\pi,\pi_{\lambda_0,\lambda_1})\cong\Ext^1_{\Lambda}(E_{\lambda_0,\lambda_1},M)^H\cong\Ext^1_{\Lambda}(E_{\lambda_0,\lambda_1},M_z)^H\oplus\Ext^1_{\Lambda}(E_{\lambda_0,\lambda_1},M_y)^H.
    \]
    Write $E=E_{\lambda_0,\lambda_1}$. We first compute $\dim_\FF\Ext^1_{\Lambda}(E,M_z)^H$. By Proposition~\ref{prop:mixed-ext-dim} and Proposition~\ref{prop:z-z-ext}, we have
    \[
    \dim_\FF\Hom_{\Lambda}(E,M_z)^H=\dim_\FF\Hom_{\Lambda}(M_z,M_z)^H=1.
    \]
    Since 
    \[
    \EP_\Lambda^H(E,M_z)=\EP_\Lambda^H(M_z,M_z)+\EP_\Lambda^H(M_y,M_z)=1,
    \]
    it suffices to compute $\dim_\FF\Ext^2_{\Lambda}(E,M_z)^H$. Put $E^\circ=M_y^\circ\oplus M_z^\circ$. Then $E^\circ/E$ can be identified with
    \[
    \left((M_y^\circ/M_y)\oplus M_z^\circ\right)/\{\left(\varphi_{\lambda_0,\lambda_1}(t),t\right):t\in M_z\}.
    \]
    Arguments similar to those in the proofs of Proposition~\ref{prop:dual} and Proposition~\ref{prop:ext2-vanish} give
    \[
    \Ext^2_{\Lambda}(E,M_z)^{\vee,H}\cong\Hom_{\Lambda}(M_z,E^\circ/E)^H.
    \]
    Applying $\Hom_{\Lambda}(M_z,\cdot)^H$ to the short exact sequence
    \[
     0\longrightarrow M_y^\circ/M_y\xrightarrow{\,\,\,\iota\,\,\,} E^\circ/E\longrightarrow M_z^\circ/M_z\longrightarrow0,
    \]
    we obtain an exact sequence
    \begin{equation}\label{eq:left-exact}
        0\longrightarrow\Hom_{\Lambda}(M_z,M_y^\circ/M_y)^H\longrightarrow\Hom_{\Lambda}(M_z,E^\circ/E)^H\longrightarrow\Hom_{\Lambda}(M_z,M_z^\circ/M_z)^H=\FF\varphi.
    \end{equation}

    We claim $\varphi$ cannot be lifted to $\Hom_{\Lambda}(M_z,E^\circ/E)^H$. Otherwise, let $\phi$ be such a lift. Let $\widetilde{\varphi}$ be the element in $\Hom_{\FF[[Z]]}(M_z,M_z^\circ)^H$ sending $e_0$ to $b_0rZ^{-r-1}e_1$ and $e_1$ to $b_1(r+1)Z^{r-p}e_0$, and $s$ be the element in $\Hom_{\FF[[Z]]}(M_z,E^\circ/E)^H$ sending $x$ to $[\left(0,\widetilde{\varphi}(x)\right)]$. Then there exists a unique element $\theta$ in $\Hom_{\FF[[Z]]}(M_z,M_y^\circ/M_y)^H$ such that $\phi=s+\iota_*\theta$. By \eqref{eq:leading-coefficient} and the proof of Lemma~\ref{lem:Hom-1-dimensional}, for $i=0,1$ we have
    \[
    Y\widetilde{\varphi}(e_i)-\widetilde{\varphi}(Ye_i)\in ZM_z,
    \]
    and hence $Ys(e_i)-s(Ye_i)=0$ because by Proposition~\ref{prop:ext-explicit} $\varphi_{\lambda_0,\lambda_1}$ factors through $M_z/ZM_z$. Since $\phi$ is $\FF[[Y]]$-linear, we obtain $Y\theta(e_i)=\theta(Ye_i)$ for $i=0,1$. Let $n$ be the minimal integer such that $Y^n\theta(e_i)=0$ for $i=0,1$. Then $n\le 1$. Otherwise there exists $i\in\{0,1\}$ such that
    \[
    0\ne Y^{n-1}\theta(e_i)=Y^{n-2}\theta(Ye_i),
    \]
    which is a contradiction because $Ye_i\in Z^2M_z$ and $Z$ raises the filtration index of $M_y^\circ/M_y$ by at least $2$. Consequently, we have $Y\theta(Ze_i)=\theta(YZe_i)=0$ for $i=0,1$. Since $\phi$ is $\FF[[Y]]$-linear, we obtain $Ys(Ze_i)=s(YZe_i)$, which implies $\varphi_{\lambda_0,\lambda_1}\left(Y\widetilde{\varphi}(Ze_i)-\widetilde{\varphi}(YZe_i)\right)=0$. By \eqref{eq:leading-coefficient}, we have 
    \[
    Y\widetilde{\varphi}(Ze_0)-\widetilde{\varphi}(YZe_0)\equiv-2ce_0\pmod{ZM_z},\quad Y\widetilde{\varphi}(Ze_1)-\widetilde{\varphi}(YZe_1)\equiv2ce_1\pmod{ZM_z}.
    \]
    Since $c\ne0$, we obtain $\lambda_0=\lambda_1=0$, which is a contradiction. As a result, by \eqref{eq:left-exact} and Proposition~\ref{prop:ext-explicit}, we have
    \[
    \dim_\FF\Ext^2_{\Lambda}(E,M_z)^H=\dim_\FF\Hom_{\Lambda}(M_z,E^\circ/E)^H=\dim_\FF\Hom_{\Lambda}(M_z,M_y^\circ/M_y)^H=2.
    \]

    Now we compute $\dim_\FF\Ext^1_{\Lambda}(E,M_y)^H$. Since $M_z$ and $M_y$ are Cohen-Macaulay of grade $2$ and essentially self-dual by assumption~\ref{it:assum-ii}, applying $\E_\Lambda^2(\cdot)\otimes(\zeta^{-1}\circ\mathrm{Nr})$ to the short exact sequence
    \[
    0\longrightarrow M_y\longrightarrow E\longrightarrow M_z\longrightarrow0
    \] 
    and then conjugating by $\varpi_D$ give a nonsplit extension
    \[
    \quad0\longrightarrow M_y\longrightarrow E^\dagger\longrightarrow M_z\longrightarrow0.
    \]
    The operation also gives an isomorphism
    \[
    \Ext^1_{\Lambda}(E,M_y)^H\cong\Ext^1_{\Lambda}(M_z,E^\dagger)^H.
    \]
    Applying $\Hom_{\Lambda}(M_z,\cdot)^H$ to the short exact sequence
    \[
    0\longrightarrow E^\dagger\longrightarrow E^{\dagger,\circ}\longrightarrow E^{\dagger,\circ}/E^\dagger\longrightarrow0
    \]
    gives a long exact sequence
    \[
    \begin{split}
         0\longrightarrow\Hom_{\Lambda}(M_z,E^\dagger)^H\longrightarrow\Hom_{\Lambda}(M_z,E^{\dagger,\circ})^H\longrightarrow\Hom_{\Lambda}(M_z,E^{\dagger,\circ}/E^\dagger)^H\\
         \longrightarrow\Ext^1_{\Lambda}(M_z,E^\dagger)^H\longrightarrow\Ext^1_{\Lambda}(M_z,E^{\dagger,\circ})^H.
    \end{split}
    \]
    We have $\Hom_{\Lambda}(M_z,E^\dagger)^H=0$ because $E^\dagger$ is nonsplit and $\End_\Lambda(M_z)^H=\FF$. We have proved $\dim_\FF\Hom_{\Lambda}(M_z,E^{\dagger,\circ}/E^\dagger)^H=2$. Since $E^{\dagger,\circ}=M_z^\circ\oplus M_y^\circ$, the proof of Proposition~\ref{prop:ext-explicit} gives 
    \[
    \Hom_{\Lambda}(M_z,E^{\dagger,\circ})^H=\FF,\quad\Ext^1_{\Lambda}(M_z,E^{\dagger,\circ})^H=0,
    \]
    which implies $\dim_\FF\Ext^1_{\Lambda}(M_z,E^\dagger)^H=1$. This gives the result.
\end{proof}

\section{The reducible case}

Let $\pi$ be an admissible smooth representation of $D^{\times}$ over $\FF$ with central character $\zeta$. Put $M=\pi^{\vee}$. We make the following assumptions throughout the section:
\begin{enumerate}[label=(\roman*')]
\item\label{it:assum-i-red} There exist an integer $1\le r\le p-4$, a character $\chi$ of $H$ satisfying $\chi^\sigma=\chi\alpha^{r+2}$, and an isomorphism of graded $\gr_{\mathfrak{m}_D}(\Lambda)$-modules with compatible $H$-action 
\begin{equation}\label{eq:red-gr}
\gr_{\mathfrak{m}_D}(M)\cong\left(\chi^\vee\otimes\FF[y,z]/(yz)\right)\oplus\left((\chi^\sigma)^\vee\otimes\FF[y,z]/(yz)\right),
\end{equation}
where on the right $\varpi_D$ permutes the two direct summands and permutes $y$ and $z$.
\item\label{it:assum-ii-red} There is a $D^{\times}$-equivariant isomorphism of $\Lambda$-modules
\begin{equation}\label{eq:self-duality-red}
\E^2_\Lambda(M)\otimes(\zeta^{-1}\circ\mathrm{Nr})\cong M.
\end{equation}
\end{enumerate}

\subsection{A short exact sequence}

By \eqref{eq:red-gr}, we have
\[
\gr_{\mathfrak m_D}(\mathfrak m_D M)
\cong
\bigoplus_{\tau\in\{\chi,\chi^\sigma\}}
\left(
(\tau^\vee\alpha^{-1}\otimes\FF[z])
\oplus(\tau^\vee\alpha\otimes\FF[y])
\right).
\]
By Lemmas~\ref{lem:decomposition-criterion} and~\ref{lem:mixed-gr-Ext},
there is a unique closed filtered $\Lambda$-module decomposition
\[
\mathfrak m_D M=(\mathfrak m_D M)_z\oplus(\mathfrak m_D M)_y
\]
lifting the two graded summands in the directions $z$ and $y$. By uniqueness, the decomposition is compatible with the $H$-action,
and $\varpi_D$ permutes the two summands. Put
\[
M_z=M/(\mathfrak m_D M)_y,\quad
M_y=M/(\mathfrak m_D M)_z.
\]
Then we have
\[
\gr_{\mathfrak m_D}(M_z)\cong(\chi^\vee\otimes\FF[z])
\oplus((\chi^\sigma)^\vee\otimes\FF[z]),\quad\gr_{\mathfrak m_D}(M_y)\cong(\chi^\vee\otimes\FF[y])
\oplus((\chi^\sigma)^\vee\otimes\FF[y]),
\]
and $\varpi_D$ permutes $M_y$ and $M_z$. Moreover, we have a short exact sequence of $\Lambda$-modules with compatible $D^\times$-action
\begin{equation}\label{eq:reducible-exact-sequence}
    0\longrightarrow M\xrightarrow{\,\,\iota\,\,}M_z\oplus M_y
\xrightarrow{\,\,q\,\,}Q\longrightarrow0,
\end{equation}
where $Q=\chi^\vee\oplus(\chi^\sigma)^\vee$, $\iota$ is the pair of quotient maps and $q(x,y)=\overline{x}-\overline{y}$.

Choose $H$-eigenvectors $e_0,e_1\in M$ lifting its two graded generators,
with
\[
\varpi_D e_0=e_1,\quad
\varpi_D e_1=\zeta^{-1}(p)e_0.
\]
We use the same notation for their images in $M_z$ and $M_y$. Then
\[
M_z=\FF[[Z]]e_0\oplus\FF[[Z]]e_1,\quad
M_y=\FF[[Y]]e_0\oplus\FF[[Y]]e_1,
\]
and the quotient maps identify
\[
(\mathfrak m_D M)_z\cong ZM_z,\quad
(\mathfrak m_D M)_y\cong YM_y.
\]
%If $q_j=q(e_j,0)$, the action on $Q$ is
%\[
%\varpi_D q_0=-q_1,\quad
%\varpi_D q_1=-\zeta^{-1}(p)q_0.
%\]

\begin{lemma}\label{lem:reducible-dual}
The isomorphism \eqref{eq:self-duality-red} identifies $\E^2_\Lambda(M_z)\otimes(\zeta^{-1}\circ\mathrm{Nr})$ with $ZM_z$ and $\E^2_\Lambda(M_y)\otimes(\zeta^{-1}\circ\mathrm{Nr})$ with $YM_y$.
\end{lemma}

\begin{proof}
By Lemma~\ref{lem:y-z-resolution} and the proof of
Lemma~\ref{lem:z-dual-z}, $M_z$ and $M_y$ are Cohen-Macaulay of grade $2$. The proof of Lemma~\ref{lem:Mz-dual} gives $\E^3_\Lambda(Q)^\vee\cong Q$ and $\E^i_\Lambda(Q)=0$ for $i\ne3$. Applying $\E^i_\Lambda(\cdot)$ to the short exact sequence \eqref{eq:reducible-exact-sequence} gives
\[
0\longrightarrow\E^2_\Lambda(M_z)\oplus\E^2_\Lambda(M_y)
\longrightarrow\E^2_\Lambda(M)
\longrightarrow\E^3_\Lambda(Q)\longrightarrow0.
\]
Twist by $\zeta^{-1}\circ\mathrm{Nr}$ and identify the middle term with $M$ via the isomorphism \eqref{eq:self-duality-red}. Since $\E^3_\Lambda(Q)$ is killed by $\mathfrak m_D$, the image of the second arrow contains $\mathfrak m_D M$, which has codimension $2$. Hence the image is precisely $\mathfrak m_D M$. Since $\zeta\circ\mathrm{Nr}|_H=\chi\chi^\sigma$, the calculation
in Lemma~\ref{lem:z-dual-z} gives
\[
\gr_{\mathfrak m_D}
\left(\E^2_\Lambda(M_z)\right)\otimes(\zeta^{-1}\circ\mathrm{Nr})\cong(\chi^\vee\alpha^{-1}\otimes\FF[z])
\oplus\left((\chi^\sigma)^\vee\alpha^{-1}\otimes\FF[z]\right),
\]
\[
\gr_{\mathfrak m_D}
\left(\E^2_\Lambda(M_y)\right)\otimes(\zeta^{-1}\circ\mathrm{Nr})\cong(\chi^\vee\alpha\otimes\FF[y])
\oplus\left((\chi^\sigma)^\vee\alpha\otimes\FF[y]\right).
\]
The result follows from the uniqueness of the decomposition of $\mathfrak m_D M$.
\end{proof}

\subsection{Leading coefficients}
Put $\wet(e_0)=0,~\wet(e_1)=r+1-\frac{p-1}{2}$ and
\[
\cF^nM_z
=Z^{\max\{n,0\}}\FF[[Z]]e_0
\oplus Z^{\max\{n-\wet(e_1),0\}}\FF[[Z]]e_1.
\]
The same character calculation as in \eqref{eq:base-change} and \eqref{eq:leading-coefficient} shows that this filtration
is independent of the eigenbasis and makes $M_z$ a filtered
$(\Lambda,\cF)$-module. The same argument as in the proof of Lemma~\ref{lem:leading-linear} gives
\begin{equation}\label{eq:leading-coefficient-red}
\begin{split}
    &Y(Z^ke_0)\equiv b_0(k-\beta_0)Z^{k+p-r-2}e_1\pmod{Z^{k+p}M_z},\\
    &Y(Z^ke_1)\equiv b_1(k-\beta_1)Z^{k+r+1}e_0\pmod{Z^{k+p}M_z}
\end{split}
\end{equation}
for $k\ge0$, where 
\[
b_0b_1=2c\ne0,\quad\beta_0-\beta_1=r+1-\frac{p-1}{2}
\]
and $\beta_0,\beta_1$ are independent of the choice of eigenbasis.

\begin{proposition}\label{prop:reducible-leading-coefficients}
$\beta_0=\frac{r+2}{2},~\beta_1=-\frac{r+1}{2}$.
\end{proposition}

\begin{proof}
The calculation in the proof of Lemma~\ref{lem:grM_z-dual} (with exponents $p-r-2$ and $r+1$) gives the parameters
\[
\beta_0^{\d}=-\beta_0+r+1,\quad
\beta_1^{\d}=-\beta_1-r-2
\]
for $\E^2_\Lambda(M_z)\otimes(\zeta^{-1}\circ\mathrm{Nr})$. With respect to the basis $(Ze_0,Ze_1)$, the parameters of $ZM_z$ are
$\beta_0-1,\beta_1-1$. Since the parameters are independent of the eigenbasis, by Lemma~\ref{lem:reducible-dual} we have
\[
-\beta_0+r+1=\beta_0-1,\quad
-\beta_1-r-2=\beta_1-1,
\]
which give the result.
\end{proof}

\subsection{Computation of
\texorpdfstring{$\Ext^i_\Lambda(M_z,M_z)^H$}{Exti(Mz,Mz)}}

Put $M_z'=\cF^{\max\{r+1-(p-1)/2,0\}}M_z$. The same argument as in the proof of Proposition~\ref{prop:further-decomposition}
gives a unique closed filtered $\Lambda$-module decomposition
\[
M_z'=N_+\oplus N_-
\]
compatible with the $H'$-action, whose two summands are permuted by any element in $H\setminus H'$. Both summands are free of rank one over $\FF[[Z]]$ and their generators have the same $H'$-character. 
The same calculation as in the proof of Lemma~\ref{lem:mixed-N-Ext} gives
\[
\dim_\FF\Ext^i_\Lambda(N_+,N_-)^{H'}
=\begin{cases}
1,&i=1,2,~r=\frac{p-5}{2}~\text{or}~\frac{p-1}{2};\\
0,&\text{otherwise}.
\end{cases}
\]
The same arguments as in the proofs of Propositions~\ref{prop:Hom-vanish},
\ref{prop:ext2-vanish} and~\ref{prop:ext-N} give
\[
\Hom_\Lambda(N_\pm,N_\pm)^{H'}=\FF,
\quad
\Ext^i_\Lambda(N_\pm,N_\pm)^{H'}=0,\quad i\ge1.
\]
Consequently, we have
\[
\Hom_\Lambda(M_z',M_z')^H=\FF,
\quad \EP_\Lambda^H(M_z',M_z')=1.
\]
Put
\[
M_z^\circ=\FF((Z))\otimes_{\FF[[Z]]}M_z,\quad N_\pm^\circ=\FF((Z))\otimes_{\FF[[Z]]}N_\pm,
\]
on which the actions of $\Lambda$ extend uniquely and \eqref{eq:leading-coefficient-red} remains valid for every $k\in\mathbb Z$ by Proposition~\ref{prop:Laurent-extension}.

\begin{lemma}\label{lem:reducible-principal-parts}
Every $\varphi\in\Hom_\Lambda(M_z,M_z^\circ/M_z)^H$ has the form
\[
\varphi(e_0)=b_0a_0[Z^{-r-2}e_1],\quad
\varphi(e_1)=b_1a_1[Z^{r+1-p}e_0]
\]
for some $a_0,a_1\in\FF$.
\end{lemma}

\begin{proof}
Put $\ell=|r+1-(p-1)/2|$. Let $N=N_+$ or $N_-$ with generator $e$. By Proposition~\ref{prop:Hom-vanish}, applying $\Hom_\Lambda(N,\cdot)^{H'}$ to the short exact sequence
\[
0\longrightarrow Z^{-\ell}N/N\longrightarrow N^\circ/N
\longrightarrow N^\circ/Z^{-\ell}N\longrightarrow0
\]
gives an injection 
\[
\Hom_\Lambda(N,N^\circ/Z^{-\ell}N)^{H'}\hookrightarrow\Ext^1_\Lambda(N,Z^{-\ell}N/N)^{H'}.
\]
By the proof of Lemma~\ref{lem:T}, we have
\[
\dim_\FF\Hom_\Lambda(N,N^\circ/Z^{-\ell}N)^{H'}
\le\dim_\FF\Ext^1_\Lambda(N,Z^{-\ell}N/N)^{H'}\le1.
\]
If $\ell=0$, then $\Hom_\Lambda(N,N^\circ/Z^{-\ell}N)^{H'}=0$. If $\ell>0$, consider the $\FF[[Z]]$-linear map given by
\[
e\longmapsto[Z^{-(p+1)/2}e].
\]
It is $H'$-equivariant and $\Lambda$-linear because 
\[
[Y,Z^{-(p+1)/2}](Z^ke)
\equiv\mp b\frac{p+1}{2}Z^{k-1}e\pmod{Z^{k+(p-1)/2}N_\pm}
\]
for $k\ge0$ by \eqref{eq:leading-coefficient-N}. Therefore it spans the whole space.

For a nonzero $H'$-equivariant $\Lambda$-linear map from
$N_+$ to $N_-^\circ/Z^{-\ell}N_-$, let $m$ be
the smallest index such that the image of $e$ is killed by $Z^m$. By \eqref{eq:leading-coefficient-N}, the leading coefficients of $[Y,Z]$ on $N_+$ and $N_-^\circ$ have opposite signs, so $[Y,Z]$ kills the image of $e$, which implies $m\le\frac{p+1}{2}$. Since the map is $H'$-equivariant, we must have $\frac{p+1}{2}\mid m+\ell$. Since $\ell<\frac{p+1}{2}$, we have $m+\ell=\frac{p+1}{2}$. The same assertion holds after interchanging $N_+$ and $N_-$.

Now consider the composition
\[
M_z'=N_+\oplus N_-\hookrightarrow M_z\xrightarrow{\,\varphi\,}M_z^\circ/M_z\twoheadrightarrow M_z^\circ/Z^{-\ell}M_z'=(N_+^\circ/Z^{-\ell}N_+)\oplus (N_-^\circ/Z^{-\ell}N_-).
\]
By the discussion above, the image of the composition is killed by $Z^{(p+1)/2-\ell}$. Therefore the image of $\varphi$ is killed by $Z^p$ because $\ell\le\frac{p-1}{2}$. This gives the form of $\varphi$ by comparing $H$-characters.
\end{proof}

\begin{proposition}\label{prop:reducible-directional-Ext}
We have
\[
\Hom_\Lambda(M_z,M_z)^H=\FF,
\quad
\Ext^1_\Lambda(M_z,M_z)^H
=\Ext^2_\Lambda(M_z,M_z)^H=0.
\]
The same assertions hold for $M_y$.
\end{proposition}

\begin{proof}
Let $\varphi$ be as in Lemma~\ref{lem:reducible-principal-parts}. Since $\varphi$ is $\Lambda$-linear, we have $Y\varphi(e_i)=\varphi(Ye_i)$. Comparing the coefficients of the leading terms $[Z^{-1}e_0]$ and $[Z^{-1}e_1]$ by \eqref{eq:leading-coefficient-red}, we obtain
\[
\begin{pmatrix}
-r-2-\beta_1&\beta_0\\
\beta_1&r+1-p-\beta_0
\end{pmatrix}
\binom{a_0}{a_1}=0.
\]
By Proposition~\ref{prop:reducible-leading-coefficients}, the matrix is invertible, which implies $\Hom_\Lambda(M_z,M_z^\circ/M_z)^H=0$. By Lemma~\ref{lem:Mz-dual}, we have $\Ext^2_\Lambda(M_z,M_z)^H=0$. The same argument as in the proof of Proposition~\ref{prop:z-z-ext} gives $\Hom_\Lambda(M_z,M_z)^H=\FF$ and 
\[
\EP_\Lambda^H(M_z,M_z)=\EP_\Lambda^H(M_z',M_z')=1,
\]
which implies $\Ext^1_\Lambda(M_z,M_z)^H=0$.
Conjugating by $\varpi_D$ gives the assertions for $M_y$.
\end{proof}

\subsection{The main result}

Put 
\[
M_y^\circ=\FF((Y))\otimes_{\FF[[Y]]}M_y,\quad M^\circ=M_z^\circ\oplus M_y^\circ.
\]
We regard $M$ as a submodule of $M^\circ$ via \eqref{eq:reducible-exact-sequence}. 

\begin{proposition}\label{prop:reducible-Hom-explicit}
    We have $\Hom_\Lambda(M,M^\circ/M)^{H,\varpi_D}
=\FF\gamma\oplus\FF\varphi_{1,0}\oplus\FF\varphi_{0,1}$, where the three elements are determined by
\[
\gamma(e_0)=[(e_0,0)],\quad
\gamma(e_1)=-[(e_1,0)],
\]
\[
\varphi_{1,0}(e_0)
=\left[(b_0Z^{-r-2}e_1,
0)\right],\quad
\varphi_{1,0}(e_1)
=\left[\left(0,\zeta^{-1}(p)b_0Y^{-r-2}e_0\right)\right],
\]
\[
\varphi_{0,1}(e_0)
=\left[\left(0,\zeta(p)b_1Y^{r+1-p}e_1\right)\right],\quad
\varphi_{0,1}(e_1)
=\left[(b_1Z^{r+1-p}e_0,0)\right].
\]
\end{proposition}

\begin{proof}
Identify $\mathfrak m_D M$ with $ZM_z\oplus YM_y$. Then $\gamma$ is the composition
\[
M\longrightarrow M/\mathfrak m_D M
\cong Q\hookrightarrow M^\circ/M.
\]
Define $H$-equivariant $\FF[[Z]]$-linear maps
$\varphi_{1,0},\varphi_{0,1}:ZM_z\longrightarrow M^\circ/M$ by
\[
\varphi_{1,0}(Ze_0)=\left[(b_0Z^{-r-1}e_1,-c(r+1)Y^{-1}e_0)\right],\quad\varphi_{1,0}(Ze_1)=\left[(0,-c(r+3)Y^{-1}e_1)\right],
\]
\[
\varphi_{0,1}(Ze_0)=\left[(0,crY^{-1}e_0)\right],\quad\varphi_{0,1}(Ze_1)=\left[(b_1Z^{r+2-p}e_0,c(r+2)Y^{-1}e_1)\right].
\]
By \eqref{eq:leading-coefficient-red} and Proposition~\ref{prop:reducible-leading-coefficients}, we have
\[
Y\varphi_{1,0}(Ze_0)-\varphi_{1,0}(YZe_0)=\left[(-c(r+1)e_0,-c(r+1)e_0)\right]=0,
\]
\[
Y\varphi_{1,0}(Ze_1)-\varphi_{1,0}(YZe_1)=\left[(-c(r+3)e_1,-c(r+3)e_1)\right]=0,
\]
\[
Y\varphi_{0,1}(Ze_0)-\varphi_{0,1}(YZe_0)=\left[(cr e_0,cr e_0)\right]=0,
\]
\[
Y\varphi_{0,1}(Ze_1)-\varphi_{0,1}(YZe_1)=\left[(c(r+2)e_1,c(r+2)e_1)\right]=0,
\]
which imply that both $\varphi_{1,0}$ and $\varphi_{0,1}$ are $\Lambda$-linear.
Extend $\varphi_{1,0}$ and $\varphi_{0,1}$ to $\mathfrak m_D M=ZM_z\oplus YM_y$ via $\varpi_D$-equivariance. Explicitly, on $YM_y$ we have 
\[
\varphi_{1,0}(Ye_0)=\left[(-c(r+3)Z^{-1}e_0,0)\right],\quad\varphi_{1,0}(Ye_1)=\left[(-c(r+1)Z^{-1}e_1,\zeta^{-1}(p)b_0Y^{-r-1}e_0)\right],
\]
\[
\varphi_{0,1}(Ye_0)=\left[(c(r+2)Z^{-1}e_0,\zeta(p)b_1Y^{r+2-p}e_1)\right],\quad\varphi_{0,1}(Ye_1)=\left[(crZ^{-1}e_1,0)\right].
\]
Then $\varphi_{1,0}$ and $\varphi_{0,1}$ are still $\Lambda$-linear. Identify $M$ with its image in $M_z\oplus M_y$ via $\iota$ in \eqref{eq:reducible-exact-sequence}. Extend $\varphi_{1,0}$ and $\varphi_{0,1}$ to $M$ by putting
\[
\varphi_{1,0}(e_0)
=\left[(b_0Z^{-r-2}e_1,
0)\right],\quad
\varphi_{1,0}(e_1)
=\left[\left(0,\zeta^{-1}(p)b_0Y^{-r-2}e_0\right)\right],
\]
\[
\varphi_{0,1}(e_0)
=\left[\left(0,\zeta(p)b_1Y^{r+1-p}e_1\right)\right],\quad
\varphi_{0,1}(e_1)
=\left[(b_1Z^{r+1-p}e_0,0)\right].
\]
Then $\varphi_{1,0}$ and $\varphi_{0,1}$ are $H$-equivariant and $\varpi_D$-equivariant. Again by \eqref{eq:leading-coefficient-red} and Proposition~\ref{prop:reducible-leading-coefficients}, we have
\[
\varphi_{1,0}(Ye_0)=\varphi_{1,0}(-b_0\beta_0Z^{p-r-2}e_1,Ye_0)=\left[\left(-c(r+3)Z^{-1}e_0,0\right)\right]=Y\varphi_{1,0}(e_0),
\]
\[
\varphi_{1,0}(Ye_1)=\varphi_{1,0}(-b_1\beta_1Z^{r+1}e_0,Ye_1)=\left[\left(0,\zeta^{-1}(p)b_0Y^{-r-1}e_0\right)\right]=Y\varphi_{1,0}(e_1),
\]
\[
\varphi_{1,0}(Ze_0)=\zeta(p)\varpi_D\varphi_{1,0}(Ye_1)=\left[(b_0Z^{-r-1}e_1,0)\right]=Z\varphi_{1,0}(e_0),
\]
\[
\varphi_{1,0}(Ze_1)=\varpi_D\varphi_{1,0}(Ye_0)=\left[\left(0,-c(r+3)Y^{-1}e_1\right)\right]=Z\varphi_{1,0}(e_1),
\]
which imply that $\varphi_{1,0}$ is $\Lambda$-linear. Similarly, $\varphi_{0,1}$ is also $\Lambda$-linear.

It remains to show $\dim_\FF\Hom_\Lambda(M,M^\circ/M)^{H,\varpi_D}\le3$. Applying $\Hom_{\Lambda}(M,\cdot)^{H,\varpi_D}$ to
\[
0\longrightarrow Q\longrightarrow M^\circ/M
\longrightarrow (M_z^\circ/M_z)\oplus(M_y^\circ/M_y)
\longrightarrow0,
\]
we obtain
\[
0\longrightarrow\Hom_{\Lambda}(M,Q)^{H,\varpi_D}\longrightarrow\Hom_{\Lambda}(M,M^\circ/M)^{H,\varpi_D}\longrightarrow\Hom_{\Lambda}(M,M_z^\circ/M_z)^H.
\]
Any $\Lambda$-linear map from $M$ to $Q$ factors through $M/\mathfrak m_D M
\cong Q$, and therefore the $H$-equivariance and $\varpi_D$-equivariance imply $\Hom_{\Lambda}(M,Q)^{H,\varpi_D}=\FF\gamma$. Apply $\Hom_{\Lambda}(\cdot,M_z^\circ/M_z)^H$ to
\[
0\longrightarrow(\mathfrak m_D M)_y=YM_y\longrightarrow M\longrightarrow M_z\longrightarrow0.
\]
By Proposition~\ref{prop:reducible-directional-Ext} and Lemma~\ref{lem:Mz-dual}, we have
\[
\Hom_\Lambda(M_z,M_z^\circ/M_z)^H\cong\Ext^2_\Lambda(M_z,M_z)^{\vee,H}=0,
\]
\[
\Ext^1_\Lambda(M_z,M_z^\circ/M_z)^H\cong\Ext^1_\Lambda(M_z,M_z)^{\vee,H}=0.
\]
The same argument as in the proof of Proposition~\ref{prop:mixed-ext-dim} gives
\[
\Hom_\Lambda(M,M_z^\circ/M_z)^H
\cong
\Hom_\Lambda(YM_y,M_z^\circ/M_z)^H\cong\Ext^2_\Lambda(M_z,YM_y)^{\vee,H}=\FF^2.
\]
Therefore $\dim_\FF\Hom_\Lambda(M,M^\circ/M)^{H,\varpi_D}\le3$, which implies the result.
\end{proof}

\begin{theorem}\label{thm:reducible-self-extensions}
We have
\[
\dim_\FF\Ext^i_{D^\times,\zeta}(\pi,\pi)
=\begin{cases}
1,&i=0,\\
3,&i=1,2,\\
0,&i\ge3.
\end{cases}
\]
\end{theorem}

\begin{proof}
    We first prove $\Hom_\Lambda(M,M_z\oplus M_y)^{H,\varpi_D}=\FF\iota$, where $\iota$ is the map in \eqref{eq:reducible-exact-sequence}. In fact, by Frobenius reciprocity we have
    \[
    \Hom_\Lambda(M,M_z\oplus M_y)^{H,\varpi_D}\cong\Hom_\Lambda(M,M_z)^H.
    \]
    By Proposition~\ref{prop:mixed-ext-dim}, any $\Lambda$-linear map from $M$ to $M_z$ factors through $M/(\mathfrak m_DM)_y=M_z$. The assertion follows by Proposition~\ref{prop:reducible-directional-Ext}. This gives the result for $i=0$. The argument of Lemma~\ref{lem:Mz-dual}, together with Proposition~\ref{prop:reducible-Hom-explicit}, gives
    \[
    \dim_\FF\Ext^2_{\Lambda}(M,M)^{H,\varpi_D}=\dim_\FF\Hom_\Lambda(M,M^\circ/M)^{H,\varpi_D}=3.
    \]
    Applying $\Hom_\Lambda(M,\cdot)^{H,\varpi_D}$ to
    \[
    0\longrightarrow M\longrightarrow M^\circ\longrightarrow M^\circ/M\longrightarrow0,
    \]
    we obtain an exact sequence
    \[
    \begin{split}
        0\longrightarrow\FF\longrightarrow\Hom_\Lambda(M,M^\circ)^{H,\varpi_D}\longrightarrow\Hom_\Lambda(M,M^\circ/M)^{H,\varpi_D}\\
        \longrightarrow\Ext^1_{\Lambda}(M,M)^{H,\varpi_D}\longrightarrow\Ext^1_{\Lambda}(M,M^\circ)^{H,\varpi_D}.
    \end{split}
    \]
    The same argument as in the proof of Proposition~\ref{prop:dual} gives $\Ext^i_{\Lambda}(Q,M^\circ)=0$ for $i\ge0$. Therefore, we have
    \[
    \Ext^i_{\Lambda}(M,M^\circ)^{H,\varpi_D}\cong\Ext^i_{\Lambda}(M_z\oplus M_y,M^\circ)^{H,\varpi_D}\cong\Ext^i_{\Lambda}(M_z,M_z^\circ\oplus M_y^\circ)^H.
    \]
    The same argument as in the proof of Proposition~\ref{prop:ext-explicit} gives
    \[
    \Hom_\Lambda(M,M^\circ)^{H,\varpi_D}=\FF,\quad\Ext^1_{\Lambda}(M,M^\circ)^{H,\varpi_D}=0.
    \]
    Therefore, we have
    \[
    \Ext^1_{\Lambda}(M,M)^{H,\varpi_D}\cong\Hom_\Lambda(M,M^\circ/M)^{H,\varpi_D}.
    \]
    In particular, $\dim_\FF\Ext^1_{\Lambda}(M,M)^{H,\varpi_D}=3$. The vanishing for $i\ge3$ follows from $\mathrm{pd}_\Lambda(M_y)=\mathrm{pd}_\Lambda(M_z)=2$, $\mathrm{pd}_\Lambda(Q)=3$ and \eqref{eq:reducible-exact-sequence}.
\end{proof}

Recall $\Hom_\Lambda(M,M^\circ/M)^{H,\varpi_D}=\FF\gamma\oplus\FF\varphi_{1,0}\oplus\FF\varphi_{0,1}$. For $\lambda_0,\lambda_1\in\FF$, put $\varphi_{\lambda_0,\lambda_1}=\lambda_0\varphi_{1,0}+\lambda_1\varphi_{0,1}$ and
\[
E_{\lambda_0,\lambda_1}=\{(t_1,t_2)\in M^\circ\oplus M:\varphi_{\lambda_0,\lambda_1}(t_2)=[t_1]\}.
\]
It is an extension of $M$ by itself, split if and only if $(\lambda_0,\lambda_1)=(0,0)$. Similarly, put
\[
E_\gamma=\{(t_1,t_2)\in M^\circ\oplus M:\gamma(t_2)=[t_1]\}.
\]
Put $\pi_{\lambda_0,\lambda_1}=E_{\lambda_0,\lambda_1}^\vee,~\pi_\gamma=E_\gamma^\vee$.

\begin{proposition}
    We have
    \[
    \pi_{1,0}[\mathfrak m_D]=\pi_{0,1}[\mathfrak m_D]=\chi\oplus\chi^\sigma,\quad\pi_\gamma[\mathfrak m_D]
    =\chi^{\oplus2}\oplus(\chi^\sigma)^{\oplus2}.
    \]
\end{proposition}

\begin{proof}
    Identify $M$ with its image in $M_z\oplus M_y$
    via $\iota$ in \eqref{eq:reducible-exact-sequence}.
    By the proof of Proposition~\ref{prop:reducible-Hom-explicit},
    we have $\varphi_{1,0}(Z^pM_z\oplus Y^pM_y)=0$, and therefore the quotient $E_{1,0}/\mathfrak m_DE_{1,0}$ factors through
    \[
    \overline{E}:=
    \left\{
    (t_1,t_2)\in
    (M^\circ/\mathfrak m_DM)
    \oplus
    \left(M/(Z^{p+1}M_z\oplus Y^{p+1}M_y)\right):
    \varphi_{1,0}(t_2)=[t_1]
    \right\}.
    \]
    Consider the two elements in $\overline{E}$
    \[
    \left(\left[\left(b_0Z^{-r-1}e_1,-c(r+1)Y^{-1}e_0\right)\right],[(Ze_0,0)]\right),\quad\left(0,[(Z^{p-r-2}e_1,0)]\right).
    \]
    By \eqref{eq:leading-coefficient-red} and
    Proposition~\ref{prop:reducible-leading-coefficients},
    we have
    \[
    Y\left(\left[\left(b_0Z^{-r-1}e_1,-c(r+1)Y^{-1}e_0\right)\right],[(Ze_0,0)]\right)=\left(-c(r+1)[(e_0,e_0)],-b_0\frac{r}{2}[(Z^{p-r-1}e_1,0)]\right),
    \]
    \[
    Z\left(0,[(Z^{p-r-2}e_1,0)]\right)
    =
    \left(0,[(Z^{p-r-1}e_1,0)]\right).
    \]
    Since the matrix $
    \begin{pmatrix}
    -c(r+1)&-\frac{b_0r}{2}\\
    0&1
    \end{pmatrix}$
    is invertible, the class of
    $\left([(e_0,e_0)],0\right)$ is zero in
    $E_{1,0}/\mathfrak{m}_DE_{1,0}$.
    The class of $\left([(e_1,e_1)],0\right)$ is also zero by the $\varpi_D$-action. Consequently, the quotient $E_{1,0}/\mathfrak m_DE_{1,0}$ factors through $E_{1,0}/M\cong M$ and
    \[
    E_{1,0}/\mathfrak{m}_DE_{1,0}\cong M/\mathfrak{m}_DM=\FF[e_0]\oplus\FF[e_1]=\chi^\vee\oplus(\chi^\sigma)^\vee.
    \]
    A similar argument gives
    \[
    E_{0,1}/\mathfrak{m}_DE_{0,1}=\chi^\vee\oplus(\chi^\sigma)^\vee,\quad E_\gamma/\mathfrak{m}_DE_\gamma=(\chi^\vee)^{\oplus2}
    \oplus\left((\chi^\sigma)^\vee\right)^{\oplus2}.
    \]
    The result follows from taking Pontryagin duals.
\end{proof}

\begin{proposition}\label{prop:reducible-application}
    Let $\pi_{\lambda_0,\lambda_1}$ be as above. Suppose $(\lambda_0,\lambda_1)\ne(0,0)$. Then
    \[
    \dim_\FF\Ext^1_{D^\times,\zeta}(\pi,\pi_{\lambda_0,\lambda_1})\le3.
    \]
\end{proposition}

\begin{proof}
    By Pontryagin duality, Euler-Poincar\'e characteristic and $\Hom_\Lambda(E_{\lambda_0,\lambda_1},M)^{H,\varpi_D}=\FF$, it suffices to prove
    \[
    \dim_\FF\Ext^2_\Lambda(E_{\lambda_0,\lambda_1},M)^{H,\varpi_D}=\dim_\FF\Hom_\Lambda\left(M,(M^\circ\oplus M^\circ)/E_{\lambda_0,\lambda_1}\right)^{H,\varpi_D}\le4.
    \]
    As in the proof of Proposition~\ref{prop:application}, we have an exact sequence
    \[
    \begin{split}
        0\longrightarrow\Hom_\Lambda(M,M^\circ/M)^{H,\varpi_D}\longrightarrow\Hom_\Lambda\left(M,(M^\circ\oplus M^\circ)/E_{\lambda_0,\lambda_1}\right)^{H,\varpi_D}\\
        \longrightarrow\Hom_\Lambda(M,M^\circ/M)^{H,\varpi_D}.
    \end{split}
    \]
    We need to show at most $1$-dimension of $\Hom_\Lambda(M,M^\circ/M)^{H,\varpi_D}$ on the right can be lifted to the middle term.

    Suppose that 
    \[
    \varphi=a\gamma+a_0\varphi_{1,0}+a_1\varphi_{0,1},\quad a,a_0,a_1\in\FF
    \]
    has a lift $\phi$ to $\Hom_\Lambda\left(M,(M^\circ\oplus M^\circ)/E_{\lambda_0,\lambda_1}\right)^{H,\varpi_D}$. We claim $a_1\lambda_0=a_0\lambda_1$. The argument is similar to the proof of Proposition~\ref{prop:application}. Consider the pushout diagram
    \[
    \begin{tikzcd}
	0 & {M^\circ/M} & {(M^\circ\oplus M^\circ)/E_{\lambda_0,\lambda_1}} & {M^\circ/M} & 0 \\
	0 & {M_z^\circ/M_z} & {\mathcal{L}_{\lambda_0,\lambda_1}} & {M^\circ/M} & 0
	\arrow[from=1-1, to=1-2]
	\arrow[from=1-2, to=1-3]
	\arrow["{p_z}", from=1-2, to=2-2]
	\arrow[from=1-3, to=1-4]
	\arrow[from=1-3, to=2-3]
	\arrow[from=1-4, to=1-5]
	\arrow[equals, nfold, from=1-4, to=2-4]
	\arrow[from=2-1, to=2-2]
	\arrow["{\iota_z}", from=2-2, to=2-3]
	\arrow[from=2-3, to=2-4]
	\arrow[from=2-4, to=2-5]
    \end{tikzcd}
    \]
    where 
    \[
    \cL_{\lambda_0,\lambda_1}=\left((M_z^\circ/M_z)\oplus M^\circ\right)/\{\left(p_z\circ\varphi_{\lambda_0,\lambda_1}(t),t\right):t\in M\}.
    \]
    The lift $\phi$ induces a lift of $\varphi|_{YM_y}$ to $\Hom_\Lambda(YM_y,\cL_{\lambda_0,\lambda_1})^H$, and we still denote it by $\phi$. Put
    \[
    u_0=c\left(-a_0(r+3)+a_1(r+2)\right),\quad u_1=c\left(-a_0(r+1)+a_1r\right)
    \]
    and let $\widetilde{\varphi}$ be the element in $\Hom_{\FF[[Y]]}(YM_y,M^\circ)^H$ sending 
    \[
    Ye_0\mapsto\left(u_0Z^{-1}e_0,\zeta(p)a_1b_1Y^{r+2-p}e_1\right),\quad Ye_1\mapsto\left(u_1Z^{-1}e_1,\zeta^{-1}(p)a_0b_0Y^{-r-1}e_0\right).
    \]
    Let $s$ be the element in $\Hom_{\FF[[Y]]}(YM_y,\cL_{\lambda_0,\lambda_1})^H$ sending $v$ to $[\left(0,\widetilde{\varphi}(v)\right)]$. Then there exists a unique element $\theta$ in $\Hom_{\FF[[Y]]}(YM_y,M_z^\circ/M_z)^H$ such that $\phi=s+\iota_{z,*}\theta$. Since $\phi$ is $\Lambda$-linear, for any $X\in\Lambda$ and $v\in YM_y$ we have
    \begin{equation}\label{eq:Lambda-linear}
        X\theta(v)-\theta(Xv)=p_z\circ\varphi_{\lambda_0,\lambda_1}\left(X\widetilde{\varphi}(v)-\widetilde{\varphi}(Xv)\right).
    \end{equation}
    Using \eqref{eq:leading-coefficient-red}, we have
    \[
    Z\widetilde{\varphi}(Ye_0)\in\left(u_0e_0,ca_1(r+2)e_0+\FF Y^{r+2}e_1+Y^{p+1}M_y\right),
    \]
    \[
    \widetilde{\varphi}(ZYe_0)\in\left(\zeta(p)b_1u_1(1-\beta_1)Y^{r+1}Z^{-1}e_1+Z^{p-1}M_z,ca_0(r+3)e_0+\FF Y^{r+2}e_1+Y^{p+1}M_y\right),
    \]
    \[
    Z\widetilde{\varphi}(Ye_1)\in\left(u_1e_1,-ca_0(r+1)e_1+\FF Y^{p-r-1}e_0+Y^{p+1}M_y\right),
    \]
    \[
    \widetilde{\varphi}(ZYe_1)\in\left(\zeta^{-1}(p)b_0u_0(1-\beta_0)Y^{p-r-2}Z^{-1}e_0+Z^{p-1}M_z,-ca_1re_1+\FF Y^{p-r-1}e_0+Y^{p+1}M_y\right).
    \]
    Therefore, by the proof of Proposition~\ref{prop:reducible-Hom-explicit}, we have
    \[
    Z\theta(Ye_0)-\theta(ZYe_0)=p_z\circ\varphi_{\lambda_0,\lambda_1}(u_0e_0)=\lambda_0u_0b_0[Z^{-r-2}e_1],
    \]
    \[
    Z\theta(Ye_1)-\theta(ZYe_1)=p_z\circ\varphi_{\lambda_0,\lambda_1}(u_1e_1)=\lambda_1u_1b_1[Z^{r+1-p}e_0].
    \]
    Recall by assumption~\ref{it:assum-i-red} we have $1\le r\le p-4$, so both $[Z^{-r-2}e_1]$ and $[Z^{r+1-p}e_0]$ are killed by $Z^{p-2}$. Let $n$ be the minimal integer such that $Z^n\theta(Ye_i)=0$ for $i=0,1$. Then $n\le p-1$. Otherwise there exists $i\in\{0,1\}$ such that
    \[
    0\ne Z^{n-1}\theta(Ye_i)=Z^{n-2}\theta(ZYe_i),
    \]
    which is a contradiction because $ZYe_i\in Y^3M_y$ and $Y$ raises the filtration index of $M_z^\circ/M_z$ by at least $2$. Consequently, by \eqref{eq:leading-coefficient-red} for $i=0,1$ we have $Y^2\theta(Ye_i)=0$, which implies $\theta(ZYe_i)=\theta([Y,Z]Ye_i)=0$. By comparing $H$-characters, we obtain
    \[
    \theta(Ye_0)=\lambda_0u_0b_0[Z^{-r-3}e_1]+\omega_0[Z^{-1}e_0],\quad\theta(Ye_1)=\lambda_1u_1b_1[Z^{r-p}e_0]+\omega_1[Z^{-1}e_1]
    \]
    for some $\omega_0,\omega_1\in\FF$. By \eqref{eq:leading-coefficient-red}, we have
    \[
    [Y,Z]\theta(Ye_0)=2\lambda_0u_0c[Z^{-1}e_0],\quad[Y,Z]\theta(Ye_1)=2\lambda_1u_1c[Z^{-1}e_1]
    \]
    and
    \[
    [Y,Z]\widetilde{\varphi}(Ye_0)\in(u_0b_0Z^{p-r-2}e_1+Z^pM_z,-2ca_1Ye_0+Y^{r+3}M_y),
    \]
    \[
    \widetilde{\varphi}\left([Y,Z]Ye_0\right)\in\left(Z^pM_z,-2ca_0Ye_0+Y^{r+3}M_y\right),
    \]
    \[
    [Y,Z]\widetilde{\varphi}(Ye_1)\in\left(u_1b_1Z^{r+1}e_0+Z^pM_z,-2ca_0Ye_1+Y^{p-r}M_y\right),
    \]
    \[
    \widetilde{\varphi}\left([Y,Z]Ye_1\right)\in\left(Z^pM_z,-2ca_1Ye_1+Y^{p-r}M_y\right).
    \]
    By \eqref{eq:Lambda-linear} and the proof of Proposition~\ref{prop:reducible-Hom-explicit}, we also have
    \[
    [Y,Z]\theta(Ye_0)=-2\lambda_0c^2(r+3)(a_0-a_1)[Z^{-1}e_0]-2\lambda_1c^2a_0[Z^{-1}e_0],
    \]
    \[
    [Y,Z]\theta(Ye_1)=-2\lambda_1c^2r(a_0-a_1)[Z^{-1}e_1]-2\lambda_0c^2a_1[Z^{-1}e_1].
    \]
    Expanding and simplifying, we obtain $a_1\lambda_0=a_0\lambda_1$.

    Finally, we show that $\gamma$ cannot be lifted to $\Hom_\Lambda\left(M,(M^\circ\oplus M^\circ)/E_{\lambda_0,\lambda_1}\right)^{H,\varpi_D}$, which implies the whole result. Suppose otherwise, and let $\widetilde{\gamma}$ be such a lift. By the proof of Proposition~\ref{prop:reducible-Hom-explicit}, $\widetilde{\gamma}|_{\mathfrak m_D M}$ takes values in $M^\circ/M$ and 
    \[
    \Hom_\Lambda(YM_y,M_z^\circ/M_z)^H=\FF\left(Ye_0\mapsto[Z^{-1}e_0]\right)\oplus\FF\left(Ye_1\mapsto[Z^{-1}e_1]\right).
    \]
    Since
    \[
    p_z\circ\varphi_{1,0}|_{YM_y}(Ye_0)=-c(r+3)[Z^{-1}e_0],\quad p_z\circ\varphi_{1,0}|_{YM_y}(Ye_1)=-c(r+1)[Z^{-1}e_1],
    \]
    \[
    p_z\circ\varphi_{0,1}|_{YM_y}(Ye_0)=c(r+2)[Z^{-1}e_0],\quad p_z\circ\varphi_{0,1}|_{YM_y}(Ye_1)=cr[Z^{-1}e_1],
    \]
    and the matrix $\begin{pmatrix}
    -c(r+3) & -c(r+1)\\
    c(r+2) & cr
    \end{pmatrix}$ is invertible, we can assume $p_z\circ\widetilde{\gamma}|_{YM_y}=0$ by adding some linear combination of $\varphi_{1,0}$ and $\varphi_{0,1}$ to $\widetilde{\gamma}$. Then $\widetilde{\gamma}|_{YM_y}$ takes values in $\ker p_z=M_y^\circ/YM_y$. The same arguments as in the proofs of Lemma~\ref{lem:reducible-principal-parts} and Proposition~\ref{prop:reducible-directional-Ext} give
    \[
    \Hom_\Lambda(YM_y,M_y^\circ/YM_y)^H\cong\Hom_\Lambda(ZM_z,M_z^\circ/ZM_z)^H=0.
    \]
    Here the parameters of $ZM_z$ are $\beta_0-1$ and $\beta_1-1$, and the matrix in the proof of Proposition~\ref{prop:reducible-directional-Ext} is still invertible. Therefore $\widetilde{\gamma}|_{YM_y}=0$. Conjugating by $\varpi_D$ gives $\widetilde{\gamma}|_{ZM_z}=0$.

    Identify $(M^\circ\oplus M^\circ)/E_{\lambda_0,\lambda_1}$ with 
    \[
    \left((M^\circ/M)\oplus M^\circ\right)/\{\left(\varphi_{\lambda_0,\lambda_1}(t),t\right):t\in M\}
    \]
    and write
    \[
    \widetilde{\gamma}(e_0)=[\left(u,(e_0,0)\right)],\quad u=[(u_z,u_y)]\in M^\circ/M.
    \]
    Using \eqref{eq:leading-coefficient-red} and the proof of Proposition~\ref{prop:reducible-Hom-explicit}, we have
    \[
    Zu=\varphi_{\lambda_0,\lambda_1}(Ze_0,0)=\left[\left(\lambda_0b_0Z^{-r-1}e_1,c\left(-\lambda_0(r+1)+\lambda_1r\right)Y^{-1}e_0\right)\right],
    \]
    \[
    Yu=\varphi_{\lambda_0,\lambda_1}(Ye_0,0)=\left[\left(-\lambda_1c(r+2)Z^{-1}e_0,0\right)\right].
    \]
    Projecting to $M_y^\circ/M_y$, we obtain
    \[
    Z[u_y]=Y[u_y]=0\implies\lambda_0(r+1)=\lambda_1r.
    \]
    Projecting to $M_z^\circ/M_z$, we obtain
    \[
    Y[u_z]=-\lambda_1c(r+2)[Z^{-1}e_0],\quad Z[u_z]=\lambda_0b_0[Z^{-r-1}e_1]\implies[u_z]=\lambda_0b_0[Z^{-r-2}e_1].
    \]
    But by \eqref{eq:leading-coefficient-red} and Proposition~\ref{prop:reducible-leading-coefficients}, we have
    \[
    Y[u_z]=\lambda_0b_0b_1(-r-2-\beta_1)[Z^{-1}e_0]=-\lambda_0c(r+3)[Z^{-1}e_0].
    \]
    Therefore $\lambda_1(r+2)=\lambda_0(r+3)$, which implies $\lambda_0=\lambda_1=0$. This is a contradiction.
\end{proof}

\section{Deformations and a commutativity criterion}\label{sec:deformations}

Let $\pi$ be an admissible smooth representation of $D^{\times}$ over $\FF$ with central character $\zeta$. Put $M=\pi^{\vee}$. In this section, we assume that $\pi$ satisfies either assumptions~\ref{it:assum-i}\ref{it:assum-ii} or assumptions~\ref{it:assum-i-red}\ref{it:assum-ii-red}. Following \cite{Pas13}, we prove a criterion for the commutativity of the universal deformation ring of $\pi$.

Let $\cO$ be the ring of integers of a finite extension of $\Qp$, with uniformizer $\varpi$ and residue field $\FF$. Fix a continuous character
\[
\widetilde\zeta:Z_D\longrightarrow\cO^\times
\]
lifting $\zeta$. Let $\Art_{\cO}$ be the category of local (possibly noncommutative) artinian $\cO$-algebras with residue field $\FF$ and central $\cO$-action. The morphisms of $\Art_{\cO}$ are local $\cO$-algebra homomorphisms inducing the identity on $\FF$. Following \cite[Definition 3.21]{Pas13}, a deformation of $M$ to $A\in\Art_{\cO}$ with central character $\widetilde{\zeta}^{-1}$ is a pair $(M_A,\iota_A)$, where $M_A$ is a flat compact left $A$-module with a continuous $A$-linear $D^\times$-action such that $Z_D$ acts via $\widetilde\zeta^{-1}$, and
\[
\iota_A:\FF\widehat\otimes_A M_A\xrightarrow{\sim}M
\]
is a $D^\times$-equivariant isomorphism. An isomorphism $(M_A,\iota_A)\cong(M'_A,\iota'_A)$ is a continuous $A[D^\times]$-linear isomorphism $f:M_A\to M'_A$ satisfying
\[
\iota'_A\circ(\id_\FF\widehat\otimes_A f)=\iota_A.
\]
Write $\operatorname{Def}_{\pi,\widetilde\zeta}(A)$ for the set of isomorphism classes of deformations of $M$ to $A$.

Let $\cC_M(\cO)$ be the full subcategory of compact $\cO$-modules with continuous $D^\times$-action, central character $\widetilde\zeta^{-1}$ and a finite filtration with successive quotients isomorphic to $M$. Let $\cC_M(\FF)$ be the full subcategory killed by $\varpi$. By Theorem~\ref{thm:self-extension} and Theorem~\ref{thm:reducible-self-extensions}, we have $\End_{D^\times}(M)=\FF$. Hence by \cite[Theorem 2.5]{Eno21}, $\cC_M(\cO)$ is an abelian length category with $M$ as its only simple object. Moreover, for $N_1,N_2\in\cC_M(\FF)$, Pontryagin duality gives
\[
\Ext^1_{\cC_M(\FF)}(N_1,N_2)
\cong\Ext^1_{D^\times,\zeta}(N_2^\vee,N_1^\vee).
\]
By \cite[IV, \S4, Th\'eor\`eme 4 and Corollaire 1]{Gab62}, there is a local pseudocompact $\cO$-algebra $(R_{\pi,\widetilde\zeta},\mathfrak{m}_R)$ with residue field $\FF$ and an $\cO$-linear equivalence between $\cC_M(\cO)$ and the category of finite-length continuous right $R_{\pi,\widetilde\zeta}$-modules, taking $M$ to $\FF$. By \cite[Lemma 3.16]{Pas13} and Theorems~\ref{thm:self-extension}, \ref{thm:reducible-self-extensions},
\[
\Hom(\mathfrak{m}_R/\overline{\mathfrak{m}_R^2},\FF)\cong\Ext^1_{R_{\pi,\widetilde\zeta}}(\FF,\FF)\cong\Ext^1_{\cC_M(\cO)}(M,M)
\]
is finite-dimensional. By the topological Nakayama lemma, the topology on $R_{\pi,\widetilde\zeta}$ is the $\mathfrak{m}_R$-adic topology.

By \cite[Lemma 3.22]{Pas13}, every deformation of $M$ to $A\in\Art_{\cO}$ belongs to $\cC_M(\cO)$. Under the above equivalence, it corresponds to an $(A,R_{\pi,\widetilde\zeta})$-bimodule which is free of rank one as a left $A$-module. As in the proof of \cite[Theorem 3.26]{Pas13}, this gives a natural bijection
\[
\operatorname{Def}_{\pi,\widetilde\zeta}(A)
\cong
\Hom_{\cO,\mathrm{cont}}(R_{\pi,\widetilde\zeta},A)/A^\times,
\]
where the homomorphisms are local and $A^\times$ acts by conjugation. By \cite[Lemma 3.31]{Pas13}, this property determines $R_{\pi,\widetilde\zeta}$ up to isomorphism. Moreover, by \cite[Corollary 3.27]{Pas13}, the maximal commutative quotient $R_{\pi,\widetilde\zeta}^{\mathrm{ab}}$
pro-represents the restriction of $\operatorname{Def}_{\pi,\widetilde\zeta}$ to commutative algebras.

\begin{proposition}\label{prop:deformation-commutativity}
    Suppose that there is a continuous local $\cO$-algebra surjection
    \[
    R_{\pi,\widetilde\zeta}\twoheadrightarrow
    \cO[[x_1,x_2,x_3]].
    \]
    Then this map is an isomorphism. In particular, $R_{\pi,\widetilde\zeta}$ is commutative.
\end{proposition}

\begin{proof}
    By \cite[Lemma 3.29]{Pas13}, Theorem~\ref{thm:self-extension} and Theorem~\ref{thm:reducible-self-extensions}, the surjection induces an isomorphism on tangent spaces. Therefore the map
    \[
    \Hom_{\cO,\mathrm{cont}}
    \left(R_{\pi,\widetilde\zeta},\FF[t]/(t^3)\right)
    \longrightarrow
    \Hom_{\cO,\mathrm{cont}}
    \left(R_{\pi,\widetilde\zeta},\FF[t]/(t^2)\right)
    \]
    is surjective. By Proposition~\ref{prop:application} and Proposition~\ref{prop:reducible-application}, there exists a $2$-dimensional $\FF$-subspace $V$ of $\Ext^1_{D^{\times},\zeta}(\pi,\pi)$ such that for any $0\ne\cE\in V$, we have 
    \[
    \dim_\FF\Ext^1_{\cC_M(\FF)}(\pi_\cE^\vee,M)
    =\dim_\FF\Ext^1_{D^\times,\zeta}(\pi,\pi_\cE)\le3.
    \]
    Here $\pi_\cE$ is an extension representing the class $\cE$. By \cite[Lemmas 3.42 and 3.43]{Pas13}, the same inequality holds for any nonzero class $\cE\in\Ext^1_{D^{\times},\zeta}(\pi,\pi)$, and the inequality is in fact an equality. The result follows by \cite[Theorem 3.39]{Pas13}.
\end{proof}

\section{Application to the mod \texorpdfstring{$p$}{p} Jacquet-Langlands correspondence}\label{sec:application}

We apply the results established in the previous sections to the representations occurring in the mod $p$ Jacquet-Langlands correspondence. 

Let $G_{\Qp}=\Gal(\Qpbar/\Qp)$ be the absolute Galois group of $\Qp$ and $I_{\Qp}$ be the inertia subgroup of $G_{\Qp}$. Let $\omega$ be the mod $p$ cyclotomic character and $\omega_2$ be Serre's fundamental character of niveau $2$. Let
\[
\rhobar:G_{\Qp}\longrightarrow\GL_2(\FF)
\]
be a continuous representation. We say that $\rhobar$ is generic if one of the following conditions holds:
\begin{enumerate}[label=\alph*)]
    \item $\rhobar$ is absolutely irreducible and up to twist
    \[
    \rhobar|_{I_{\Qp}}\cong
    \begin{pmatrix}
    \omega_2^{r+1}&0\\
    0&\omega_2^{p(r+1)}
    \end{pmatrix}
    \]
    for some $2\le r\le p-3$;
    \item $\rhobar$ is reducible and up to twist
    \[
    \rhobar|_{I_{\Qp}}\cong
    \begin{pmatrix}
    \omega^{r+1}&*\\
    0&1
    \end{pmatrix}
    \]
    for some $1\le r\le p-4$.
\end{enumerate}

Let $\pi(\rhobar)$ be the smooth admissible representation of $\GL_2(\Qp)$ attached to $\rhobar$ by the mod $p$ local Langlands correspondence. Its central character is $\zeta=\det(\rhobar)\omega^{-1}$, viewed as a character of $\Qp^\times$ via local class field theory. For $i\ge0$, let
\[
\cS^i(\pi(\rhobar)):=H^i_{\mathrm{\acute et}}\left(\PP^1_{\CC_p},\cF_{\pi(\rhobar)}\right),
\]
where $\{\cS^i\}$ are Scholze's functors constructed in \cite{Sch18}. By \cite[Theorem 1.1]{Sch18} and \cite[Lemma 7.3]{DPS23}, each $\cS^i(\pi(\rhobar))$ is an admissible smooth representation of $D^\times$ with central character $\zeta$, carrying a continuous and commuting action of $G_{\Qp}$. Now assume that $\rhobar$ is generic. By \cite[Theorem 1.4]{HW24a}, there is a $G_{\Qp}\times D^\times$-equivariant isomorphism
\begin{equation}\label{eq:JL-isotypy}
\cS^1(\pi(\rhobar))\cong\rhobar(-1)\boxtimes\mathrm{JL}(\rhobar)
\end{equation}
for some admissible smooth representation $\mathrm{JL}(\rhobar)$ of $D^\times$ with central character $\zeta$.

Next we verify our assumptions~\ref{it:assum-i}\ref{it:assum-ii} or \ref{it:assum-i-red}\ref{it:assum-ii-red} for $\mathrm{JL}(\rhobar)$.

\begin{proposition}\label{prop:JL-graded}
    Suppose that $\rhobar$ is generic. If $\rhobar$ is absolutely irreducible, then assumption~\ref{it:assum-i} holds for $\mathrm{JL}(\rhobar)$. If $\rhobar$ is reducible, then assumption~\ref{it:assum-i-red} holds for $\mathrm{JL}(\rhobar)$.
\end{proposition}

\begin{proof}
    This follows from \cite[Theorem 4.17]{HW24a} and \cite[Proposition 6.1]{HW24b}.
\end{proof}

\begin{proposition}\label{prop:JL-self-duality}
    Suppose that $\rhobar$ is generic. Then there is a $D^\times$-equivariant isomorphism of $\Lambda$-modules
    \[
    \E^2_\Lambda\left(\mathrm{JL}(\rhobar)^\vee\right)
    \otimes(\zeta^{-1}\circ\mathrm{Nr})
    \cong\mathrm{JL}(\rhobar)^\vee.
    \]
\end{proposition}

\begin{proof}
    By \cite[Theorem 1.3(ii)(a)]{HW24b}, we can assume $\rhobar=\rhobar^{\ss}.$ Write $\cS^i_{\GL_2(\Qp)}$ and $\cS^i_{D^\times}$ for the higher smooth duality functors of \cite{Koh17}. By \cite[Proposition 5.4 and Theorem 5.13]{Koh17} (see also \cite[Theorem 1.1]{Wit18} for the supersingular case), $\pi(\rhobar)$ is Cohen-Macaulay of dimension $1$, and
    \[
    \cS^1_{\GL_2(\Qp)}(\pi(\rhobar))
    \cong\pi(\rhobar)\otimes(\zeta^{-1}\circ\det).
    \]
    Since $\rhobar$ is generic, we have $\cS^2(\pi(\rhobar))=0$ by \cite[Theorem 1.2]{HW24b} when $\rhobar$ is absolutely irreducible and by \cite[Theorem 4.6]{Lud17} when $\rhobar$ is reducible. This together with \cite[Theorem 3.2 and Proposition 4.7]{Sch18} gives $\cS^i(\pi(\rhobar))=0$ for $i\ne1$.
    Using \eqref{eq:JL-isotypy}, \cite[Proposition 4.5]{HM26} and \cite[Lemma 7.4]{DPS23}, we obtain an isomorphism of $D^\times$-representations
    \[
    \left(\cS^1_{D^\times}(\mathrm{JL}(\rhobar))\right)^{\oplus2}
    \cong\left(\mathrm{JL}(\rhobar)\otimes(\zeta^{-1}\circ\mathrm{Nr})\right)^{\oplus2}.
    \]
    Both representations are admissible with central character $\zeta^{-1}$. Applying the Krull-Schmidt theory for finitely generated modules over complete semilocal noetherian rings \cite[Theorem B]{Row87} to their Pontryagin duals, we obtain
    \[
    \cS^1_{D^\times}(\mathrm{JL}(\rhobar))
    \cong\mathrm{JL}(\rhobar)\otimes(\zeta^{-1}\circ\mathrm{Nr}).
    \]
    By \cite[Corollary 3.15 and Theorems 4.1, 5.1]{Koh17}, we have a $D^\times$-equivariant isomorphism
    \[
    \left(\cS^1_{D^\times}(\mathrm{JL}(\rhobar))\right)^\vee
    \cong\E^2_\Lambda\left(\mathrm{JL}(\rhobar)^\vee\right).
    \]
    The result follows from taking Pontryagin duals. 
\end{proof}

\begin{theorem}\label{thm:JL-self-extensions}
    Suppose that $\rhobar$ is generic. Then
    \[
    \dim_\FF\Ext^i_{D^\times,\zeta}
    \left(\mathrm{JL}(\rhobar),\mathrm{JL}(\rhobar)\right)
    =\begin{cases}
    1,&i=0;\\
    3,&i=1,2;\\
    0,&i\ge3.
    \end{cases}
    \]
\end{theorem}

\begin{proof}
    This follows by Propositions~\ref{prop:JL-graded}, \ref{prop:JL-self-duality} and Theorems~\ref{thm:self-extension}, \ref{thm:reducible-self-extensions}.
\end{proof}

\begin{remark}
    \begin{enumerate}
        \item The result for $i=0$ was obtained earlier by Yongquan Hu and Vytautas Pa\v{s}k\={u}nas \cite{HP}.
        \item The same argument as in the proof of \cite[Proposition 8.1]{Pas10} gives a short exact sequence
        \[
        0\longrightarrow\Ext^1_{D^\times,\zeta}(\mathrm{JL}(\rhobar),\mathrm{JL}(\rhobar))\longrightarrow\Ext^1_{D^\times}(\mathrm{JL}(\rhobar),\mathrm{JL}(\rhobar))\longrightarrow\Hom(Z_D,\FF)\longrightarrow0.
        \]
        In particular, we have $\dim_\FF\Ext^1_{D^\times}(\mathrm{JL}(\rhobar),\mathrm{JL}(\rhobar))=5$ in the category of smooth representations without a fixed central character.
    \end{enumerate}
\end{remark}

\begin{theorem}\label{thm:commutativity}
    Use the notation in Section~\ref{sec:deformations}. Suppose that $\rhobar$ is generic and there is a continuous local $\cO$-algebra surjection
    \[
    R_{\mathrm{JL}(\rhobar),\widetilde\zeta}\twoheadrightarrow
    \cO[[x_1,x_2,x_3]].
    \]
    Then this map is an isomorphism. In particular, $R_{\mathrm{JL}(\rhobar),\widetilde\zeta}$ is commutative.
\end{theorem}

\begin{proof}
    This follows by Propositions~\ref{prop:JL-graded}, \ref{prop:JL-self-duality} and Proposition~\ref{prop:deformation-commutativity}.
\end{proof}

\bibliography{1}
\bibliographystyle{alpha}

\end{document}